\documentclass[11pt,a4paper]{article}

\usepackage[english]{babel}
\usepackage{amsmath}
\usepackage{amssymb}
\usepackage{graphicx}  

\usepackage{epstopdf}    
\usepackage{ulem}

\usepackage{amsfonts}
\usepackage{color}
\usepackage{verbatim}
\usepackage{caption}
\usepackage{subcaption}
\usepackage{fullpage}

\usepackage{array}

\usepackage{epsfig}
\usepackage{float}

\usepackage{algpseudocode}
\usepackage{url}

\usepackage{theorem}

\usepackage{fullpage}

\usepackage{hyperref}

\usepackage{tikz}

\newtheorem{rem}{Remark}[section]
\newtheorem{lem}{Lemma}[section]
\newtheorem{thm}{Theorem}[section]

\newenvironment{pf}{\noindent{\bf Proof. \/}\noindent%
}{\hfill\EndProofMarker}
\newcommand{\EndProofMarker}{$\Box$}

\def\Xint#1{\mathchoice
{\XXint\displaystyle\textstyle{#1}}%
{\XXint\textstyle\scriptstyle{#1}}%
{\XXint\scriptstyle\scriptscriptstyle{#1}}%
{\XXint\scriptscriptstyle\scriptscriptstyle{#1}}%
\!\int}
\def\XXint#1#2#3{{\setbox0=\hbox{$#1{#2#3}{\int}$ }
\vcenter{\hbox{$#2#3$ }}\kern-.6\wd0}}

\def\dashint{\Xint-}

\graphicspath{{../Imgs/}}

\begin{document}
  
\title{Strict separation and numerical approximation for a non--local Cahn--Hilliard equation with single--well potential}
\author{Abramo Agosti$^{\sharp}$\footnote{Corresponding author. E-mail address: {\tt abramo.agosti@unipv.it} \newline \textit{Email addresses:} {\tt abramo.agosti@unipv.it} (A. Agosti), {\tt elisabetta.rocca@unipv.it} (E. Rocca), {\tt luca.scarpa@polimi.it} (L. Scarpa)}, Elisabetta Rocca$^{\sharp,\S}$, Luca Scarpa$^\dag$}

\maketitle 

\begin{center}
{\small $^\sharp$  Dipartimento di Matematica \\
        Universit\`a degli studi di Pavia and IMATI-C.N.R.\\
        via Ferrata, 5 - 27100 Pavia}
\end{center}
\begin{center}
{\small $^\S$ IMATI-C.N.R.\\
        via Ferrata, 5 - 27100 Pavia}
\end{center}
\begin{center}
{\small $^\dag$  Dipartimento di Matematica \\
        Politecnico di Milano \\
        via Bonardi 9, 20133 Milano, Italy}
\end{center}

\date{}

\begin{abstract}
In this paper we study a non--local Cahn--Hilliard equation with singular single--well potential and degenerate mobility. This results as a particular case of a more general model derived for a binary, saturated, closed and incompressible mixture, composed by a tumor phase and a healthy phase, evolving in a bounded domain. The general system couples a Darcy-type evolution for the average velocity field with a convective reaction--diffusion type evolution for the nutrient concentration and a non--local convective Cahn--Hilliard equation for the tumor phase.
The main mathematical difficulties are related to the proof of the separation property for the tumor phase in the Cahn-Hilliard equation:
up to our knowledge, such problem is {indeed} open in the literature. For this reason,
in {the present} contribution we restrict the analytical study to the Cahn--Hilliard equation only.
For the non--local Cahn--Hilliard equation with singular single--well potential and degenerate mobility, {we study the existence and uniqueness of weak solutions for spatial dimensions $d\leq 3$. After showing existence, we prove the strict separation property in three spatial dimensions, implying the same property also for lower spatial dimensions, which opens the way to the proof of uniqueness of solutions}. Finally, {we propose a well posed and gradient stable continuous finite element approximation of the model {for $d\leq 3$}, which preserves the physical properties of the continuos solution and which is computationally efficient, and we show simulation results in two spatial dimensions which prove the consistency of the proposed scheme and which describe the phase ordering dynamics associated to the system.} 
\end{abstract}


\vspace*{0.5cm}


\section{Introduction}
\label{intro}

In this paper we first present a model derivation for a binary, closed, incompressible mixture of tumor cells with volume fraction $\phi_c$ and healthy cells {and liquid}, with volume fraction $\phi_l$, in presence of a nutrient species, with concentration $n$, evolving in a bounded regular domain $\Omega$ in {$\mathbb{R}^d$, where $d\leq 3$ is the spatial dimension}. 

\begin{figure}[h!]
\includegraphics[width=0.9\linewidth]
{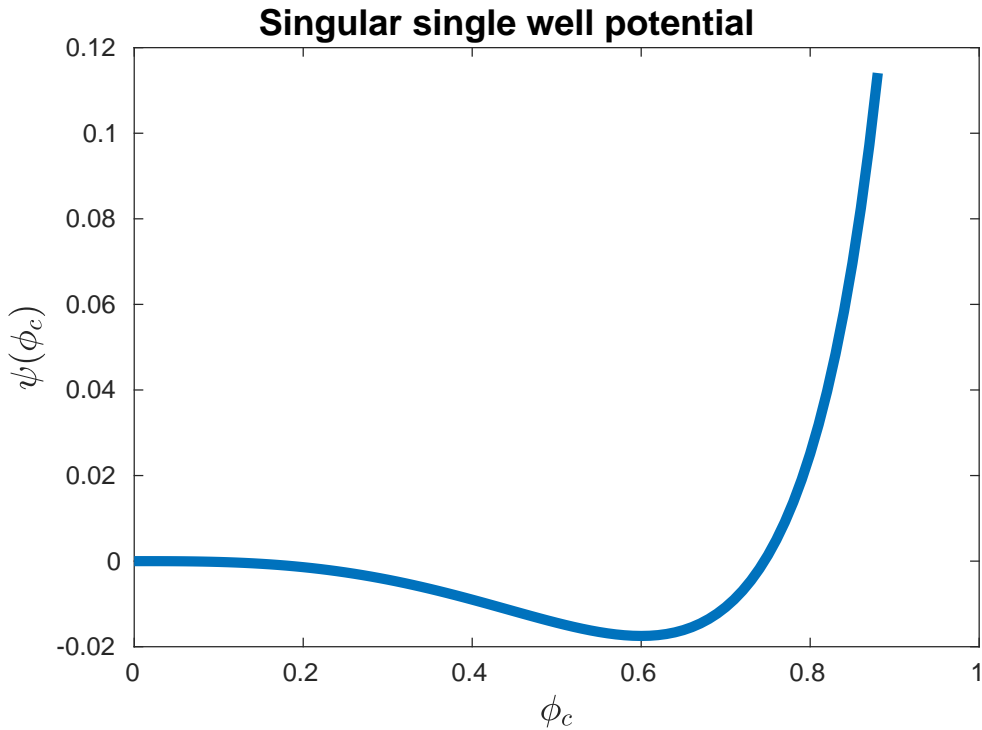}
\centering
\caption{Plot of the single well potential \eqref{eqn:4} corresponding to the value $\bar{\phi}=0.6$.}
\label{fig:1}
\end{figure}

Every component satisfies a continuity equation, where the nutrient is assumed to be advected by the mixture velocity field ${\bf v}$. 
We postulate the following form of the free energy of the system:
\begin{align}
\label{eqn:3}
&E(\phi_c,n)=\int_{\Omega}e(\phi_c,n)\,d\mathbf{x}=\\
\notag &\int_{\Omega}\Pi\biggl(\frac{\psi(\phi_c)}{\epsilon}+\frac{\chi_n}{2}|n|^2+\chi_cn(1-\phi_c)\biggr)d\mathbf{x}+\frac{\Pi\epsilon}{4}\int_{\Omega}\int_{\Omega}J(\mathbf{x},\mathbf{y})\left(\phi_c(\mathbf{x})-\phi_c(\mathbf{y})\right)^2d\mathbf{x}\,d\mathbf{y},
\end{align}
where $e(\phi_c,n)$ is the free energy per unit volume, and $\Pi$ is proportional to the surface tension, with units ({for $d=3$}) of $[N/m]$. The parameter $\epsilon$ has units of $[m]$ and is related to the interface thickness. Here, $\psi(\phi_c)$ is the bulk energy due to cell-cell mechanical interactions, which has the form of a single well potential of the \textit{Lennard--Jones} type (cf.~\cite{agosti2}):
\begin{equation}
\label{eqn:4}
\psi(\phi_c)=-(1-\bar{\phi})\log(1-\phi_c)-\frac{\phi_c^3}{3}-(1-\bar{\phi})\left(\frac{\phi_c^2}{2}+\phi_c\right),
\end{equation}
where $\bar{\phi}$ represents the equilibrium value of the cell concentration at which no interacting force is exerted between the cells. In Figure \ref{fig:1} we show a plot of the single well potential \eqref{eqn:4} corresponding to the value $\bar{\phi}=0.6$.
{We observe from Figure \ref{fig:1} that for $\phi_c<\bar{\phi}$, i.e. at a moderate cell volume fraction, cell-cell interactions are attractive, while for $\phi_c>\bar{\phi}$, i.e. at a high volume fraction, they are repulsive, with an unstable zero for $\phi_c=0$ (corresponding to the absence of cells) and an infinite cell-cell repulsion at the volume saturation level $\phi_c=1$. Moreover, the condition $\phi_c<\hat{\phi}_c$, where $\hat{\phi}_c$ is such that $\psi''(\hat{\phi}_c)=0$, defines the metastable domain of the potential \eqref{eqn:4}.
}

\noindent
The terms $\frac{\chi_n}{2}|n|^2$ and $\chi_cn(1-\phi_c)$ represent a mass term for the nutrient and an interaction term between the cells and the nutrient inducing chemotactic effects in the mixture dynamics.
A non--local term is used to represent long--range interactions between the phases, where $J:\Omega\times\Omega \to \mathbb{R}$ is a localized symmetric kernel of the form (cf., e.g., \cite{giacomin1})
\begin{equation}
\label{eqn:5}
J(\mathbf{x},\mathbf{y})=\frac{1}{\epsilon^{d+2}}J\left(\frac{|\mathbf{x}-\mathbf{y}|}{\epsilon}\right).
\end{equation}
{For biological cells dynamics and in particular for tumor growth, long range interactions are appropriate to describe cells-cells and cells-matrix adhesion, representing the sensing range over which cells can detect their surrounding \cite{armstrong}, especially during invasion and metastasis in the presence of stem cells mutations \cite{hillen}.}
The system modelling the evolution of the phases and the velocity field is then obtained in Section~\ref{model} via a generalized principle of least dissipation and it results in the following nonlinear PDE system
\begin{align}
\label{eqn:11intro}
&\mathbf{v}=-k\left(\nabla \bar{p} - \left(\mu+\chi_c n\right) \nabla \phi_c\right),\\
&\text{div}\mathbf{v}=0,\\
&\frac{\partial \phi_c}{\partial t}+\mathbf{v}\cdot \nabla \phi_c -
\text{div}\biggl(b(\phi_c)\mathbf{T}\nabla \mu\biggr)=\frac{\Gamma_c}{\gamma} ,\\
&\mu:=\frac{\delta E}{\delta \phi_c}=\frac{\psi'(\phi_c)}{\epsilon}-\chi_c n+\epsilon \left(J\star 1\right)\phi_c-\epsilon J\star \phi_c,\\
&\frac{\partial n}{\partial t}+\mathbf{v}\cdot \nabla n - \text{div}\left(D(\phi_c)\mathbf{D}\nabla (\chi_n n+\chi_c(1-\phi_c))\right)= S,
\end{align} 
where $k$ is a friction parameter, $b(\phi_c), D(\phi_c)$ are positive mobilities, assumed to be dependent only on $\phi_c$, $\mathbf{T}, \mathbf{D}$ are symmetric positive definite tensors, representing anisotropy in the nutrient diffusion and in the cells motility, respectively. {Here, we indicate the convolution of the kernel $J$ with a function $f:\Omega\to \mathbb{R}$ as $J\star f:\Omega\to \mathbb{R}$, where
\[
(J\star f)(\mathbf{x}):=\int_{\Omega}J(\mathbf{x}-\mathbf{y})f(\mathbf{y})d\mathbf{y}.
\]
}
The source terms $\Gamma_c, \Gamma_l$ represent cells proliferation and death and conversion of mass between the two phases, while $S$ includes source and consumption terms for the nutrient. In order to study the well-posedness of this PDE system coupled with suitable initial and boundary conditions, we start with the analysis of a reduced case, where ${\bf T}\equiv{\bf D}\equiv {\bf Id}$, ${\bf v}\equiv \boldsymbol{0}$, $\Gamma_c\equiv n\equiv 0$, that is the following non--local Cahn--Hilliard equation 
\begin{equation}
\label{eqn:22intro}
\begin{cases}
\frac{\partial \phi_c}{\partial t}-\text{div}\biggl(b(\phi_c)\nabla \mu\biggr)=0 ,\\
\mu=\frac{\psi'(\phi_c)}{\epsilon}+\epsilon \left(J\star 1\right)\phi_c-\epsilon J\star \phi_c,
\end{cases}
\end{equation} 
endowed with the homogeneous boundary condition
\begin{equation}
\label{eqn:23intro}
b(\phi_c)\nabla \mu \cdot \boldsymbol{\nu}|_{\partial \Omega}=0.
\end{equation}
For modelling reasons (cf. Section~\ref{model}) we consider here  a general expression for the mobility 
\begin{equation}
\label{eqn:22aintro}
b(\phi_c)=\frac{\phi_c^\alpha(1-\phi_c)^2}{M},\quad\alpha\geq 0,
\end{equation} 
where $M$ denotes a generic friction constant. 
{This particular
form of the mobility will be derived in Section~\ref{model}, where it will be shown that different values of $\alpha$ correspond to different empirical laws for the filtration processes between the phases in the mixture. We observe that the case with $\alpha=1$ was considered in \cite{chat} in the framework of the  application of mixture model for solid tumors. As we will show, this case corresponds to a filtration process described by a Kozeny--Karman permeability law.}

Equation \eqref{eqn:22intro} models, indeed, the evolution of biological cells such as solid tumors. The degeneracy set of the the mobility and set of singularities of the cellular potential do not coincide, and
the absence of cells is an unstable potential equilibrium configuration.
This feature introduces a non-trivial difference with respect to the standard Cahn--Hilliard equation analyzed in the literature.  
Regarding the case of standard non--local Cahn--Hilliard equation with double--well type potentials, we can refer to \cite{gz}, where a model of
Cahn--Hilliard type for phase separation in a two--phase system
involving non--local interactions has been presented and to \cite{ckrs1} where an abstract equation including the Cahn--Hilliard type  one, possibly also with non--smooth potential, has been studied. 

The Cahn--Hilliard model itself goes back
to \cite{cahnhill} and a fairly complete review on the recent
related literature can be found e.\,g.~in \cite{mi}. 
In these papers, the nonlinearity $\psi$ is always represented by a double--well smooth potential 
{or of logarithmic type}. However, in modeling tumor
growth, this form of the potential seems unphysical for biological cells, since it has been observed that cell--cell interactions are
attractive at a moderate cell volume fraction  and repulsive at a high volume
fraction, with a zero for $\phi_c=0$ and an infinite cell--cell repulsion as $\phi_c$ approaches the value 1, (c.f., \cite{bp}).
Hence, in this paper, we propose to use a single--well potential of Lennard--Jones type (cf.~\eqref{eqn:4}). As reference for this case, we can quote the paper \cite{agosti2}, where existence results for different classes of weak solutions have been given for the case of the local Cahn--Hilliard type equation and  a continuous finite element approximation of the problem, where the positivity of the solution is enforced through a discrete variational inequality, has been provided.

Here we are interested in the non--local Cahn--Hilliard type equation \eqref{eqn:22intro} with singular single--well type potential \eqref{eqn:4}, degenerate mobility \eqref{eqn:23intro} and convolution kernel \eqref{eqn:5}, because for this type of equation we can prove the so--called ``strict--separation property'', which guarantees that the phase parameter $\phi_c$ stays confined in the interval $(\delta,1-\delta)$ for some positive constant $\delta$ in case $\alpha\geq 2$. Up to our knowledge, this is the first result in the framework of the non--local Cahn--Hilliard equation with single--well potential. Notice that this result opens the way to other investigations, like the ones related to regularity of solutions or to the associated optimal control analysis, due to the fact that the singularity of the potential does not apply anymore in the interval $[\delta, 1-\delta]$. 

Regarding instead results on the separation property for the case of non--local Cahn--Hilliard equation with double--well potential, we can quote the 
{
references \cite{ggg,GGG2, frigeri1, frigeri2, frigeri3}, as well as the recent 
contributions \cite{G_sep, P_sep} in dimension three.}

{We observe that the strict separation property implies that pure phases of the mixture are never attained if the system starts its evolution from an initial condition with no pure phases. In the context of the dynamics of a binary mixture of biological cells (such as a mixture of cancer and healthy cells) with no source or sink terms this means that a cells species cannot saturate all the available volume at the expense of the other species, due to mechanical repulsion and friction forces between the cells.

We are also interested in the design of a continuous finite element approximation of the model which is well posed and gradient stable and which preserves the physical properties of the continuous solutions.
}

The plan of the paper is the following. In the next Section~\ref{notation} we introduce the notation and the functional setting. In Section~\ref{model} we introduce a general model for a binary, closed, incompressible mixture of tumor cells in presence of a nutrient species, which is assumed to be advected by the mixture velocity field. In Section~\ref{exi} we prove existence of weak solutions for  the single non--local Cahn--Hilliard equation \eqref{eqn:22intro} in case of singular single--well potential and degenerate mobility, while, in Section~\ref{separation}, we prove the main result of the paper, which is the separation property (and consequently uniqueness of solutions) based on a Moser iteration scheme. {The latter property is derived in three spatial dimensions, which implies the same result also for lower spatial dimensions.} {In Section~\ref{numerics}, we propose a well posed and gradient stable continuous finite element approximation of the model {for $d\leq 3$} which preserves the physical properties of the continuos solution and which is computationally efficient. Finally, in Section~\ref{Simulations} we report the results of numerical simulations for different test cases in two spatial dimensions, which prove the consistency of the proposed numerical scheme and which describe the dynamics of the spinodal decomposition for the considered model.}

\section{Notation and functional setting}
\label{notation}
Let {$\Omega\subseteq \mathbb{R}^d$, $d\leq 3$}, be an open bounded domain. Let $T > 0$ denote
some final time, and set $\Omega_T:=\Omega\times (0,T)$. We indicate with $L^p(\Omega)$, $W^{m,p}(\Omega)$, $H^m(\Omega)=W^{m,2}(\Omega)$ and $L^p((0,T);X)$ the usual Lebesgue, Sobolev and Bochner spaces, for $p\in [1,\infty]$ and $m \in \mathbb{N}$. We moreover set $H:=L^2(\Omega)$, $V:=H^1(\Omega)$.  For a normed space $X$, the associated norm is denoted by $||\cdot||_X$.  When $X=L^2(\Omega)$, we denote by $(\cdot,\cdot)$ and $||\cdot||$ the standard $L^2$ inner product and induced norm respectively. If $X$ is a Banach space, we denote as $X'$ its dual space. When $X=V$, we denote by $<\cdot,\cdot>$ the duality pairing between $V'$ and $V$. We also indicate with $C(\bar{\Omega})$ the space of continuous functions from $\bar{\Omega}$ to $\mathbb{R}$.

Furthermore, $C$ denotes throughout a generic positive constant independent of
the unknown variables, the discretization and the regularization parameters, the value of which might change from line to line; $C_1,C_2,\dots$ indicate generic positive constants whose particular value must be tracked through the calculations; $C(a)$ denotes a constant depending on the non-negative parameter $a$, such that, for $C_1>0$, if $a\leq C_1$, there exists a $C_2>0$ such that $C(a)\leq C_2$.

It is useful to introduce the \textit{inverse Laplacian} operator $\mathcal{G}:V_0'\rightarrow V_0$  such that
\begin{equation}
\label{eqn:greencont}
(\nabla \mathcal{G}v,\nabla \eta)=<v,\eta > \quad \forall \eta \in V,
\end{equation}
where $V_0':=\{v\in V':<v,1>=0\}$ and $V_0=\{v\in V:(v,1)=0\}$. The existence and uniqueness of an element $\mathcal{G}v\in V_0$, for any $v\in V_0'$, follows from the Lax-Milgram Theorem and the Poincar\'e's inequality.
\\
We can define a norm on $V_0'$ by setting
\begin{equation}
\label{eqn:normf}
||v||_{-1}:=||\nabla \mathcal{G}v||\equiv <v,\mathcal{G}v>^{1/2} \quad \forall v\in V_0'.
\end{equation}
We recall the Gagliardo-Nirenberg inequality {(see e.g. \cite{brezis,gagliardo,nirenberg})}.
\begin{lem}
\label{lem:gagliardoniremberg}
Let {$\Omega \subset \mathbb{R}^d$, $d\leq 3$}, be a bounded domain with Lipschitz boundary and $f\in W^{m,r}\cap L^q$, $q\geq 1$, $r\leq \infty$, where $f$ can be a function with scalar, vectorial or tensorial values. For any integer $j$ with $0 \leq j < m$, suppose there is $\alpha \in \mathbb{R}$ such that
\[
j-\frac{d}{p}=\left(m-\frac{d}{r}\right)\alpha+(1-\alpha)\left(-\frac{d}{q}\right), \quad \frac{j}{m}\leq \alpha \leq 1.
\]
Then, there exists a positive constant $C$ depending on $\Omega$, d, m, j, q, r and $\alpha$ such that
\begin{equation}
\label{eqn:gagliardoniremberg}
||D^jf||_{L^p(\Omega)}\leq C||f||_{W^{m,r}(\Omega)}^{\alpha}||f||_{L^q(\Omega)}^{1-\alpha}.
\end{equation}
\end{lem}

Let $\mathcal{T}_{h}$ be a quasi-uniform conforming decomposition of $\Omega$ into $d-$simplices $K$, 
and let us introduce the following finite element spaces:
\begin{align}
\notag & S^{h} := \{\chi \in C(\bar{\Omega}):\chi |_{K}\in P^{1}(K) \; \forall K\in \mathcal{T}_{h}\}\subset V,\\
\notag & K^{h} := \{\chi \in S^{h}: \chi \geq 0\; \rm in \, \Omega\}
\end{align}
where $\mathbb{P}_{1}(K)$ indicates the space of polynomials of total order one on $K$.

Let $\mathcal{I}$ be the set of nodes of $\mathcal{T}_{h}$ and $\{\mathbf{x}_j\}_{j\in \mathcal{I}}$ be the set of their coordinates. Moreover, let $\{\chi_j\}_{j\in \mathcal{I}}$ be the Lagrangian basis functions associated with each node $j\in \mathcal{I}$.
Denoting by $\pi^h:C(\bar{\Omega})\rightarrow S^h$ the standard Lagrangian interpolation operator we define the lumped scalar product as
\begin{equation}
\label{eqn:lump}
(\eta_1,\eta_2)^h=\int_{\Omega}\pi^h(\eta_1(\mathbf{x})\eta_2(\mathbf{x}))d\mathbf{x}\equiv \sum_{j\in \mathcal{I}}(1,\chi_j)\eta_1(\mathbf{x}_j)\eta_2(\mathbf{x}_j),
\end{equation}
for all $\eta_1,\eta_2\in C(\bar{\Omega})$. We observe that, since $\chi_j\in K^h$ for all $j\in \mathcal{I}$, the lumped scalar product induces a norm on $S^h$ (while it induces a seminorm in higher order finite element spaces).
Given $\eta\in C(\bar{\Omega})$, $J(\cdot)\in C(\bar{\Omega})$, we define the approximation of the convolution function $J\star \eta$ as
\begin{equation}
    \label{eqn:lump2}
    (J\star \eta)_h:=\pi_h\left(\sum_{j\in \mathcal{I}}J(\mathbf{x}-\mathbf{x}_j)\eta(\mathbf{x}_j)(1,\chi_j)\right)=\sum_{i,j\in \mathcal{I}}J(\mathbf{x}_i-\mathbf{x}_j)\eta(\mathbf{x}_j)(1,\chi_j)\chi_i(\mathbf{x}).
\end{equation}
\begin{rem}
    \label{rem:astarpos}
    If $J(\cdot)\geq 0$ (resp. $>0$) we have that $(J\star 1)_h\geq 0$ (resp. $>0$). This is a consequence of the fact that, since we are working with first order finite elements method, $\chi_i\in K^h$.
    We observe that for higher order finite elements this property could be not satisfied. 
\end{rem}
We also introduce the lumped scalar product between a convolution function $J\star \eta_1$ and $\eta_2$ as
\begin{align}
\label{eqn:lump3}
&(J\star \eta_1,\eta_2)^{h^2}:=((J\star \eta_1)_h,\eta_2)^{h}\equiv \int_{\Omega}\int_{\Omega}\pi_x^h\{\pi_y^h[J(\mathbf{x}-\mathbf{y})\eta_1(\mathbf{y})]\eta_2(\mathbf{x})\}d\mathbf{y}d\mathbf{x}\\
& \displaystyle \notag =\sum_{i,j\in \mathcal{I}}(1,\chi_i)(1,\chi_j)J(\mathbf{x}_i-\mathbf{x}_j)\eta_1(\mathbf{x}_j)\eta_2(\mathbf{x}_i).
\end{align}
We observe that the formula \eqref{eqn:lump}-\eqref{eqn:lump3} correspond to considering an approximation of integral terms through a quadrature formula based on the nodes in $\mathcal{I}$ and on the trapezoidal integration rule (see e.g. \cite[Section 11.4]{quart}). We give the following result, which will be useful in the sequel.
\begin{lem}
    \label{lem:jhcomp}
    Let $\eta\in C(\bar{\Omega})$ and $J(\cdot)\in C(\bar{\Omega})$, with $J(\mathbf{x})=J(-\mathbf{x})$. Then, the following relation is satisfied:
    \begin{equation}
        \label{eqn:jhcomp}
      \left((J\star 1)_h \eta,\eta\right)^{h}\pm(J\star \eta,\eta)^{h^2}=\sum_{i,j\in \mathcal{I}}J(\mathbf{x}_i-\mathbf{x}_j)\left(\eta(\mathbf{x}_i)\pm\eta(\mathbf{x}_j)\right)^2(1,\chi_i)(1,\chi_j).  
    \end{equation}
\end{lem}
\begin{pf}
    Using \eqref{eqn:lump2}, \eqref{eqn:lump3} and the symmetry of the kernel $J$, we may write
    \begin{align*}
     &\left((J\star 1)_h \eta,\eta\right)^{h}\pm(J\star \eta,\eta)^{h^2} \\
     &=\sum_{i,j\in \mathcal{I}}J(\mathbf{x}_i-\mathbf{x}_j)\eta^2(\mathbf{x}_i)(1,\chi_i)(1,\chi_j)\pm \sum_{i,j\in \mathcal{I}}J(\mathbf{x}_i-\mathbf{x}_j)\eta(\mathbf{x}_i)\eta(\mathbf{x}_j)(1,\chi_i)(1,\chi_j)\\
     &=\frac{1}{2}\sum_{i,j\in \mathcal{I}}J(\mathbf{x}_i-\mathbf{x}_j)\eta^2(\mathbf{x}_i)(1,\chi_i)(1,\chi_j)+\frac{1}{2}\sum_{i,j\in \mathcal{I}}J(\mathbf{x}_i-\mathbf{x}_j)\eta^2(\mathbf{x}_j)(1,\chi_i)(1,\chi_j)\\
     & \pm\sum_{i,j\in \mathcal{I}}J(\mathbf{x}_i-\mathbf{x}_j)\eta(\mathbf{x}_i)\eta(\mathbf{x}_j)(1,\chi_i)(1,\chi_j),
    \end{align*}
    from which we get \eqref{eqn:jhcomp}.
\end{pf}

\noindent
We introduce the $L^2$ lumped projection operator 
$\hat{P}^h:L^2(\Omega)\rightarrow S^h$ defined by
\begin{align}
\label{eqn:l2proj}
& (\hat{P}^h\eta,\chi)^h=(\eta,\chi) \quad \forall \chi \in S^h.
\end{align}
We recall the following well-known results, (see e.g. \cite{barrett}).
\begin{lem}
The following properties hold
\begin{align}
\label{eqn:interp1}
&||\chi||_{L^{p_2}(\Omega)}  \leq  Ch^{-d(\frac{1}{p_1}-\frac{1}{p_2})}||\chi||_{L^{p_1}(\Omega)} \quad  \forall \chi \in S^h, \; 1\leq p_1 \leq p_2 \leq \infty; \\
\label{eqn:interp5}
&||\chi||^{2} \leq (\chi,\chi)^{h} \leq (d+2)||\chi||^{2} \quad \; \; \forall \chi \in S^h.
\end{align}
\end{lem}
Similarly to \eqref{eqn:greencont}, we define the discrete Green operator 
$\hat{\mathcal{G}}^h:\mathcal{F}^h\rightarrow V^h$ as follows
\begin{align}
\label{eqn:greendiscr}
(\nabla \hat{\mathcal{G}}^h v,\nabla \chi)&=(v,\chi )^h \quad \forall \chi \in S^h,
\end{align}
where $\mathcal{F}^h=\{v\in \bar{C}(\Omega) : (v,1)=0\}$ and $V^h=\{v^h\in S^h:(v^h,1)=0\}$.

\section{Derivation of the model}
\label{model}

We consider a binary, saturated, closed and incompressible mixture in an open bounded domain {$\Omega\subset \mathbb{R}^d$}, composed by a tumor phase $c$ with volume fraction $\phi_c$ and a phase $l$ composed by liquid and healthy cells  with volume fraction $\phi_l$. In the sequel, we will call $c$ the cells phase and $l$ the liquid phase. Both phases have a constant density $\gamma$, equal to the water density. We assume that the mixture dynamics is coupled with the evolution of a nutrient species, with concentration $n$. Every component satisfies a continuity equation,
\begin{align}
\label{eqn:1a}
\frac{\partial \phi_c}{\partial t}+\text{div}(\phi_c \mathbf{v})+\text{div}(\mathbf{J_c})=\frac{\Gamma_c(n,\phi_c)}{\gamma},\\
\label{eqn:1b}
\frac{\partial \phi_l}{\partial t}+\text{div}(\phi_l \mathbf{v})+\text{div}(\mathbf{J_l})=\frac{\Gamma_l(n,\phi_c)}{\gamma},\\
\label{eqn:1c}
\frac{\partial n}{\partial t}+\text{div}(n \mathbf{v})+\text{div}(\mathbf{J_n})=S(n,\phi_c),
\end{align}
with $\phi_c+\phi_l=1$, $\Gamma_c+\Gamma_l=0$, $\mathbf{J_c}=\phi_c(\mathbf{v_c}-\mathbf{v})$, $\mathbf{J_l}=\phi_l(\mathbf{v_l}-\mathbf{v})$. Here, $\mathbf{v}=\phi_c \mathbf{v_c}+\phi_l \mathbf{v_l}$ is the volume-averaged mixture velocity, which satisfies the incompressibility condition 
\begin{equation}
\label{eqn:2}
\text{div}\mathbf{v}=0,
\end{equation}
as a consequence of the saturation and the closedness properties of the mixture. 
The nutrient is assumed to be advected by the mixture velocity $\mathbf{v}$ and transported by the flux $\mathbf{J_n}$. 

We postulate the form \eqref{eqn:3} (with $\Pi=1$ for ease of notation) of the free energy of the system with a single well potential of the \textit{Lennard--Jones} type \eqref{eqn:4} and a non--local interactions between the phases, given by a kernel $J$ of the form \eqref{eqn:5}. 

We note that the local diffuse--interface free energy introduced in \cite{agosti1} can be obtained from \eqref{eqn:3}, \eqref{eqn:4} and \eqref{eqn:5} by considering the first order terms of a formal expansion in $\epsilon$ of the non--local kernel (see e.g. \cite{frigeri1} for details).
 \begin{rem}
 \label{rem:1}
 We can rewrite the free energy \eqref{eqn:3} as
 \begin{align}
\label{eqn:6}
&\int_{\Omega}\biggl(\frac{\psi(\phi_c)}{\epsilon}-\epsilon \left(\frac{J\star 1}{2}\right)\phi_c(1-\phi_c)+\frac{\chi_n}{2}|n|^2+\chi_cn(1-\phi_c)\biggr)d\mathbf{x}\\
& \notag +\frac{\epsilon}{2}\int_{\Omega}\int_{\Omega}J(\mathbf{x}-\mathbf{y})\phi_c(\mathbf{x})(1-\phi_c(\mathbf{y}))d\mathbf{x}\,d\mathbf{y},
\end{align}
and observe that the non--local term in \eqref{eqn:6} represents non--local adhesive interactions between the cells and the liquid phases. The term $\frac{\psi(\phi_c)}{\epsilon}-\epsilon \left(\frac{J\star 1}{2}\right)\phi_c(1-\phi_c)$ represents the phenomenological analogue of the entropy of mixing between the phases (see \cite{giacomin1}), which is strictly convex, for all $\mathbf{x}\in \Omega$, when
\begin{equation}
\label{eqn:7}
\epsilon \inf_{\mathbf{x}\in \Omega}\left(J\star 1\right)>\frac{2+(1-\bar{\phi})-3\sqrt[3]{1-\bar{\phi}}}{\epsilon}.
\end{equation}
Due to the localized form \eqref{eqn:5} of the kernel $J$, the condition \eqref{eqn:7} is realistic even for small values of $\epsilon$. Indeed, it's easy to show that $\epsilon \left(J\star 1\right)\sim \frac{1}{\epsilon}$.
 \end{rem}
We now derive the model for the dynamics of the system satisfying the second law of thermodynamics in isothermal situation and a generalized variational principle of least dissipation.\\
Let us firstly introduce the functional derivatives of the free energy $E(\phi_c,n)$:

\[\mu:=\frac{\delta E}{\delta \phi_c}=\frac{\psi'(\phi_c)}{\epsilon}-\chi_{c}n+\epsilon \left(J\star 1\right)\phi_c-\epsilon J\star \phi_c,\] 
\[\eta:=\frac{\delta E}{\delta n}=\chi_nn+\chi_c(1-\phi_c).\]
We multiply equation \eqref{eqn:1a} by $\mu$, equation \eqref{eqn:1c} by $\eta$, \eqref{eqn:2} by a multiplier $p$, integrate the resulting terms over a generic volume $R(t)\subset \Omega$ transported with the mixture velocity and sum the three contributions. We get
\begin{align*}
& \int_{R(t)}\mu \left(\frac{\partial \phi_c}{\partial t}+\mathbf{v}\cdot \nabla \phi_c + \phi_c \text{div}\mathbf{v}+\text{div}\mathbf{J}_c\right)\,d\mathbf{x}+ \int_{R(t)}\eta \left(\frac{\partial n}{\partial t}+\mathbf{v}\cdot \nabla n + n \text{div}\mathbf{v}+\text{div}\mathbf{J}_n\right)\,d\mathbf{x}\\
&+\int_{R(t)}p\text{div}\mathbf{v}\,d\mathbf{x}=\int_{R(t)}\left(\frac{\Gamma_c}{\gamma}\mu+S\eta\right)\,d\mathbf{x}.
\end{align*}
Inserting the explicit expressions of $\mu$ and $\eta$, rearranging the terms and integrating by parts we get
\begin{align*}
& \int_{R(t)}\left(\left(\frac{\psi'(\phi_c)}{\epsilon}-\chi_{c}n+\epsilon \left(J\star 1\right)\phi_c-\epsilon J\star \phi_c\right)\frac{\partial \phi_c}{\partial t}+\left(\chi_nn+\chi_c(1-\phi_c)\right)\frac{\partial n}{\partial t}\right)\,d\mathbf{x}\\
&+\int_{R(t)}\left(\left(\frac{\psi'(\phi_c)}{\epsilon}-\chi_{c}n+\epsilon \left(J\star 1\right)\phi_c- \epsilon J\star \phi_c\right)\mathbf{v}\cdot \nabla \phi_c+\left(\chi_nn+\chi_c(1-\phi_c)\right)\mathbf{v}\cdot \nabla n\right)\,d\mathbf{x}\\
&+\int_{R(t)}\left(\phi_c\mu+\chi_nn^2+\chi_cn(1-\phi_c)\right)\text{div}\mathbf{v}\,d\mathbf{x}+\int_{S(t)}\left(\mu \mathbf{J}_c+\eta\mathbf{J}_n\right)\cdot\boldsymbol{\nu}\,dS-\\
& \int_{R(t)}\left(\mathbf{J}_c\cdot \nabla \mu+\mathbf{J}_n\cdot \nabla \eta \right)\,d\mathbf{x}+\int_{R(t)}p\text{div}\mathbf{v}\,d\mathbf{x}=\int_{R(t)}\biggl(\frac{\partial e}{\partial t}+\mathbf{v}\cdot \nabla e + \\
&\left(\frac{\chi_n}{2}n^2+\chi_cn(1-\phi_c)\right)\text{div}\mathbf{v}+\left(\mu \phi_c+\frac{\chi_n}{2}n^2+p\right)\text{div}\mathbf{v}\biggr)\,d\mathbf{x}+\\
&\int_{R(t)}\left(- \frac{\epsilon}{2} \left(J\star \phi_c\right)\mathbf{v}\cdot \nabla \phi_c- \frac{\epsilon}{2} \mathbf{v}\cdot \left(\nabla J\star 1\right)\phi_c^2+ \frac{\epsilon}{2} \mathbf{v}\cdot \left(\nabla J\star \phi_c\right)\phi_c\right)\,d\mathbf{x}+\\
&\int_{S(t)}\left(\mu \mathbf{J}_c+\eta\mathbf{J}_n\right)\cdot\boldsymbol{\nu}\,dS-\int_{R(t)}\left(\mathbf{J}_c\cdot \nabla \mu+\mathbf{J}_n\cdot \nabla \eta \right)\,d\mathbf{x}=\int_{R(t)}\left(\frac{\Gamma_c}{\gamma}\mu+S\eta\right)\,d\mathbf{x}.
\end{align*}
We rewrite the terms
\begin{align*}
& \int_{R(t)}\left(- \frac{\epsilon}{2} \left(J\star \phi_c\right)\mathbf{v}\cdot \nabla \phi_c- \frac{\epsilon}{2} \mathbf{v}\cdot \left(\nabla J\star 1\right)\phi_c^2+ \frac{\epsilon}{2} \mathbf{v}\cdot \left(\nabla J\star \phi_c\right)\phi_c\right)\,d\mathbf{x}\\
&=\int_{R(t)}\left(\mathbf{v}\cdot\left(\mu+\chi_cn\right)\nabla \phi_c+\left(-\frac{\epsilon}{2}\left(J\star \phi_c\right)\phi_c+\frac{\epsilon}{2}\left(J\star 1\right)\phi_c^2+\frac{\psi(\phi_c)}{\epsilon}\right)\text{div}\mathbf{v}\right)\,d\mathbf{x}\\
&+\int_{S(t)}\left(\frac{\epsilon}{2}\left(J\star \phi_c\right)\phi_c-\frac{\epsilon}{2}\left(J\star 1\right)\phi_c^2-\frac{\psi(\phi_c)}{\epsilon}\right)\mathbf{v}\cdot \boldsymbol{\nu}\,dS,
\end{align*}
and obtain
\begin{align}
\label{eqn:8}
&\int_{R(t)}\left(\frac{\partial e}{\partial t}+\mathbf{v}\cdot \nabla e + e\text{div}\mathbf{v}+\mathbf{v}\cdot\left(\nabla\left(-p-\frac{\chi_n}{2}n^2-\mu\phi_c\right)+\left(\mu+\chi_cn\right)\nabla \phi_c\right)\right)\,d\mathbf{x}\\
&\notag+\int_{S(t)}\left(\mu \mathbf{J}_c+\eta\mathbf{J}_n+\left(p+\mu\phi_c+\eta n -e\right)\mathbf{v}\right)\cdot\boldsymbol{\nu}\,dS-\int_{R(t)}\left(\mathbf{J}_c\cdot \nabla \mu+\mathbf{J}_n\cdot \nabla \eta \right)\,d\mathbf{x}\\
&\notag=\int_{R(t)}\left(\frac{\Gamma_c}{\gamma}\mu+S\eta\right)\,d\mathbf{x}.
\end{align}

The second law of thermodynamics in isothermal situations and with source terms takes the form of the following dissipation inequality \cite{garcke1,gurtin1}
\begin{equation}
\label{eqn:9}
\frac{d}{dt}\int_{R(t)}e\leq -\int_{\partial R(t)}\mathbf{J_E}\cdot \boldsymbol{\nu}dS+\int_{R(t)}\biggl(k_c\frac{\Gamma_c}{\gamma}+k_nS\biggr)d\mathbf{x},
\end{equation}
for each material volume $R(t)\subset \Omega$, with the energy flux $\mathbf{J_E}$ and the multipliers $k_c, k_n$ to be determined.
By comparing \eqref{eqn:8} with \eqref{eqn:9}, the 
 constitutive assumptions can be made in order for the system to fulfill the dissipation inequality \eqref{eqn:9}
\begin{align}
\label{eqn:10a}
&\bar{p}=p-\frac{\chi_n}{2}n^2-\mu \phi_c,\\
\label{eqn:10b}
&\mathbf{v}=-k\left(\nabla \bar{p}-\left(\mu+\chi_c n\right)\nabla \phi_c\right),\\
\label{eqn:10c}
&\mathbf{J_c}=-b(\phi_c)\mathbf{T}\nabla \mu,\\
\label{eqn:10d}
&\mathbf{J_n}=-D(\phi_c)\mathbf{D}\nabla (\chi_nn+\chi_c(1-\phi_c)),\\
\label{eqn:10e}
&\mathbf{J_E}=\mu \mathbf{J}_c+\eta\mathbf{J}_n+\left(p+\mu\phi_c+\eta n -e\right)\mathbf{v},\\
&k_c=\mu, \quad k_n=\chi_nn+\chi_c(1-\phi_c),
\end{align}
where $k$ is a friction parameter, $b(\phi_c), D(\phi_c)$ are positive mobilities, assumed to be dependent only on $\phi_c$, $\mathbf{T}, \mathbf{D}$ are symmetric positive definite tensors, representing anisotropy in the nutrient diffusion and in the cells motility respectively.

Inserting \eqref{eqn:10a}--\eqref{eqn:10d} in \eqref{eqn:1a}--\eqref{eqn:1c}, we get the following system of equations
\begin{equation}
\label{eqn:11}
\begin{cases}
\mathbf{v}=-k\left(\nabla \bar{p} - \left(\mu+\chi_c n\right) \nabla \phi_c\right),\\
\text{div}\mathbf{v}=0,\\
\frac{\partial \phi_c}{\partial t}+\mathbf{v}\cdot \nabla \phi_c + \text{div}\biggl(b(\phi_c)\mathbf{T}\nabla \mu\biggr)=\frac{\Gamma_c}{\gamma} ,\\
\mu=\frac{\psi'(\phi_c)}{\epsilon}-\chi_c n+\epsilon \left(J\star 1\right)\phi_c-\epsilon J\star \phi_c,\\
\frac{\partial n}{\partial t}+\mathbf{v}\cdot \nabla n - \text{div}\left(D(\phi_c)\mathbf{D}\nabla (\chi_n n+\chi_c(1-\phi_c))\right)= S,
\end{cases}
\end{equation} 
in $\Omega$, which we endow with the homogeneous boundary conditions
\begin{equation}
\label{eqn:11bis}
b(\phi_c)\mathbf{T}\nabla \mu \cdot \boldsymbol{\nu}|_{\partial \Omega}=D(\phi_c)\mathbf{D}\nabla (\chi_n n+\chi_c(1-\phi_c)) \cdot \boldsymbol{\nu}|_{\partial \Omega}=\mathbf{v}|_{\partial \Omega} =0 \to \mathbf{J_E}\cdot \boldsymbol{\nu}|_{\partial \Omega}=0.
\end{equation}
A solution of system \eqref{eqn:11}, supplemented with the boundary conditions \eqref{eqn:11bis}, formally satisfies the following energy equality
\begin{align}
\label{eqn:12}
&\frac{dE}{dt}+\frac{1}{k}\int_{\Omega}\mathbf|v|^2 d\mathbf{x} +\int_{\Omega}b(\phi_c)\mathbf{T}\nabla \mu\cdot \nabla \mu d\mathbf{x} \\
& \notag +\int_{\Omega}D(\phi_c) \mathbf{D}\nabla (\chi_nn+\chi_c(1-\phi_c))\cdot \nabla (\chi_nn+\chi_c(1-\phi_c)) d\mathbf{x}\\
& \notag =\int_{\Omega}\frac{\Gamma_c}{\gamma}\mu d\mathbf{x}+\int_{\Omega}S(\chi_nn+\chi_c(1-\phi_c)) d\mathbf{x},
\end{align}
which is obtained by multiplying the first equation of \eqref{eqn:11} by $\frac{1}{k}\mathbf{v}$, the second equation by $\bar{p}$, the third equation by $\mu$, the fourth equation by $\frac{\partial \phi_c}{\partial t}$, the fifth equation by $\chi_n n+\chi_c(1-\phi_c)$, integrating over $\Omega$ and summing all the contributions, observing also that
\[
-\chi_c\int_{\Omega}\left(1-\phi_c\right)\mathbf{v}\cdot \nabla n\,d\mathbf{x}=-\chi_c\int_{\partial \Omega}\left(n(1-\phi_c)\right)\mathbf{v}\cdot \boldsymbol{\nu}\,dS-\chi_c\int_{\Omega}n\mathbf{v}\cdot \nabla \phi_c.
\]
In order to close the system \eqref{eqn:11} and determine particular forms for the mobility functions $b(\phi_c)$ and $D(\phi_c)$, we apply the Onsager Variational Principle (OVP) \cite{onsager1}, which defines the irreversible non--equilibrium dynamics for near--equilibrium thermodynamic systems in terms of linear fluxes--forces balance equations. This principle has been widely applied for the derivation of continuum phenomenological models of soft matter, see e. g. \cite{onsager2,onsager3,onsager4}, based on the observation that the macroscopic behavior of soft matter is driven by the interactions between its {constituents} units, at the constituent scale, which drive the system locally in a metastable regime out of equilibrium, and by the dissipation mechanisms operating in the system itself. In isothermal situations, the OVP takes the following form: given a set of slow state variables $\mathbf{x}_i, i=1, \dots, n$, the dynamics of the system is described by the thermodynamic fluxes which minimize the Onsager functional $\mathcal{O}(\dot{\mathbf{x}}_i)=\Phi(\dot{\mathbf{x}}_i)+\dot{E}(\mathbf{x}_i,\dot{\mathbf{x}}_i)$, where $\Phi$ is the dissipation function, which is quadratic in $\dot{\mathbf{x}}_i$ as a near--equilibrium approximation, and $E$ is the free energy of the system. In our case, we minimize \eqref{eqn:12} with respect to the variables $\mathbf{v}_c$ and $\mathbf{v}_l$, thus obtaining the momentum balance equations for the two phases of the mixture as linear fluxes--forces relations, to be supplemented to the mass balance equations \eqref{eqn:1a}-\eqref{eqn:1c} in the mixture dynamics description. We thus rewrite \eqref{eqn:12} as 
\begin{align}
\label{eqn:13}
&\int_{\Omega}\frac{\delta E}{\delta \phi_c}\biggl(\frac{\Gamma_c}{\gamma}-\text{div}(\phi_c \mathbf{v}_c)\biggr)d\mathbf{x}+\int_{\Omega}\frac{\delta E}{\delta n}\biggl(S-\text{div}(n \mathbf{v})-\text{div}\mathbf{J_n}\biggr)d\mathbf{x}+\frac{1}{k}\int_{\Omega}\mathbf|v|^2 d\mathbf{x} \\
& \notag+\int_{\Omega}b(\phi_c)\mathbf{T}\nabla \mu\cdot \nabla \mu d\mathbf{x} +\int_{\Omega}D(\phi_c) \mathbf{D}\nabla \eta \cdot \nabla \eta d\mathbf{x}=\int_{\Omega}\frac{\Gamma_c}{\gamma}\mu d\mathbf{x}+\int_{\Omega}S\eta d\mathbf{x},
\end{align}
which gives 
\begin{equation}
\label{eqn:14}
\int_{\Omega}\mu \left(-\text{div}(\phi_c \mathbf{v}_c)\right)d\mathbf{x}+\int_{\Omega}\eta \left(-\text{div}(n \mathbf{v})\right)d\mathbf{x}+\frac{1}{k}\int_{\Omega}\mathbf|v|^2 d\mathbf{x} +\int_{\Omega}b(\phi_c)\mathbf{T}\nabla \mu\cdot \nabla \mu d\mathbf{x}=0.
\end{equation}
As in the framework of the OVP, we assume that the dissipation term due to the viscous interactions between the phases can be assumed to depend only on the drag between the phases and can be written as a quadratic term in the relative velocity between the two phases
\begin{equation}
\label{eqn:15}
\int_{\Omega}b(\phi_c)\mathbf{T}\nabla \mu \cdot \nabla \mu d\mathbf{x}\equiv \int_{\Omega}K(\phi_c)\mathbf{T}^{-1}(\mathbf{v_c}-\mathbf{v_l})\cdot (\mathbf{v_c}-\mathbf{v_l}) d\mathbf{x} \to
\mathbf{v_c}-\mathbf{v_l}=-\sqrt{\frac{b(\phi_c)}{K(\phi_c)}}\mathbf{T}\nabla \mu,
\end{equation}
where $K(\phi_c)$ is a friction function to be specified later. Sobstituting \eqref{eqn:15} in \eqref{eqn:14}, we find that a minimum of \eqref{eqn:14} with respect to the phase velocities $\mathbf{v_c}$ and $\mathbf{v_l}$ satisfies the following first-order conditions
\begin{equation}
\label{eqn:16}
\begin{cases}
\phi_c \nabla \mu +\phi_c n \nabla \eta + 2\frac{\phi_c}{k} \mathbf{v}+K(\phi_c)\mathbf{T}^{-1}(\mathbf{v_c}-\mathbf{v_l})=0,\\
\phi_l n \nabla \eta + 2\frac{\phi_l}{k} \mathbf{v}-K(\phi_c)\mathbf{T}^{-1}(\mathbf{v_c}-\mathbf{v_l})=0,
\end{cases}
\end{equation}
which, after the elimination of the variable $\mathbf{v}$, gives
\begin{equation}
\label{eqn:17}
\mathbf{v_c}-\mathbf{v_l}=-\frac{\phi_c(1-\phi_c)}{K(\phi_c)}\mathbf{T}\nabla \mu.
\end{equation}
Comparing \eqref{eqn:17} with \eqref{eqn:15}, we get that 
\begin{equation}
\label{eqn:18}
b(\phi_c)=\frac{\phi_c^2(1-\phi_c)^2}{K(\phi_c)}.
\end{equation}
Note that the OVP does not impose any constraint on the mobility function $D(\phi_c)$. 

For what concerns the friction function $K(\phi_c)$, it can be expressed as 
\begin{equation}
\label{eqn:19}
K(\phi_c)=\frac{(1-\phi_c)\rho_c}{k(\phi_c)},
\end{equation}
where $\rho_c$ is the viscosity of the liquid phase and $k(\phi_c)$ is the intrinsic permeability of the tumor phase. Indeed, we can consider the model dynamics as describing a Darcy--like flow of the viscous fluid associated to the phase $l$ through the porous--permeable solid matrix associated to the soft material of the phase $c$ \cite{holmes}, and \eqref{eqn:15} represents the momentum exchange between the two phases. A possible expression for the intrinsic permeability can be derived by assuming that the tumor tissue consists of homogeneous and isotropic parallel cylindrical pores, and using the Poiseuille formula for a capillary tube \cite{holmes}
\begin{equation}
\label{eqn:20}
k(\phi_c)=\frac{r^2(1-\phi_c)}{8\delta^2},
\end{equation}
where $r$ is the effective radius of the pores and $\delta$ the tortuosity factor. This gives $K(\phi_c)=M$, where $M$ is a generic friction constant, and 
\begin{equation}
\label{eqn:21}
b(\phi_c)=\frac{\phi_c^2(1-\phi_c)^2}{M}.
\end{equation} 
\begin{rem}
\label{rem:2}
We observe that in \cite{agosti2} the form $K(\phi_c)=M\phi_c$ was assumed, which lead to 
\begin{equation}
\label{eqn:22a}
b(\phi_c)=\frac{\phi_c(1-\phi_c)^2}{M}.
\end{equation} 
This form for the friction function $K$ can be obtained by considering a variable fluid viscosity $\rho_c \propto (1-\phi_c)$ in \eqref{eqn:19} and a Kozeny--Carman permeability--porosity law \cite{kozeny,carman}
\[
k(\phi_c)\propto \frac{(1-\phi_c)^2}{\phi_c},
\] 
which can be obtained from \eqref{eqn:20} by substituting the effective radius $r$  with the hydraulic radius $r_h=\frac{1-\phi_c}{s_h \phi_c}$, where $s_h$ is the ratio of the exposed interfacial surface with the solid volume. The Kozeny--Carman law is tipically used to model the Darcy flow of a fluid in a porous permeable rock, where the hydraulic radius is an effective definition for the representative radius of the pores.
\end{rem}
On account of the possible forms of the mobility function \eqref{eqn:21}{-}\eqref{eqn:22a}, obtained by considering different modelizations of the friction function $K(\phi_c)$, in the following we will consider a general expression for the mobility of the form
\begin{equation}
\label{eqn:23a}
b(\phi_c)=\frac{\phi_c^{\alpha}(1-\phi_c)^2}{M}, \quad \alpha \geq 0.
\end{equation}
 {We observe that the expression \eqref{eqn:23a} with $\alpha=1$ is the standard form for the mobility considered for the description of the tumor growth dynamics in literature (see e.g. \cite{chat,agosti2}). Here we generalize this expression by introducing a parameter $\alpha\geq 0$, related to mechanical parameters which enter in the constitutive laws for the filtration processes between the mixture's components and which may vary for different type of biological cells.}
 
In the following, we will consider the analysis of the simplified model problem when $\mathbf{v}\equiv 0$, $\Gamma_c\equiv 0$ and $n\equiv 0$, i.e.
\begin{equation}
\label{eqn:22}
\begin{cases}
\frac{\partial \phi}{\partial t}-\text{div}\biggl(b(\phi)\nabla \mu\biggr)=0 ,\\
\mu=\frac{\psi'(\phi)}{\epsilon}+\epsilon \left(J\star 1\right)\phi-\epsilon J\star \phi,
\end{cases}
\end{equation} 
in $\Omega$, where we omit to write the subscript $c$ for the cells phase and for simplicity we {have chosen} $\mathbf{T}=\mathbf{Id}$, endowed with the homogeneous boundary condition
\begin{equation}
\label{eqn:23}
b(\phi)\nabla \mu \cdot \boldsymbol{\nu}|_{\partial \Omega}=0.
\end{equation}
Moreover, for simplicity we will take $\epsilon=M=1$ in the analysis.
{Let us notice that the techniques we employ in the next sections in order to prove the existence of solutions and the separation property for equation \eqref{eqn:22} are a starting point in order then to study the complete coupled system \eqref{eqn:11} at least under proper assumptions on the source terms $\Gamma_c$ and $S$.}

\section{Existence of a weak solution}
\label{exi}

We recall the following existence theorem for the non--degenerate mobility case with smooth potential reported e. g. in \cite{frigeri1}, {valid for any spatial dimension $d\leq 3$,} {which will be useful in the proof of existence of weak solutions in case of singular potential and degenerate mobility (cf.~Thm.~\ref{thm:2}).}
\begin{thm}
\label{thm:1}
Let $\phi_0\in H$ such that $\psi(\phi_0)\in L^1(\Omega)$, and suppose that the following properties are satisfied:
\begin{itemize}
\item[$(H_1)$] $b\in C_{\text{loc}}^{0,1}(\mathbb{R})$, and there exist $b_1, b_2>0$ such that
\[
b_1\leq b(s)\leq b_2, \quad \forall s \in \mathbb{R};
\]
\item[$(H_2)$] $J(\cdot - \mathbf{x})\in W^{1,1}(\Omega)$ for almost any $\mathbf{x}\in \Omega$ and satisfies
\[
J(\mathbf{x})=J(-\mathbf{x}), \quad a(\mathbf{x}):=J\star 1\geq 0, \; \text{a.e.} \, \mathbf{x}\in \Omega,
\]
\[
a^*:=\sup_{\mathbf{x}\in \Omega}\int_{\Omega}|J(\mathbf{x}-\mathbf{y})|\,d\mathbf{y}<+\infty, \quad b^*:=\sup_{\mathbf{x}\in \Omega}\int_{\Omega}|\nabla J(\mathbf{x}-\mathbf{y})|\,d\mathbf{y}<+\infty;
\]
\item[$(H_3)$] $\psi\in C_{\text{loc}}^{2,1}(\mathbb{R})$, and there exist $c_0>0$ such that
\[
\psi^{''}(s)+a(\mathbf{x})\geq c_0, \quad \forall s \in \mathbb{R}, \; \text{a.e.} \, \mathbf{x}\in \Omega;
\]
\item[$(H_4)$] There exist $c_1>\frac{a^{*}-a_{*}}{2}$, where $a_{*}:=\inf_{\mathbf{x}\in \Omega}J\star 1$, and $c_2\in \mathbb{R}$ such that
\[
\psi(s)\geq c_1 s^2 - c_2, \; \forall s \in \mathbb{R};
\]
\item[$(H_5)$] There exist $c_3>0$, $c_4\geq 0$ and $r\in (1,2]$ such that
\[
|\psi'(s)|^r\leq c_3 |\psi(s)| + c_4, \; \forall s \in \mathbb{R}.
\]
\end{itemize}
Then, for every given $T>0$, there exists a weak solution $\phi$ to \eqref{eqn:22}-\eqref{eqn:23}, with
\begin{align}
\label{eqn:23bis}
& \phi \in L^{\infty}(0,T;H) \cap L^2(0,T;V), \quad \phi_t \in L^2(0,T;V');\\
\notag &  \mu:=a\phi - J\star \phi + \psi'(\phi) \in L^2(0,T;V),
\end{align}
satisfying the following weak formulation
\begin{equation}
\label{eqn:24}
<\phi_t,\xi>+\left(b(\phi) \nabla \mu, \nabla \xi \right)=0, \quad \forall \xi \in V, \; for \; \text{a.e.} \; t \in (0,T),
\end{equation}
together with the initial condition $\phi(\mathbf{x}, t=0)=\phi_0(\mathbf{x})$ for $\mathbf{x}\in \Omega$. Moreover, the following energy inequality is satisfied, for every $t \in [0,T]$,
\begin{equation}
\label{eqn:25}
E(\phi(t))+\int_0^t||\sqrt{b(\phi)}\nabla \mu||^2\,d\tau \leq E(\phi_0),
\end{equation}
where 
\begin{equation}
\label{eqn:25bis}
E(\phi(t))=\frac{1}{4}\int_{\Omega}\int_{\Omega}J(\mathbf{x}-\mathbf{y})\left(\phi_c(\mathbf{x})-\phi_c(\mathbf{y})\right)^2d\mathbf{x}\,d\mathbf{y}+\int_{\Omega}\psi(\phi_c)\,d\mathbf{x}.
\end{equation}
\end{thm}
\begin{rem}
Taking $\xi=\phi$ in \eqref{eqn:24} and using \eqref{eqn:23bis}, it's possible to show (following the same calculations as in \cite[Corollary 1]{frigeri1}) that the weak solution constructed in Theorem \ref{thm:1} satisfies the Energy equality
\begin{align}
\label{eqn:25tris}
& \displaystyle \notag \frac{1}{2}||\phi||^2+\int_0^t\int_{\Omega}b(\phi) \psi''(\phi) |\nabla \phi|^2+\int_0^t\int_{\Omega}b(\phi) a |\nabla \phi|^2+ \int_0^t\int_{\Omega}b(\phi) \left(\phi \nabla a - \nabla J \star \phi \right)\cdot \nabla \phi=\\
& \displaystyle\frac{1}{2}||\phi_0||^2.
\end{align}
\end{rem}
We now give the existence result for solutions of equations \eqref{eqn:22}-\eqref{eqn:23} in the case of degenerate mobility and singular single-well potential. Since in the degenerate case the a-priori estimates do not provide a bound in some $L^p$ space for the gradient of the chemical potential, the weak solutions of the equations will be defined to satisfy a primal formulation in which the variable $\mu$ does not appear. This is standard for degenerate Cahn--Hilliard equations, both in the local and in the non-local case (see e.g. \cite{garcke,frigeri1}). Since the proof of the existence result will make use of Entropy estimates, in particular to obtain compactness for approximation sequences of the solution, we introduce the Entropy function $\Phi \in C^2(0,1)$, defined by the relation $\Phi''(r)=\frac{1}{b(r)}$ for all $r\in (0,1)$, with $\Phi'(A)=\Phi(A)=0$, $A\in (0,1)$. 
{
Given \eqref{eqn:23a}, we observe that
\begin{equation}
    \label{eqn:25Phi0}
    \Phi(r)\sim 
    \begin{cases}
        r^{2-\alpha} \quad \text{if}\; \alpha\neq 2,\\
        \log(r)\quad \text{if}\; \alpha= 2,
    \end{cases}
    \text{as}\; r\to 0,
\end{equation}
while
\begin{equation}
  \label{eqn:25Phi1}
  \Phi(r)\sim \log(1-r)\quad \text{as}\; r\to 1.
\end{equation}
}
We give the following theorem.
\begin{thm}
\label{thm:2}
{Let $d\leq 3$.} Assume that $(H_2)$ is satisfied, that $b(\phi)$ has the particular form \eqref{eqn:23a} (with $M=1$) and that $\psi(\phi)$ is a single-well potential of the form \eqref{eqn:4}. Assume moreover that \eqref{eqn:7} is satisfied.
Let $\phi_0\in L^{\infty}(\Omega)$ such that $\psi(\phi_0)\in L^1(\Omega)$ and $\Phi(\phi_0)\in L^1(\Omega)$. 
Then, for every given $T>0$, there exists a weak solution $\phi$ to \eqref{eqn:22}-\eqref{eqn:23}, with
\begin{align}
\label{eqn:37}
& \phi \in L^{\infty}(0,T;H) \cap L^2(0,T;V), \quad \phi_t \in L^2(0,T;V');\\
\notag &  \phi \in L^{\infty}(\Omega \times (0,T)), \quad 0\leq \phi(\mathbf{x},t) < 1 \quad \text{a. e.} \; (\mathbf{x},t) \in \Omega \times (0,T);\\
\notag & \int_{\Omega}\phi(\mathbf{x},t)=\int_{\Omega}\phi(\mathbf{x},t=0);
\end{align}
satisfying the following weak formulation, $\forall \xi \in V, \; for \; \text{a.e.} \; t \in (0,T)$,
\begin{equation}
\label{eqn:38}
<\phi_t,\xi>+\left(b(\phi) \psi''(\phi) \nabla \phi, \nabla \xi \right)+\left(b(\phi) a \nabla \phi, \nabla \xi \right) + \left(b(\phi) \left(\phi \nabla a - \nabla J \star \phi \right), \nabla \xi \right)=0,
\end{equation}
together with the initial condition $\phi(\mathbf{x}, t=0)=\phi_0(\mathbf{x})$ for $\mathbf{x}\in \Omega$. Moreover, the weak solution satisfies the Energy equality \eqref{eqn:25tris}.
\end{thm}
\begin{pf}
The proof follows the calculations in the proof of \cite[Theorem 2]{frigeri1}. Here, we adapt the proof to the specific case of a degenerate mobility of the form \eqref{eqn:23a} and a singular single-well potential of the form \eqref{eqn:4}. We observe that the latter potential does not satisfy the hypotheses $\mathbf{A}_2,\mathbf{A}_3$ in \cite[Theorem 2]{frigeri1}, which are standard hypotheses for double-well potentials. Hence, in the following we will report the details of the proof which necessitate an ad-hoc treatment due to the single-well form of the potential, referring the reader to \cite{frigeri1} for other details which follow closely the calculations therein.

We start by introducing a proper regularization of system \eqref{eqn:22}-\eqref{eqn:23}. Given a regularization parameter $\lambda \in (0,1)$, we define a regularized mobility
\begin{equation}
\label{eqn:26}
b_{\lambda}(r):=
\begin{cases}
b(\lambda) \quad \text{for} \; r\leq \lambda, \\
b(r) \quad \text{for} \;  \lambda < r <1-\lambda,\\
b(1-\lambda) \quad \text{for} \; r\geq 1-\lambda.
\end{cases}
\end{equation}
We introduce the following  convex splitting of the potential $\psi(\phi)=\psi_1(\phi)+\psi_2(\phi)$, where
\begin{equation}
\label{eqn:27}
\psi_1(\phi)=-(1-\bar{\phi})\log(1-\phi), \quad \psi_2(\phi)=-\frac{\phi^3}{3}-(1-\bar{\phi})\frac{\phi^2}{2}-(1-\bar{\phi})\phi,
\end{equation}
where the singular part is contained in the convex component, and moreover we define its regularization
\begin{equation}
\label{eqn:28}
\psi''_{1,\lambda}(r):=
\begin{cases}
\psi''_{1}(\lambda) \quad \text{for} \; r\geq 1-\lambda, \\
\psi''_{1}(r) \quad \text{for} \; r <1-\lambda,
\end{cases}
\end{equation}
with $\psi'_{1,\lambda}(1-\lambda)=\psi'_{1}(1-\lambda)$ and $\psi_{1,\lambda}(1-\lambda)=\psi_{1}(1-\lambda)$. Integrating \eqref{eqn:28} two times we thus get
\begin{equation}
\label{eqn:29}
\psi_{1,\lambda}(r)=
\begin{cases}
-(1-\bar{\phi})\log(\lambda)+\frac{3}{2}(1-\bar{\phi})-\frac{2}{\lambda}(1-\bar{\phi})(1-r)+\frac{1-\bar{\phi}}{2\lambda^2}(1-r)^2 \quad \text{for} \; r\geq 1-\lambda, \\
\psi_{1}(r) \quad \text{for} \; r <1-\lambda,
\end{cases}
\end{equation}
and
\begin{equation}
\label{eqn:30}
\psi'_{1,\lambda}(r)=
\begin{cases}
\frac{2}{\lambda}(1-\bar{\phi})-\frac{1-\bar{\phi}}{\lambda^2}(1-r) \quad \text{for} \; r\geq 1-\lambda, \\
\psi'_{1}(r) \quad \text{for} \; r <1-\lambda,
\end{cases}
\end{equation}
Finally, we define the continuous extension $\bar{\psi}_2\in C^2(\mathbb{R})$ of $\psi_2\in C^2([0,1])$ given by
\begin{equation}
\label{eqn:31}
\bar{\psi}_{2}(r)=
\begin{cases}
\psi_2(1)+\psi'_2(1)(r-1)+\frac{1}{2}\psi''_2(1)(r-1)^2 \quad \text{for} \; r\geq 1, \\
\psi_{2}(r) \quad \text{for} \; r <1,
\end{cases}
\end{equation}
and set $\psi_{\lambda}(r):=\psi_{1,\lambda}(r)+\bar{\psi}_2(r)$ for $r\in \mathbb{R}$. 

We then consider the regularized problem
\begin{equation}
\label{eqn:32}
\begin{cases}
\frac{\partial \phi_{\lambda}}{\partial t}-\text{div}\left(b_{\lambda}(\phi_{\lambda})\nabla \mu_{\lambda}\right)=0 ,\\
\mu_{\lambda}=\psi_{\lambda}'(\phi_{\lambda})+ \left(J\star 1\right)\phi_{\lambda}- J\star \phi_{\lambda},
\end{cases}
\end{equation} 
in $\Omega$, endowed with the homogeneous boundary condition
\begin{equation}
\label{eqn:33}
\nabla \mu_{\lambda} \cdot \boldsymbol{\nu}|_{\partial \Omega}=0,
\end{equation}
and initial condition 
\begin{equation}
\label{eqn:34}
\phi_{\lambda}(\mathbf{x},t=0)=\phi_0(\mathbf{x}), \quad \forall\, \mathbf{x}\in \Omega.
\end{equation}
From \eqref{eqn:29} and \eqref{eqn:31} we observe that
\begin{equation}
\label{eqn:29+31}
\psi_{\lambda}(r)=-(1-\bar{\phi})\log(\lambda)-\frac{1}{3}+\left(\frac{2}{\lambda}(1-\bar{\phi})+(2\bar{\phi}-3)\right)(r-1)+\left(\frac{1-\bar{\phi}}{2\lambda^2}+\frac{\bar{\phi}-3}{2}\right)(r-1)^2,
\end{equation}
for $r\geq 1$, with the terms in the round brackets being positive if
\[
\lambda \leq \lambda_0:=\min\left(\sqrt{\frac{1-\bar{\phi}}{3-\bar{\phi}}},\frac{2-2\bar{\phi}}{3-2\bar{\phi}}\right).
\] 
We also observe that, for $r\leq 0$, $\psi_{\lambda}(r)\geq 0$. Moreover, for $r\leq 0$,
\begin{align*}
& \psi_{\lambda}(r)=-(1-\bar{\phi})\log(1+|r|)+\frac{|r|^3}{3}-(1-\bar{\phi})\frac{r^2}{2}+(1-\bar{\phi})|r|\geq \\
&-(1-\bar{\phi})|r|+\frac{|r|^3}{3}-(1-\bar{\phi})\frac{|r|^3}{3}-\frac{1-\bar{\phi}}{6}+(1-\bar{\phi})|r|\geq \bar{\phi}|r|^3-C \geq \frac{3\bar{\phi}}{2\gamma}r^2-C-\frac{\bar{\phi}}{2\gamma^3},
\end{align*}
for a generic $\gamma>0$ (where we have used the Young inequality). Hence, for sufficiently small $\gamma$, we have that
\begin{equation}
\label{eqn:34tris}
\psi_{\lambda}(r)\geq c_1 r^2-c_2, \quad \forall\, r \in \mathbb{R},\quad \forall\, \lambda \leq \lambda_0,
\end{equation}
with $c_1>\frac{a^*-a_*}{2}$, $c_2\in \mathbb{R}$. Moreover, thanks to \eqref{eqn:30} and \eqref{eqn:31}, we have that there exist $c_3>0$, $c_4\geq 0$, such that
\[
|\psi'_{\lambda}(r)|^{\frac{3}{2}}\leq c_3|\psi_{\lambda}(r)|+c_4, \quad \forall r \in \mathbb{R}.
\]
We also observe that, for $\lambda$ sufficiently small, $\psi_{\lambda}(r)\leq \psi(r)$, and that, thanks to assumption \eqref{eqn:7}, there exists $c_0>0$ such that 
\begin{equation}
\label{eqn:34bis}
\psi_{\lambda}^{''}(s)+a(\mathbf{x})\geq c_0, \quad \forall s \in \mathbb{R}, \; \text{a.e.} \, \mathbf{x}\in \Omega.
\end{equation}
Hence, the hypothesis $(H_1)-(H_5)$ of Theorem \eqref{thm:1} are satisfied, with also $\psi_{\lambda}(\phi_0)\in L^1(\Omega)$, and thus there exists a weak solution $\phi_{\lambda}$ of \eqref{eqn:32}--\eqref{eqn:33} satisfying \eqref{eqn:23bis}, \eqref{eqn:24} and the Energy inequality \eqref{eqn:25}.
\begin{rem}
\label{rem:3}
The hypothesis $(H_5)$ in Theorem \ref{thm:1}, which is needed to control the $L^1(\Omega)$ norm of $\mu$ and to pass to the limit in the $\psi'$ term in $\mu$ in the context of a Galerkin approximation of the problem, can be substituted in our framework by the property that
\[
|\psi'_{\lambda}(r)|\leq C + C|r|^2, \quad \forall\, r \in \mathbb{R},
\] 
and by the fact that, from the energy inequality, $\phi_{\lambda}\in L^{\infty}(0,T;H)$.
\end{rem}
From the Energy inequality \eqref{eqn:25} and from \eqref{eqn:34tris} we obtained the uniform (in $\lambda$) bounds
\begin{align}
\label{eqn:boundphih}
&||\phi_{\lambda}||_{L^{\infty}(0,t;H)}\leq C,\\
\label{eqn:boundmu}
&||\sqrt{b_{\lambda}(\phi_{\lambda})}\mu_{\lambda}||_{L^{2}(0,t;H)}\leq C.
\end{align}

We now introduce the Entropy function $\Phi_{\lambda} \in C^2(\mathbb{R})$ such that $\Phi''_{\lambda}(r)=\frac{1}{b_{\lambda}(r)}$, with $\Phi'_{\lambda}(A)=\Phi_{\lambda}(A)=0$, $A\in (0,1)$. We note that
{
\begin{equation}
\label{eqn:34_4}
\Phi_{\lambda}(r)\leq \Phi(r), \quad \forall r \in [0,1], 
\end{equation}
for $\lambda$ sufficiently small. Moreover, 
\begin{equation}
    \label{eqn:34_5}
    \Phi_{\lambda}(r)\geq 0, \;\; \forall r\in \mathbb{R}, \quad \Phi_{\lambda}(r)\equiv \Phi(r), \;\; \text{for} \; r\in (\delta,1-\delta).
\end{equation}}
We also note that
\begin{align}
\label{eqn:35}
&\Phi_{\lambda}(r)\geq \frac{1}{2}\frac{1}{b(1-\lambda)}(r-1)^2, \quad \text{for} \; r\geq 1,\\
& \notag \Phi_{\lambda}(r)\geq \frac{1}{2}\frac{1}{b(\lambda)}r^2, \quad \text{for} \; r\leq 0.
\end{align}
Taking $\xi \equiv  \Phi'_{\lambda}(\phi_{\lambda})$ in \eqref{eqn:24} (which is a suitable test function due to the boundedness of $\Phi''_{\lambda}$ and to \eqref{eqn:23bis}), and using \eqref{eqn:23bis}, we obtain that
\begin{align}
\label{eqn:36}
& \displaystyle \notag \frac{d}{dt}\int_{\Omega}\Phi_{\lambda}(\phi_{\lambda}(\mathbf{x},t))\,d\mathbf{x}+\int_{\Omega}\nabla \mu_{\lambda}\cdot \nabla \phi_{\lambda}=\frac{d}{dt}\int_{\Omega}\Phi_{\lambda}(\phi_{\lambda}(\mathbf{x},t))\,d\mathbf{x}+
\int_{\Omega}\left(a+\psi''(\phi_{\lambda})\right) |\nabla \phi_{\lambda}|^2+\\
& \displaystyle \int_{\Omega}\left(\phi_{\lambda} \nabla a - \nabla J \star \phi_{\lambda} \right)\cdot \nabla \phi_{\lambda}=0.
\end{align}
Then, thanks to \eqref{eqn:34bis}, to the hypotesis $(H_2)$ and \eqref{eqn:boundphih}, using moreover the Cauchy--Schwarz and the Young inequalities, we get
\begin{equation}
\label{eqn:36bis}
\displaystyle \frac{d}{dt}\int_{\Omega}\Phi_{\lambda}(\phi_{\lambda}(\mathbf{x},t))\,d\mathbf{x}+\frac{c_0}{2}||\nabla \phi_{\lambda}||^2\leq C||\phi_{\lambda}||^2\leq C.
\end{equation}
Given \eqref{eqn:34_4} and the assumption that $\Phi(\phi_0)\in L^1(\Omega)$, integrating in time the inequality \eqref{eqn:36bis} between $0$ and $t$, for any $t\in (0,T)$, we obtain the uniform (in $\lambda$) bounds
\begin{align}
\label{eqn:boundphiv}
&||\phi_{\lambda}||_{L^{2}(0,t;V)}\leq C,\\
\label{eqn:boundPhi}
&||\Phi_{\lambda}(\phi_{\lambda})||_{L^{\infty}(0,t;L^1(\Omega))}\leq C.
\end{align}
It's easy to deduce from \eqref{eqn:24} by comparison, using the bound \eqref{eqn:boundmu}, that 
\begin{align}
\label{eqn:boundphip}
&||\partial_t\phi_{\lambda}||_{L^{2}(0,t;V')}\leq C.
\end{align}
Collecting the results \eqref{eqn:boundphih}, \eqref{eqn:boundmu}, \eqref{eqn:boundphiv} and \eqref{eqn:boundphip}, which are uniform in $\lambda$, from the Banach--Alaoglu and the Aubin--Lions lemma, we finally obtain the convergence properties, up to subsequences of the solutions, which we still label by the index $\lambda$, as follows:
\begin{align}
\label{eqn:conv3d1} & \phi_{\lambda} \overset{\ast}{\rightharpoonup} \phi \quad \text{in} \quad L^{\infty}(0,T;H),\\
\label{eqn:conv3d2} & \phi_{\lambda} {\rightharpoonup} \phi \quad \text{in} \quad L^{2}(0,T;V),\\
\label{eqn:conv3d3} & \partial_t \phi_{\lambda} \rightharpoonup \partial_t \phi \quad \text{in} \quad L^{2}(0,T;V'),\\
\label{eqn:conv3d4} & \phi_{\lambda} \to \phi \quad \text{in} \quad C^{0}(0,T;H) \cap L^{2}(0,T;L^p(\Omega)), \;\; \text{and} \;\; \text{a.e. in} \; \; \Omega_T,
\end{align}
as $\lambda \to 0$, {with $p\in [1,6)$ for $d=3$, $p\in [1,\infty)$ for $d=2$. For $d=1$, $\phi_{\lambda} \to \phi$ uniformly in $\Omega_T$ (see e.g. \cite{agosti1})}.
Similarly to \cite{frigeri1}, we may use \eqref{eqn:boundPhi}, \eqref{eqn:35} and the Lebesgue convergence theorem to deduce that
\begin{equation}
\label{eqn:conv3d5}
0\leq \phi \leq 1 \quad \text{a.e.}\; (\mathbf{x},t)\in \Omega \times (0,T).
\end{equation}
Finally, thanks to \eqref{eqn:conv3d1}- \eqref{eqn:conv3d5} we can pass to the limit for $\lambda \to 0$ in the weak formulation of the regularized problem \eqref{eqn:32}--\eqref{eqn:33}, following the same procedure as in \cite{frigeri1}, and obtain that the limit satisfies \eqref{eqn:38} and, taking $\xi=\phi$ in \eqref{eqn:38}, that the limit solution satisfies the Energy equality \eqref{eqn:25tris}.

We are now left to prove that $\phi(\mathbf{x},t)<1$ for a.e. $(\mathbf{x},t)\in \Omega \times (0,T)$. We proceed by contradiction, following the same arguments as in \cite{agosti1}. Suppose that there exists a set $\mathcal{V}\subset \Omega \times (0,T)$, with positive measure, such that $\phi(x,t)\equiv 1$ on $\mathcal{V}$. Recalling the strong convergence \eqref{eqn:conv3d4}, we can apply the Egorov Theorem and obtain that, for each $\delta > 0$, there exists a subset $\mathcal{V}'$, with $|\mathcal{V}\setminus \mathcal{V}'|<\delta$, on which (a subsequence of) $\phi_{\lambda}$ converges uniformly to $\phi$. Hence we get that there exists a modulus of continuity $\delta_1(\epsilon)$, with $\delta_1(0)=0$, and $\lambda_0>0$ such that 
\[
\phi_{\lambda}\geq 1 -\delta_1(\epsilon) \quad \text{on} \; \mathcal{V}', \; \lambda\leq \lambda_0.
\]
From \eqref{eqn:29} we have that
\[
\psi _{1,\lambda}(\phi_{\lambda}(x,t))\geq -(1-\bar{\phi})\log \delta_1(\epsilon) -\frac{2}{\delta_1(\epsilon)}(1-\bar{\phi}), \; \forall (x,t) \in \mathcal{V}'.
\]
Using the Beppo Levi Theorem and the fact that $\psi_{1,\lambda}(r)\geq 0$ for $r\geq 0 $ if $\lambda$ is sufficiently small, we get that
\[
\limsup_{\lambda \to 0}\int_{\Omega}\psi_{1,\lambda}(\phi_{\lambda}(x,\bar{t}))dx\geq \limsup_{\lambda \to 0}\biggl(-(1-\bar{\phi})\log \delta_1(\epsilon) -\frac{2}{\delta_1(\epsilon)}(1-\bar{\phi})\biggr)|\mathcal{V}'|\to \infty,
\]
which contradicts \eqref{eqn:25}. Hence, $0\leq \phi < 1$ almost everywhere in $\Omega\times (0,T)$.
\end{pf}
 \begin{rem}
 \label{rem:4}
 {
 In the case $\alpha \geq 2$, we observe form \eqref{eqn:25Phi0}-\eqref{eqn:25Phi1} that the assumption $\Phi(\phi_0)\in L^1(\Omega)$ of Theorem \ref{thm:2} implies that $0<\phi_0<1$, i.e. that the initial condition corresponds to a mixed configuration with no pure phases. Moreover, if $\alpha \geq 2$ we can prove that $0<\phi(\mathbf{x},t)$ for 
 almost every $(\mathbf{x},t) \in \Omega \times [0,T]$. 
 Indeed, using \eqref{eqn:34_5}, the Fatou lemma and \eqref{eqn:boundPhi} we have that 
 \[
 \int_{\Omega}\lim_{\lambda \to 0}\Phi_{\lambda}(\phi_{\lambda})\leq \lim_{\lambda \to 0}\int_{\Omega}\Phi_{\lambda}(\phi_{\lambda})\leq C,
 \]
 for a.a. $t\in (0,T)$.
 Moreover, thanks to \eqref{eqn:34_5} and \eqref{eqn:25Phi0}, we have that
 \[
 \Phi_{\lambda}(\phi_{\lambda})\geq \min\{\Phi(\lambda),\Phi(\phi_{\lambda})\}\to \infty,
 \]
 if $\lim_{\lambda\to 0}\phi_{\lambda}(t,\mathbf{x})=1$, hence from the latter relations we conclude that the set $\{\mathbf{x}| \phi(\mathbf{x},t)=0\}$ has zero measure for almost all $t\in [0,T]$.}
 \end{rem}

 \section{Strict separation property and uniqueness result}
 \label{separation}
 
 We now give a strict separation result, obtained from a Moser--Alikakos iteration argument. As a consequence of the latter result, we will prove the uniqueness of the weak solution defined in Theorem \ref{thm:2}. {This result is obtained in the most restrictive case of $d=3$, but the Moser--Alikakos iteration argument remains valid also for $d\leq 3$, using in particular the less restrictive Sobolev embeddings \eqref{eqn:gagliardoniremberg} for $d<3$ in the process.}. In the following, we will prove the separation result by a formal argument, justifying the calculations rigorously in a remark. 
 \begin{thm}
\label{thm:3}
Let {$d=3$} and the hypothesis of Theorem \ref{thm:2} be satisfied. {Let moreover $0<\phi_0<1$}, and let $\phi$ be a weak solution of \eqref{eqn:38}. Moreover, let us assume that $\alpha \geq 2$ and that $a_{*}>1-\bar{\phi}$. Then, there exists a positive constant $\delta$ such that
\begin{equation}
\label{eqn:39}
\delta \leq \phi(\mathbf{x},t)\leq 1-\delta, \quad \text{for a.~e.} \; (\mathbf{x},t)\in \Omega \times [0,T].
\end{equation}
\end{thm}
\begin{pf}
The proof is based on a Moser--Alikakos iteration argument, following \cite{laurencot}. In particular, we derive the proof for the strict separation result $\delta \leq \phi(\mathbf{x},t)$. Once this result is obtained, we use a result reported in \cite{frigeri2} to show that $\phi(\mathbf{x},t)\leq 1-\delta$.

We start by deriving formal estimates for inverse powers of $\phi$ starting from $\eqref{eqn:38}$ and taking $\xi=-\phi^{-l-2}$, with $l\in [0,\infty)$. Rigorous estimates could be obtained by considering proper truncations of the test functions, as described in Remark \ref{rem:rig}.
Recalling that
\[
\psi_1''(\phi)=\frac{1-\bar{\phi}}{(1-\phi)^2}, \quad \psi_2''(\phi)=-2\phi-(1-\bar{\phi}),
\]
and defining $h:=1-\bar{\phi}$, we have
\begin{align}
\label{eqn:40}
&<\phi_t,-\phi^{-l-2}>+h(l+2)\left(\phi^{\alpha-l-3} \nabla \phi, \nabla \phi \right)=(l+2)\left((1-\phi)^2 \phi^{\alpha-l-3} \nabla J\star \phi, \nabla \phi \right) \\
& \notag + 2(l+2) \left((1-\phi)^2 \phi^{\alpha-l-2} \nabla \phi, \nabla \phi \right) + h(l+2)\left((1-\phi)^2\phi^{\alpha-l-3} \nabla \phi, \nabla \phi \right) \\
& \notag-(l+2)\left((1-\phi)^2 a \phi^{\alpha-l-3} \nabla \phi, \nabla \phi \right)-(l+2)\left((1-\phi)^2\phi^{\alpha-l-2} \nabla J\star 1, \nabla \phi \right).
\end{align}
Since $-\phi^{-l-2}$ is monotone and due to \eqref{eqn:37}, we have that
\[<\phi_t,-\phi^{-l-2}>=\frac{1}{l+1}\frac{d}{dt}(\phi^{-l-1},1).\]
Expanding the factor $(1-\phi)^2$ in the third and fourth terms on the right hand side of \eqref{eqn:40}, integrating in time from $0$ to $T$ and rearranging terms, we get
\begin{align}
\label{eqn:41}
&(\phi^{-l-1},1)+\frac{4a_*(l+1)(l+2)}{(\alpha-l-1)^2}\int_0^T\int_{\Omega}\left|\nabla \phi^{\frac{\alpha-l-1}{2}}\right|^2\,d\mathbf{x}\,dt+\\
& \notag\frac{4(a_*-h)(l+1)(l+2)}{(\alpha-l+1)^2}\int_0^T\int_{\Omega}\left|\nabla \phi^{\frac{\alpha-l+1}{2}}\right|^2\,d\mathbf{x}\,dt\leq\\
& \notag (\phi_0^{-l-1},1)+\underbrace{\frac{2(l+1)(l+2)}{\alpha-l-1}\int_0^T\int_{\Omega}(1-\phi)^2\phi^{\frac{\alpha-l-3}{2}}\nabla J \star \phi \cdot \nabla \phi^{\frac{\alpha-l-1}{2}}\,d\mathbf{x}\,dt}_{I_1}\\
& \notag + \underbrace{\frac{8(l+1)(l+2)}{(\alpha-l)^2}\int_0^T\int_{\Omega}\left((1-\phi)^2+(a-h)\right)\left|\nabla \phi^{\frac{\alpha-l}{2}}\right|^2\,d\mathbf{x}\,dt}_{I_2}\\
& \notag - \underbrace{\frac{2(l+1)(l+2)}{\alpha-l-1}\int_0^T\int_{\Omega}(1-\phi)^2\phi^{\frac{\alpha-l-1}{2}}\nabla J \star 1 \cdot \nabla \phi^{\frac{\alpha-l-1}{2}}\,d\mathbf{x}\,dt}_{I_3}.
\end{align}
Using \eqref{eqn:37}, $(H_2)$, the fact that $\alpha\geq 2$, the Cauchy--Schwarz and Young inequalities
\begin{align*}
& I_1\leq \frac{2b(l+1)(l+2)}{\alpha-l-1}\int_0^T\left(\phi^{\alpha-l-3},1\right)^{\frac{1}{2}}  \left(\nabla \phi^{\frac{\alpha-l-1}{2}},\nabla \phi^{\frac{\alpha-l-1}{2}}\right)^{\frac{1}{2}}\,d\mathbf{x}\,dt\\
& \leq \frac{2b(l+1)(l+2)}{\alpha-l-1}\int_0^T\left(\phi^{-l-1},1\right)^{\frac{1}{2}}  \left(\nabla \phi^{\frac{\alpha-l-1}{2}},\nabla \phi^{\frac{\alpha-l-1}{2}}\right)^{\frac{1}{2}}\,d\mathbf{x}\,dt\\
& \leq \frac{a_*(l+1)(l+2)}{(\alpha-l-1)^2}\int_0^T\int_{\Omega}\left|\nabla \phi^{\frac{\alpha-l-1}{2}}\right|^2\,d\mathbf{x}\,dt+\frac{b^2(l+1)(l+2)}{a_*}\int_0^T\left(\phi^{-l-1},1\right)\,dt,
\end{align*}
and
\[
|I_3|\leq \frac{a_*(l+1)(l+2)}{(\alpha-l-1)^2}\int_0^T\int_{\Omega}\left|\nabla \phi^{\frac{\alpha-l-1}{2}}\right|^2\,d\mathbf{x}\,dt+\frac{b^2(l+1)(l+2)}{a_*}\int_0^T\left(\phi^{-l-1},1\right)\,dt.
\]
Using \eqref{eqn:37} and $(H_2)$ we moreover have
\[
I_2\leq \frac{C(l+1)(l+2)}{(\alpha-l)^2}\int_0^T\int_{\Omega}\left|\nabla \phi^{\frac{\alpha-l}{2}}\right|^2\,d\mathbf{x}\,dt.
\]
We thus obtain
\begin{align}
\label{eqn:42}
&\left(\frac{1}{\phi}^{l+1},1\right)+\frac{2a_*(l+1)(l+2)}{(\alpha-l-1)^2}\int_0^T\int_{\Omega}\left|\nabla \phi^{\frac{\alpha-l-1}{2}}\right|^2\,d\mathbf{x}\,dt\\
& \notag+\frac{4(a_*-h)(l+1)(l+2)}{(\alpha-l+1)^2}\int_0^T\int_{\Omega}\left|\nabla \phi^{\frac{\alpha-l+1}{2}}\right|^2\,d\mathbf{x}\,dt\\
& \notag \leq \left(\frac{1}{\phi_0}^{l+1},1\right)+C(l+1)(l+2)\int_0^T\left(\frac{1}{\phi}^{l+1},1\right)\,dt\\
& \notag+\frac{C(l+1)(l+2)}{(\alpha-(l-1)-1)^2}\int_0^T\int_{\Omega}\left|\nabla \phi^{\frac{\alpha-(l-1)-1}{2}}\right|^2\,d\mathbf{x}\,dt.
\end{align}
\begin{rem}
\label{rem:5}
The inequality \eqref{eqn:42} gives, for $l=0$,
\begin{align}
\label{eqn:43}
& \left(\frac{1}{\phi},1\right)+\frac{4a_*}{(\alpha-1)^2}\int_0^T\int_{\Omega}\left|\nabla \phi^{\frac{\alpha-1}{2}}\right|^2\,d\mathbf{x}\,dt \leq \left(\frac{1}{\phi_0},1\right)+C\int_0^T\left(\frac{1}{\phi},1\right)\,dt\\
& \notag+\frac{C}{4}\int_0^T\int_{\Omega}\phi^{\alpha-2}|\nabla \phi|^2\,d\mathbf{x}\,dt.
\end{align}
Hence, since $\alpha\geq 2$ and thanks to \eqref{eqn:37}, a Gronwall argument gives that
\begin{equation}
\label{eqn:44}
\frac{1}{\phi}\in L^{\infty}(0,T;L^1(\Omega)), \quad \phi^{\frac{\alpha-1}{2}}\in L^2(0,T;V).
\end{equation}
We note also that, for $l=-1$, the inequality \eqref{eqn:42} gives that $-\log \phi \in L^{\infty}(0,T;L^1(\Omega))$.
Then, starting from \eqref{eqn:43} and iterating over the index $l$ in \eqref{eqn:42}, for $l$ finite, we get that 
\begin{equation}
\label{eqn:45}
\frac{1}{\phi}\in L^{\infty}(0,T;L^p(\Omega)), \quad \phi^{\frac{\alpha-p}{2}}\in L^2(0,T;V), \quad \forall p\in(1,\infty).
\end{equation}
\end{rem}
We continue from \eqref{eqn:42}, rewriting the second term on the right hand side as 
\[
C(l+1)(l+2)\int_0^T\left(\frac{1}{\phi}^{l+1},1\right)\,dt=C(l+1)(l+2)\int_0^T\left(\left(\frac{1}{\phi}\right)^{\frac{\alpha}{2}}\left(\frac{1}{\phi}\right)^{\frac{l+1}{2}}\left(\frac{1}{\phi}\right)^{\frac{l+1-\alpha}{2}},1\right)\,dt.
\]
Then, using a trilinear H{\"o}lder inequality, \eqref{eqn:45} with $p=3\alpha$, the Gagliardo--Niremberg inequality \eqref{eqn:gagliardoniremberg} and the Young inequality, the fact that $0\leq \phi \leq 1$ and the inequality $\sqrt{a+b}\leq \sqrt{a}+\sqrt{b}$ for each $a,b \in \mathbb{R}^+$, we get
\begin{align*}
& C_1(l+1)(l+2)\int_0^T\left(\frac{1}{\phi}^{l+1},1\right)dt\leq \\
& \notag C(l+1)(l+2)\int_0^T\left|\left|\left(\frac{1}{\phi}\right)^{\frac{\alpha}{2}}\right|\right|_{L^6(\Omega)}\left(\int_{\Omega}\left(\frac{1}{\phi}\right)^{3(l+1-\alpha)}\right)^{\frac{1}{6}}\left(\int_{\Omega}\left(\frac{1}{\phi}\right)^{\frac{3(l+1)}{4}}\right)^{\frac{2}{3}}dt \leq \\
&C(l+1)(l+2)\int_0^T\left(\int_{\Omega}\left(\phi^{\alpha-l-1}+\left|\nabla \phi^{\frac{\alpha-l-1}{2}}\right|^2\right)\right)^{\frac{1}{2}}\left(\int_{\Omega}\left(\frac{1}{\phi}\right)^{\frac{3(l+1)}{4}}\right)^{\frac{2}{3}}dt\leq \\
& \notag\frac{C_1}{2}(l+1)(l+2)\int_0^T\left(\frac{1}{\phi}^{l+1},1\right)dt+\frac{C^2}{2C_1}(l+1)(l+2)\int_0^T\left(\int_{\Omega}\left(\frac{1}{\phi}\right)^{\frac{3(l+1)}{4}}\right)^{\frac{4}{3}}\,dt+\\
&\frac{a_*(l+1)(l+2)}{2(\alpha-l-1)^2}\int_0^T\int_{\Omega}\left|\nabla \phi^{\frac{\alpha-l-1}{2}}\right|^2\,d\mathbf{x}\,dt\\
&+\frac{C^2(l+1)(l+2)(\alpha-l-1)^2}{2a_*}\int_0^T\left(\int_{\Omega}\left(\frac{1}{\phi}\right)^{\frac{3(l+1)}{4}}\right)^{\frac{4}{3}}\,dt.
\end{align*}
Finally we obtain
\begin{align}
\label{eqn:46}
&\left(\frac{1}{\phi}^{l+1},1\right)+\frac{a_*(l+1)(l+2)}{(\alpha-l-1)^2}\int_0^T\int_{\Omega}\left|\nabla \phi^{\frac{\alpha-l-1}{2}}\right|^2\,d\mathbf{x}\,dt+\\
& \notag \frac{4(a_*-h)(l+1)(l+2)}{(\alpha-l+1)^2}\int_0^T\int_{\Omega}\left|\nabla \phi^{\frac{\alpha-l+1}{2}}\right|^2\,d\mathbf{x}\,dt\leq \left(\frac{1}{\phi_0}^{l+1},1\right)+\\
& \notag C(l+1)(l+2)(\alpha-l-1)^2\int_0^T\left(\int_{\Omega}\left(\frac{1}{\phi}\right)^{\frac{3(l+1)}{4}}\right)^{\frac{4}{3}}\,dt+\\
& \notag \frac{C(l+1)(l+2)}{(\alpha-l)^2}\int_0^T\int_{\Omega}\left|\nabla \phi^{\frac{\alpha-l}{2}}\right|^2\,d\mathbf{x}\,dt\leq\\
& \notag \left(\frac{1}{\phi_0}^{l+1},1\right)+C(l+1)^4\int_0^T\left(\int_{\Omega}\left(\frac{1}{\phi}\right)^{\frac{3(l+1)}{4}}\right)^{\frac{4}{3}}\,dt+\frac{Cl(l+1)}{(\alpha-l)^2}\int_0^T\int_{\Omega}\left|\nabla \phi^{\frac{\alpha-l}{2}}\right|^2\,d\mathbf{x}\,dt.
\end{align}
Iterating \eqref{eqn:46}, and considering that $\frac{1}{\phi_0},\frac{1}{\phi}>1$ for a. e. $(\mathbf{x},t)\in \Omega \times [0,T]$, we get
\begin{align}
\label{eqn:46bis}
&\left(\frac{1}{\phi}^{l+1},1\right)+\frac{a_*(l+1)(l+2)}{(\alpha-l-1)^2}\int_0^T\int_{\Omega}\left|\nabla \phi^{\frac{\alpha-l-1}{2}}\right|^2\,d\mathbf{x}\,dt+\\
& \notag \frac{4(a_*-h)(l+1)(l+2)}{(\alpha-l+1)^2}\int_0^T\int_{\Omega}\left|\nabla \phi^{\frac{\alpha-l+1}{2}}\right|^2\,d\mathbf{x}\,dt\leq \left(\frac{1}{\phi_0}^{l+1},1\right)+\\
& \notag C(l+1)^4\int_0^T\left(\int_{\Omega}\left(\frac{1}{\phi}\right)^{\frac{3(l+1)}{4}}\right)^{\frac{4}{3}}\,dt+
\frac{C_1l(l+1)}{(\alpha-l)^2}\int_0^T\int_{\Omega}\left|\nabla \phi^{\frac{\alpha-l}{2}}\right|^2\,d\mathbf{x}\,dt\leq \\
& \notag \left(\frac{1}{\phi_0}^{l+1},1\right)+C(l+1)^4\int_0^T\left(\int_{\Omega}\left(\frac{1}{\phi}\right)^{\frac{3(l+1)}{4}}\right)^{\frac{4}{3}}\,dt +\frac{C_1}{a_*}\left(\frac{1}{\phi_0}^{l},1\right)+\\
& \notag \frac{C_1C}{a_*}l^4\int_0^T\left(\int_{\Omega}\left(\frac{1}{\phi}\right)^{\frac{3l}{4}}\right)^{\frac{4}{3}}\,dt+\frac{C_1}{a_*}\left(\frac{1}{\phi_0}^{l-1},1\right)+\frac{C_1C}{a_*}(l-1)^4\int_0^T\left(\int_{\Omega}\left(\frac{1}{\phi}\right)^{\frac{3(l-1)}{4}}\right)^{\frac{4}{3}}\,dt+\\
& \notag+ \dots +\frac{C_1}{a_*}\left(\frac{1}{\phi_0},1\right)+\frac{C_1C}{a_*}\int_0^T\left(\int_{\Omega}\left(\frac{1}{\phi}\right)^{\frac{3}{4}}\right)^{\frac{4}{3}}\,dt C(l+1)\left(\frac{1}{\phi_0}^{l+1},1\right)+\\
& \notag \leq C(l+1)^5\int_0^T\left(\int_{\Omega}\left(\frac{1}{\phi}\right)^{\frac{3(l+1)}{4}}\right)^{\frac{4}{3}}\,dt\leq C(l+1)^5\max\left \{R_0^{l+1},\sup_{[0,T]}\left[\int_{\Omega}\left(\frac{1}{\phi}\right)^{\frac{3(l+1)}{4}}\,d\mathbf{x}\right]^{\frac{4}{3}}\right \},
\end{align}
where $R_0:=\left(\frac{1}{\phi_0},1\right)$.
\newline
Introducing the sequence $(l_k)_{k\geq 0}$ of real numbers defined by
\[
l_0=1, \; l_{k+1}=\frac{4}{3}l_k, \; k\in \mathbb{N},
\]
with $l+1=l_{k+1}$, and defining moreover the sequence $(\gamma_k)_{k\geq 0}$ of positive real numbers as
\[
\gamma_{k+1}=\sup_{[0,T]}\left[\int_{\Omega}\left(\frac{1}{\phi}\right)^{l_{k+1}}\,d\mathbf{x}\right],
\]
we can rewrite \eqref{eqn:46bis} as
\begin{equation}
\label{eqn:47}
\gamma_{k+1}\leq C l_{k+1}^5\max \left\{R_0^{l_{k+1}},\gamma_k^{\frac{4}{3}}\right\}.
\end{equation}
Thanks to Lemma $A.1$ in \cite{laurencot}, \eqref{eqn:47} implies that $\left(\gamma_k^{1/l_k}\right)_{k\geq 0}$ is bounded, which implies that  there exists a positive constant $\delta$ such that $\delta \leq \phi(\mathbf{x},t)$ for a. e. $(\mathbf{x},t)\in \Omega \times [0,T]$.

Once that the strict separation $\delta \leq \phi(\mathbf{x},t)$ is satisfied, the property $\phi(\mathbf{x},t)\leq 1- \delta$ in \eqref{eqn:39} is satisfied as a consequence of Lemma $8.3$ in \cite{frigeri2}. Indeed, given $\phi$ a weak solution of \eqref{eqn:38}, since $\phi\geq\delta$, if condition \eqref{eqn:7} is satisfied, i.e. $a_*>h+2-3\sqrt[3]{h}$, we have that
\[
b(\phi)\left(\psi^{''}(\phi)+a\right)=\phi^{\alpha}(1-\phi)^2\left(\frac{1-\phi*}{(1-\phi)^2}-2\phi-(1-\phi*)+a\right)>0.
\] 
Hence, similar arguments as in Lemma $8.3$ in \cite{frigeri2} give an $L^{\infty}(\Omega \times (0,T))$ bound for the function $-\log(1-\phi)$, which implies that there exists a $\delta>0$ such that $\phi(\mathbf{x},t)\leq 1- \delta$ for a. e. $(\mathbf{x},t)\in \Omega \times [0,T]$. The calculations in Lemma $8.3$ in \cite{frigeri2} can be straightforwardly adapted to our case by introducing the $a\frac{\phi^2}{2}$ term in the definition of the convex potential $F(\phi)$ and considering that $|\phi \nabla a - \nabla J \star \phi|\leq C$.
\end{pf}
\begin{rem}
\label{rem:rig}
We remark that this procedure may be made rigorous by considering truncated test functions $\xi=-\left(\max \{ \theta, \phi\}\right)^{-l-2}$ for some $\theta>0$ and taking the limit for $\theta \to 0$ in the derived estimates, splitting the domain of the integrals in the estimates between the sets $\{(\mathbf{x},t):\phi(\mathbf{x},t) < \theta\}$ and $\{(\mathbf{x},t):\phi(\mathbf{x},t) \geq \theta\}$, using moreover the results in remark \ref{rem:4} and the hypothesis that $\alpha\geq 2$.
Note that the argument works for all weak solutions in the class defined by Theorem \ref{thm:2}, since it does not require any approximation of the equations.
\end{rem}
\begin{rem}
\label{rem:h}
The condition \eqref{eqn:7} for the convexity of the potential $\psi(\phi)+\frac{a}{2}\phi^2$, i.e. $a_*>h+2-3\sqrt[3]{h}$, implies the condition $a_*>h$, which is needed to guarantee the strict separation property $\delta \leq \phi(\mathbf{x},t)$, only when $\bar{\phi}>\frac{19}{27}\sim 0.7$. When $\bar{\phi}<\frac{19}{27}\sim 0.7$, the condition $a_*>h$ is stricter than the condition for the convexity of the potential $\psi(\phi)+\frac{a}{2}\phi^2$. Indeed, since a-priori we have only that $\phi^{\alpha}(1-\phi)^2(\psi''(\phi)+a)\geq 0$, a stricter condition than the condition which guarantees that $(1-\phi)^2(\psi''(\phi)+a)>0$ may be needed to obtain the separation from zero of $\phi$. This fact is different from the standard cases analyzed in literature, where $b(\phi)(\psi''(\phi)+a)$ is taken to be greater than zero by constitutive assumptions. 
\end{rem}
\begin{rem}
\label{rem:mixed}
Due to the strict separation property \eqref{eqn:39}, it's easy to prove that the weak formulation \eqref{eqn:38} is equivalent to the dual mixed formulation \eqref{eqn:23bis}--\eqref{eqn:24}.
\end{rem}
\begin{thm}
\label{thm:4}
Let the hypothesis of Theorem \ref{thm:2} and Theorem \ref{thm:3} be satisfied. Then, the weak solution $\phi$ of \eqref{eqn:38} is unique. 
\end{thm}
\begin{pf}
Following \cite{giacomin1,frigeri1}, let us define the functions
\[
\Pi(s,\mathbf{x})=\int_0^s b(r)\left(\psi''(r)+a(\mathbf{x})\right)\,dr,
\]
with 
\begin{equation}
\label{eqn:48}
\Pi'(\cdot,\mathbf{x})\geq \alpha_0>0
\end{equation} 
for each $\mathbf{x}\in \Omega$, and
\[
\Theta(s)=\int_0^s b(r)\,dr.
\]
We can rewrite \eqref{eqn:38} as
\begin{equation}
\label{eqn:49}
<\phi_t,\xi>+\left(\nabla \Pi(\phi,\mathbf{x}), \nabla \xi \right)-\left(\Theta\nabla a, \nabla \xi \right)+ \left(b(\phi) \left(\phi \nabla a - \nabla J \star \phi \right), \nabla \xi \right)=0,
\end{equation}
$\forall \xi \in V, \; for \; \text{a.e.} \; t \in (0,T)$.
Considering two solutions $\phi_1$ and $\phi_2$ of \eqref{eqn:49}, with $\phi_1(0)=\phi_2(0)$, we take the difference between the corresponding equations \eqref{eqn:49}, set $\tilde{\phi}=\phi_1-\phi_2\in V_0'$, and take $\xi=\mathcal{G}\tilde{\phi}$, where $\mathcal{G}$ is the inverse Laplacian operator defined in \eqref{eqn:greencont}. We get
\begin{align*}
&\frac{1}{2}\frac{d}{dt}||\tilde{\phi}||_{-1}^2+\left(\Pi(\phi_1)-\Pi(\phi_2),\tilde{\phi}\right)=
\left(\left(\Theta(\phi_1)-\Theta(\phi_2)\right)\nabla a, \nabla \mathcal{G}\tilde{\phi} \right)- \\
& \left(\left(b(\phi_1)-b(\phi_2)\right) \left(\phi_1 \nabla a - \nabla J \star \phi_1 \right), \nabla  \mathcal{G}\tilde{\phi}\right)-\left(b(\phi_2) \left(\tilde{\phi} \nabla a - \nabla J \star \tilde{\phi} \right), \nabla  \mathcal{G}\tilde{\phi}\right).
\end{align*} 
Using \eqref{eqn:48}, the Cauchy--Schwarz and Young inequalities, the Lipschitz continuity of the functions $\Theta$ and $b$, the fact that $0<\phi_1,\phi_2<1$ and hypothesis $(H_2)$, we obtain
\begin{align}
\label{eqn:50}
&\frac{1}{2}\frac{d}{dt}||\tilde{\phi}||_{-1}^2+\alpha_0||\tilde{\phi}||^2\leq \frac{\alpha_0}{2}||\tilde{\phi}||^2+\\
& \notag C||\tilde{\phi}||_{-1}^2+C\left(||\nabla a||_{L^{\infty}(\Omega \times (0,T))}+||\nabla J||_{L^{\infty}(\Omega \times (0,T))}\right)\leq  \frac{\alpha_0}{2}||\tilde{\phi}||^2+C||\tilde{\phi}||_{-1}^2+C,
\end{align}
from which, applying a Gronwall argument and considering that $||\tilde{\phi}(0)||_{-1}=0$, we get that $||\tilde{\phi}||_{-1}=0$ for each $t\in [0,T]$, which implies the uniqueness of the solution.
\end{pf}

\section{Continuous Galerkin Finite Element approximation}
\label{numerics}
In this section we introduce the finite element and time discretization of \eqref{eqn:22}-\eqref{eqn:23}.
The Entropy estimate \eqref{eqn:boundPhi}, which guarantees the positivity of the continuous solution $\phi\geq 0$ a.e. $(\mathbf{x},t)\in \Omega_T$ (see \eqref{eqn:conv3d5}), and $\phi> 0$ a.e. $(\mathbf{x},t)\in \Omega_T$ if $\alpha \geq 2$, is not straightforwardly available at the discrete level. Also, the uniqueness result of Theorem \ref{thm:4} for the continuous solution in the case $\alpha\geq 2$, which is a consequence of the strict separation property \eqref{eqn:39}, is not straightforwardly available at the discrete level.
Following \cite{barrett,agosti1}, we impose the positivity property of the discrete solution as a constraint through a variational inequality, formulating a well posed discrete variational inequality defined on the proper domain of the degenerate elliptic operator in \eqref{eqn:22}$_1$. This approach ensures the uniqueness of the discrete solution together with its positivity. A proper approximation scheme for the primal formulation \eqref{eqn:38} of the model, based on positivity preserving schemes for porous-medium like equations, will be the subject of a forthcoming work. Here we concentrate on the formulation of a well posed and gradient stable finite element approximation scheme for the dual formulation \eqref{eqn:22} of the model. The forthcoming analysis is an adaptation of the analysis introduced in \cite{agosti1}, concerning a local degenerate Cahn--Hilliard equation with a singular single-well potential, to the present non-local case. The numerical analysis is performed for {$d\leq3$, while, for computational issues, the numerical simulations will be performed in the case $d=2$}.

We set $\Delta t = T/N$ for a $N \in \mathbb{N}$ and $t_{n}=n\Delta t$, $n=,...,N$.
Starting from a datum $\phi_{0}\in H$ and $\phi_{h}^{0}=\hat{P}^{h}\phi_{0}$, with $0\leq \phi_{h}^{0}<1$ if $0\leq \alpha < 2$ and $0< \phi_{h}^{0}<1$ if $\alpha \geq 2$, we consider the following fully discretized problem:
\\

\textbf{Problem $\mathbf{P}^{h}$.}
\label{pbm:ph0}
For $n=1,\dots,N$, given $\phi_{h}^{n-1}\in K^{h}$, find $(\phi_{h}^{n},\mu_{h}^{n})\in K^{h}\times S^{h}$ such that for all $(\chi,\xi)\in S^h\times K^h$,
\begin{equation}
\label{pbm:ph}
\begin{cases}
\displaystyle \biggl(\frac{\phi_{h}^n-\phi_{h}^{n-1}}{\Delta t},\chi \biggr)^h+(b(\phi_{h}^{n-1})\nabla \mu_{h}^n,\nabla \chi)=0,
\\
\displaystyle \left(\frac{1}{\epsilon}\psi'_{1}(\phi_{h}^n)+\epsilon(J\star 1)_h\phi_{h}^n,\xi -\phi_{h}^n\right)^h\geq \left(\mu_{h}^n-\frac{1}{\epsilon}\psi'_{2}(\phi_{h}^{n-1}),\xi -\phi_{h}^n\right)^h+\epsilon\left(J\star \phi_{h}^{n-1}, \xi -\phi_{h}^n\right)^{h^2},
\end{cases}
\end{equation}
where we used the notations \eqref{eqn:lump}-\eqref{eqn:lump3} for the lumped mass and convolution approximations.
Here, $b(\phi)$ has the particular form \eqref{eqn:23a} (with $M=1$) and $\psi_1(\phi),\psi_2(\phi)$ are given in \eqref{eqn:27}.
\begin{rem}
\label{rmk:ineq}
The variational inequality in \eqref{pbm:ph}$_2$ guarantees the positivity of the solution $\phi_h^n$. In order to prove this, let's define $ W^n:=W(\phi_h^n,\mu_h^n) $ in such a way that
\[
(W^n,\xi)^h = \epsilon \left( (J\star 1)_h\phi_{h}^n-(J\star \phi_{h}^{n-1})_h, \xi \right)^h + \left(\frac{1}{\epsilon}\left(\psi_1^{\prime} (\phi_h^n) +  \psi_2^{\prime}(\phi_h^{n-1})\right) - \mu_h^n, \xi\right)^h,
\]
for all $ \xi\in S_h $. Choosing $ \xi=0 $ and $ \xi=2\phi_h^n $ in \eqref{pbm:ph}$_2$ we obtain that
\begin{equation}
\label{ineq1}
(W^n,\phi_h^n)^h = \sum_{j\in \mathcal{I}}(1,\chi_j)W^n(\mathbf{x}_j)\phi_h^n(\mathbf{x}_j) = 0.
\end{equation}
Thanks to \eqref{ineq1}, we observe, as a consequence of the relation  
\[
\left( W^n, \xi-\phi_h^n \right)^h \geq 0 \qquad \forall \xi\in S_h^+,
\]
that
\begin{equation}
\label{ineq2}
(W^n,\xi)^h \geq 0 \qquad \forall \xi \in S_h^+.
\end{equation}
Thus, since \eqref{ineq2} is valid for any $ \xi\in S_h^+ $, it is valid in particular for $ \xi = \chi_j $, meaning that
\begin{equation}
\label{ineq3}
(W^n,\chi_j)^h = W^n(\mathbf{x}_j)\chi_j(\mathbf{x}_j) = W^n(\mathbf{x}_j) \geq 0.
\end{equation}
Collecting \eqref{ineq1} and \eqref{ineq3} we finally deduce that, for all $ j\in \mathcal{I} $, either $ \phi_h^n(\mathbf{x}_j)=0 $ or $ \phi_h^n(\mathbf{x}_j)>0 $ and 
\begin{equation}
\label{eqn:eqmu}
\epsilon \left( (J\star 1)_h\phi_{h}^n-(J\star \phi_{h}^{n-1})_h, \chi_j \right)^h + \left(\frac{1}{\epsilon}\left(\psi_1^{\prime} (\phi_h^n) +  \psi_2^{\prime}(\phi_h^{n-1})\right) - \mu_h^n, \chi_j\right)^h = 0.
\end{equation}
Therefore, the variational inequality in \eqref{pbm:ph}$_2$ projects the solution onto the set $\phi_h^n\geq 0 $, with $ \mu_h^n $ uniquely determined if $ \phi_h^n >0 $. 
\end{rem}
Defining the discrete energy functional $F_1:S^h\rightarrow \mathbb{R}^+$ as
\begin{equation}
\label{eqn:enfuncconv}
F_1[\phi_{h}^n]=\int_{\Omega}\left\{\frac{\epsilon}{2}(J\star 1)_h(\phi_{h}^n)^2+
\frac{1}{\epsilon}\psi_1^{\prime} (\phi_h^n)+\chi_{\mathbb{R}^+}(\phi_{h}^n)\right\}dx,
\end{equation}
where $\chi_{\mathbb{R}^+}(\cdot)$ is the indicator function of the closed and convex set $\mathbb{R}^+$, and endowing the space $S^h$ with the lumped scalar product \eqref{eqn:lump}, we can rewrite \eqref{pbm:ph}$_2$ as
\begin{equation}
\label{eqn:phsubdiff}
\left(\mu_{h}^n-\frac{1}{\epsilon}\psi'_{2}(\phi_{h}^{n-1})+\epsilon(J\star \phi_{h}^{n-1})_h,\xi -\phi_{h}^n\right)^h+F_1[\phi_{h}^n]\leq F_1[\xi], \quad \forall \xi \in S^h,
\end{equation}
which is equivalent to
\begin{equation}
\label{eqn:phsubdiffbiss}
\mu_{h}^n-\frac{1}{\epsilon}\psi'_{2}(\phi_{h}^{n-1})+\epsilon(J\star \phi_{h}^{n-1})_h\in \partial F_1[\phi_{h}^n],
\end{equation}
where $\partial$ is the subdifferential of the convex and lower semi-continuous function $F_1$.
\begin{rem}
\label{rmk:01}
Given the assumption $0\leq \phi_{h}^{0}<1$, the term $(\psi'_{1}(\phi_{h}^n),\xi -\phi_{h}^n)^h$ in \eqref{pbm:ph}$_2$ is well defined, since we will show that $||\phi_{h}^{0}||_{L^{\infty}(\Omega)}<1$ implies that $||\phi_{h}^{n}||_{L^{\infty}(\Omega)}<1$ for all $n\geq 1$ (see Lemma \ref{lem:chepswhepsconv}).
\end{rem}
We now introduce the discrete Green operator of the degenerate elliptic term in \eqref{pbm:ph}$_1$, which will be used to express the chemical potential $\mu_{h}^{n}$ in terms of $\frac{\phi_{h}^n-\phi_{h}^{n-1}}{\Delta t}$ and to show the well posedness of Problem $P^{h}$. We follow the approach in \cite{barrett} to invert the degenerate elliptic form on a proper closed and convex subset of $S^h$.
\begin{rem}
    \label{rem:deggreen}
    As observed in \cite{barrett}, the lumping approximation in the mass scalar product in \eqref{pbm:ph}$_1$ is introduced to avoid the locking phenomenon, due to the degeneracy, that regions where $\phi_h^{n-1}\equiv 0$ on the support of a basis functions $\chi_j$, for some $j\in \mathcal{I}$, remain fixed regions of zero values for the variable $\phi_h^n$. Indeed, the lumping approximation requires that in these regions only the values of $\phi_h^n$ at the nodes $\mathbf{x}_j$ are fixed to zero, leaving them free to evolve on the other nodes in the regions. Therefore, the discrete solution has the property of a moving support with finite speed, proportional to $h/\Delta t$. We highlight the fact that in the present context, the lumping approximations and integration rules \eqref{eqn:lump}-\eqref{eqn:lump3} make the computations of the convolution terms fast and easily parallelizable, as it will be discussed in Section \ref{numerics}. 
\end{rem}
We subdivide the partition $\mathcal{T}^h$ of $\Omega$ into elements on which $\phi_{h}^{n-1}\equiv 0$ and elements on which $\phi_{h}^{n-1}\neq 0$.
Given $q^h\in K^h$ with $\dashint q^h\in (0,1)$, where $\dashint q^h=\frac{1}{|\Omega|}(q^h,1)$, we define the set of passive nodes $J_{0}(q^h)\subset \mathcal{I}$ by
\begin{equation}
\label{eqn:passnod}
j\in J_{0}(q^h)  \Leftrightarrow \hat{P}^h q^h(x_j)=0 \Leftrightarrow (q^h,\chi_j)=0.
\end{equation}
The nodes in the set $J_{+}(q^h)=\mathcal{I}\setminus J_{0}(q^h)$ are called active nodes; these nodes can be partitioned into mutually disjoint and maximally connected subsets $I_{m}(q^h)$ such that $J_{+}(q^h)\equiv \bigcup_{m=1}^{M}I_{m}(q^h)$.
Defining
\[
\Sigma_{m}(q^h)=\sum_{j\in I_{m}(q^h)}\chi_{j},
\]
we note that
\begin{equation}
\label{eqn:sigma1}
\Sigma_{m}(q^h)\equiv 1 \quad \text{on each element on which} \; q^h \neq 0,
\end{equation}
since all the vertices of this elements belong to $I_{m}(q^h)$. 
There are also elements on which $q^h \equiv 0$ and $\Sigma_{m}(q^h)\equiv 1$. 
Hence, on each element $K\in \mathcal{T}^h$, we have that $q^h \equiv 0$ or $\Sigma_{m}(q^h)\equiv 1$ for some $m$, except for those elements on which both $q^h \equiv 0$ and $\Sigma_{m}(q^h)\equiv 1$.
Let us define the following set:
\[
\Omega_{m}(q^h)=\left \{\bigcup_{K\in \mathcal{T}^h} \bar{K}:\Sigma_{m}(q^h)(x)=1 \; \forall x \in K\right \},
\]
i.e. the union of the maximally connected elements on which $q^h\neq 0$, or $q^h\equiv 0$ and the indexes of the vertices of the elements belong to $I_{m}(q^h)$ for a given $m$.
Finally, we introduce the space
\begin{equation}
\label{eqn:setvh}
V^h(q^h)=\{v^h\in S^h:v^h(x_j)=0 \, \forall j \in J_{0}(q^h) \; \text{and} \; (v^h,\Sigma_{m}(q^h))^h=0, \, m=1,\dots , M\},
\end{equation}
that consists of all $v^h\in S^h$ which are orthogonal (with respect to the lumped discrete scalar product \eqref{eqn:lump}) to $\chi_j$, for $j \in J_{0}(q^h)$, (see \eqref{eqn:passnod}), and which have zero average (again with respect to the scalar product \eqref{eqn:lump}) on each element which does not contain any passive node.

We recall from \cite{barrett} that any $v^h \in S^h$ can be written as
\begin{equation}
\label{eqn:vhexpansion}
v^h\equiv \bar{v}^h+\sum_{j\in J_{0}(q^h)}v^h(x_j)\chi_j +\sum _{m=1}^{M}\biggl[\dashint_{\Omega_{m}(q^h)}v^h\biggr]\Sigma_{m}(q^h),
\end{equation}
where $\bar{v}^h$ is the $\hat{P}^h$ projection of $v^h$ onto $V^h(q^h)$, and
\begin{equation}
\label{eqn:average}
\dashint_{\Omega_{m}(q^h)}v^h:=\frac{(v^h,\Sigma_{m}(q^h))^h}{(1,\Sigma_{m}(q^h))}.
\end{equation}

We can now define, for all $q^h\in K^h$ with $q^h<1$, the discrete anisotropic Green's operator $\hat{G}_{q^h}^h:V^h(q^h)\rightarrow V^h(q^h)$ as
\begin{equation}
\label{eqn:discgreendeg}
(b(q^h)\nabla \hat{G}_{q^h}^h v^h, \nabla \chi)=(v^h,\chi)^h \quad \forall \chi \in S^h.
\end{equation}
The well posedness of $\hat{G}_{q^h}^h$ is shown in \cite{barrett}.

We now introduce a regularized version of the discrete problem \eqref{pbm:ph} , in which the singularity of the potential is regularized. We will show the well posedness and gradient stability of the regularized scheme, and we will study its limit problem as the regularization parameter tends to zero.

\subsection{Regularized problem}
In order to show the well posedness of Problem $P^{h}$, we introduce the following regularized version of \eqref{pbm:ph}:

\textbf{Problem $\mathbf{P}_{\lambda}^{h}$.}
\label{pbm:phlamb}
For $n=1,\dots,N$, given $\lambda>0$ a regularization parameter and $\phi_{h}^{n-1}\in K^{h}$, with $0\leq \phi_{h}^{n-1}<1$, find $(\phi_{h,\lambda}^{n},\mu_{h,\lambda}^{n})\in K^{h}\times S^{h}$ such that for all $(\chi,\xi)\in S^h\times K^h$,
\begin{equation}
\label{eqn:phlamb}
\begin{cases}
\displaystyle \biggl(\frac{\phi_{h,\lambda}^n-\phi_{h}^{n-1}}{\Delta t},\chi \biggr)^h+(b(\phi_{h}^{n-1})\nabla \mu_{h,\lambda}^n,\nabla \chi)=0,
\\ \\
\displaystyle \left(\frac{1}{\epsilon}\psi'_{1,\lambda}(\phi_{h,\lambda}^n)+\epsilon(J\star 1)_h\phi_{h,\lambda}^n,\xi -\phi_{h,\lambda}^n\right)^h\geq \left(\mu_{h,\lambda}^n-\frac{1}{\epsilon}\bar{\psi}'_{2}(\phi_{h}^{n-1}),\xi -\phi_{h,\lambda}^n\right)^h+\\
\epsilon\left(J\star \phi_{h}^{n-1}, \xi -\phi_{h,\lambda}^n\right)^{h^2},
\end{cases}
\end{equation}
where we used \eqref{eqn:29} and \eqref{eqn:31}.
We now prove the well posedness of Problem $P_{\lambda}^{h}$.

\begin{lem}
\label{lem:discex}
Assume that $J(\cdot)\geq 0$. Then, there exists a solution $(\phi_{h,\lambda}^n,\mu_{h,\lambda}^n)$ to Problem $P_{\lambda}^{h}$. Moreover, if $J(\cdot)> 0$ the solution $\{\phi_{h,\lambda}^n\}_{n=1}^{N}$ is unique, and $\mu_{h,\lambda}^n$ is unique on $\Omega_{m}(\phi_{h}^{n-1})$, for $m=1, \dots, M$ and $n=1, \dots, N$.
\end{lem}
\begin{pf}
From the first equation in \eqref{eqn:phlamb} and from \eqref{eqn:discgreendeg} it follows that, given $\phi_{h}^{n-1}\in K^h$, $\phi_{h}^{n-1}<1$, we search for $\phi_{h,\lambda}^{n}\in K^h(\phi_{h}^{n-1})$, where
\begin{equation}
\label{eqn:khc}
K^h(\phi_{h}^{n-1})=\{\chi \in K^h \, : \, \chi-\phi_{h}^{n-1}\in V^h(\phi_{h}^{n-1})\}.
\end{equation}
Moreover, a solution $\mu_{h,\lambda}^n\in S^h$ can be expressed in terms of $\phi_{h,\lambda}^{n}-\phi_{h}^{n-1}$ through the discrete anisotropic Green operator \eqref{eqn:discgreendeg}, recalling \eqref{eqn:vhexpansion}, as
\begin{equation}
\label{eqn:wneps}
\mu_{h,\lambda}^n= -\hat{\mathcal{G}}_{\phi_{h}^{n-1}}^h
\biggl[\frac{\phi_{h,\lambda}^{n}-\phi_{h}^{n-1}}{\Delta t}\biggr]+
\sum_{j\in J_{0}(\phi_{h}^{n-1})}\nu_{j,\lambda}^n\chi_j +\sum _{m=1}^{M}\gamma_{m,\lambda}^n\Sigma_{m}(\phi_{h}^{n-1}),
\end{equation}
where $\{\nu_{j,\lambda}^n\}_{j\in J_{0}(\phi_{h}^{n-1})}$ and $\{\gamma_{m,\lambda}^n\}_{m=1}^{M}$ are constants which express the values of $\mu_{h,\lambda}^n$ on the passive nodes and its \textit{average} value on $\Omega_{m}(\phi_{h}^{n-1})$, respectively.
Hence, Problem $P_{\lambda}^{h}$ can be restated as follows: given $\phi_{h}^{n-1}\in K^{h}$, with $\phi_{h}^{n-1}<1$, find $\phi_{h,\lambda}^{n}\in K^{h}(\phi_{h}^{n-1})$ and constant Lagrange multipliers $\{\nu_{j,\lambda}^n\}_{j\in J_{0}(\phi_{h}^{n-1})}$ and $\{\gamma_{m,\lambda}^n\}_{m=1}^{M}$ such that, for all $\chi \in K^h$,
\begin{multline}
\label{eqn:phepsvarin}
\displaystyle \left(\hat{\mathcal{G}}_{\phi_{h}^{n-1}}^h
\biggl[\frac{\phi_{h,\lambda}^{n}-\phi_{h}^{n-1}}{\Delta t}\biggr]+\epsilon (J\star 1)_h\phi_{h,\lambda}^n+\frac{1}{\epsilon}\psi'_{1,\lambda}(\phi_{h,\lambda}^n),\chi -\phi_{h,\lambda}^n\right)^h\\
\displaystyle \geq \biggl(\sum_{j\in J_{0}(\phi_{h}^{n-1})}\nu_{j,\lambda}^n\chi_j +\sum _{m=1}^{M}\gamma_{m,\lambda}^n\Sigma_{m}(\phi_{h}^{n-1})+\epsilon (J\star \phi_{h}^{n-1})_h-\frac{1}{\epsilon}\bar{\psi}'_{2}(\phi_{h}^{n-1}),\chi -\phi_{h,\lambda}^n\biggr)^h.
\end{multline}
We note that \eqref{eqn:phepsvarin} represents, together with $\phi_{h,\lambda}^{n}\in K^{h}(\phi_{h}^{n-1})$, the Karush-Kuhn-Tucker optimality condition of the minimization problem
\begin{multline}
\label{eqn:lagrangian}
\displaystyle \inf_{v_{h,\lambda}\in S^h}\sup_{\nu_{j,\lambda},\gamma_{m,\lambda},\eta_{\lambda} \geq 0}\biggl\{
\epsilon\left((J\star 1)_hv_{h,\lambda}^n,v_{h,\lambda}^n\right)^h+\frac{1}{\Delta t}||[b(\phi_{h}^{n-1})]^{1/2}\nabla \hat{\mathcal{G}}_{\phi_{h}^{n-1}}^h(v_{h,\lambda}-\phi_{h}^{n-1}) ||^{2}\\
\displaystyle +\frac{2}{\epsilon}(\psi_{1,\lambda}(v_{h,\lambda}),1)^h -2\epsilon \left((J\star \phi_{h}^{n-1})_h,v_{h,\lambda}^n\right)^h+\frac{2}{\epsilon}(\bar{\psi}'_{2}(\phi_{h}^{n-1}),v_{h,\lambda})^h\\
\displaystyle
-\sum_{j\in J_{0}(\phi_{h}^{n-1})}\nu_{j,\lambda}(\chi_j,v_{h,\lambda})^h
-\sum _{m=1}^{M}\gamma_{m,\lambda}(\Sigma_{m}(\phi_{h}^{n-1}),v_{h,\lambda})^h
-(\eta_{\lambda},v_{h,\lambda})^h
\biggr\},
\end{multline}
being $\eta_{\lambda}\in K^h$ the Lagrange multiplier of the inequality constraint.
Noting the convexity and the coercivity of $\psi_{1,\lambda}(\cdot)$ (see equation \eqref{eqn:29}), the facts that $(J\star 1)_h\geq 0$ (cfr. Remark \ref{rem:astarpos})
and that $\phi_h^{n-1}\in K^h$, the primal form associated to the Lagrangian \eqref{eqn:lagrangian} is a convex, proper, lower semi continuous and coercive function from the closed and convex set $K^h(\phi_{h}^{n-1})$ to $\mathbb{R}$, and the primal problem is stable. 
Hence, from the Kuhn-Tucker Theorem (see, e.g., theorem $5.1$ in \cite{temam}), we infer the existence of $\phi_{h,\lambda}^n\in K^h(\phi_{h}^{n-1})$, solution to the primal problem, and Lagrange multipliers $\{\nu_{j,\lambda}^n\}_{j\in J_{0}(\phi_{h}^{n-1})}$, $\{\gamma_{m,\lambda}^n\}_{m=1}^{M}$ and $\eta_{\lambda}(x_i)\in -\partial \chi _{\mathbb{R}^+}(\phi_{h,\lambda}^n(x_i))$, for each $i\in \mathcal{I}$ and each $n$. Therefore, from \eqref{eqn:wneps} we have the existence of a solution $(\phi_{h,\lambda}^n,\mu_{h,\lambda}^n)_{n=1}^{N}$ to Problem $P_{\lambda}^{h}$.

Let us now prove uniqueness. If, for fixed $n\geq 1$, \eqref{eqn:phepsvarin} has two solutions \\ $(\phi_{h,\lambda}^{n,i},\{\nu_{j,\lambda}^{n,i}\}_{j\in J_{0}(\phi_{h}^{n-1})},\{\gamma_{m,\lambda}^{n,i}\}_{m=1}^{M})$, $i=1,2$, choosing $\chi = \phi_{h,\lambda}^{n,2}$ in the inequality for $\phi_{h,\lambda}^{n,1}$ and $\chi = \phi_{h,\lambda}^{n,1}$ in the inequality for $\phi_{h,\lambda}^{n,2}$, taking the difference between the former inequality, multiplied by $-1$, and the latter inequality, setting $\phi_{h,\lambda}^{n,1}-\phi_{h,\lambda}^{n,2}=d_{h,\lambda}^{n}$, we have
\begin{align*}
&\epsilon \left((J\star 1)_hd_{h,\lambda}^{n},d_{h,\lambda}^{n}\right)^h+\frac{1}{\Delta t}||[b(\phi_{h}^{n-1})]^{1/2}\nabla \hat{\mathcal{G}}_{\phi_{h}^{n-1}}^hd_{h,\lambda}^{n} ||^{2}+
\frac{1}{\epsilon}(\psi_{1,\lambda}'(\phi_{h,\lambda}^{n,1})-\psi_{1,\lambda}'(\phi_{h,\lambda}^{n,2}),d_{h,\lambda}^{n})^h\\
&-\left(\sum_{j\in J_{0}(\phi_{h}^{n-1})}\nu_{j,\lambda}^{n,1}\chi_j -\sum_{j\in J_{0}(\phi_{h}^{n-1})}\nu_{j,\lambda}^{n,2}\chi_j +\sum _{m=1}^{M}\gamma_{m,\lambda}^{n,1}\Sigma_{m}(\phi_{h}^{n-1})-\sum _{m=1}^{M}\gamma_{m,\lambda}^{n,2}\Sigma_{m}(\phi_{h}^{n-1}),d_{h,\lambda}^{n}\right)^h\leq 0.
\end{align*}
Since $d_{h,\lambda}^{n}\in V^h(\phi_{h}^{n-1})$, recalling the definition \eqref{eqn:setvh}, we have that the second line of the previous inequality is identically equal to zero. Given also the monotonicity of $\psi_{1,\lambda}'(\cdot)$, we obtain that
\begin{equation}
\label{eqn:unique1}
\epsilon \left((J\star 1)_hd_{h,\lambda}^{n},d_{h,\lambda}^{n}\right)^h+\frac{1}{\Delta t}||[b(\phi_{h}^{n-1})]^{1/2}\nabla \hat{\mathcal{G}}_{\phi_{h}^{n-1}}^hd_{h,\lambda}^{n} ||^{2}\leq 0.
\end{equation}
We observe finally that $(J\star 1)_h> 0$, as a consequence of the assumption $J(\cdot)>0$ and of the Remark \ref{rem:astarpos}, and that the lumped scalar product \eqref{eqn:lump} induces a norm on $S^h$ (as a consequence of the fact that $\chi_j\in K^h$). Hence, we deduce from \eqref{eqn:unique1} the uniqueness of $\phi_{h,\lambda}^{n}$.
Choosing $\chi = \phi_{h,\lambda}^{n}\pm \delta \pi_h [\phi_{h,\lambda}^{n}\Sigma_{m}(\phi_{h}^{n-1})]$ in \eqref{eqn:phepsvarin}, for any $\delta \in (0,1)$ and $m=1, \dots, M$, gives an equation with terms depending on $ \gamma_{m,\lambda}^n$ and $\phi_{h,\lambda}^{n}$ and concentrated on the support of $ \Sigma_m(\phi_{h}^{n-1})$.
This yields uniqueness of the Lagrange multiplier $\gamma_{m,\lambda}^n$, as a consequence of the uniqueness of $\phi_{h,\lambda}^{n}$. Hence the uniqueness of $\mu_{h,\lambda}^{n}$ on $ \Omega_{m}(\phi_{h}^{n-1})$ follows from \eqref{eqn:wneps}. The proof is complete.
\end{pf}
\begin{rem}
    \label{rem:convdisc}
    Observing from \eqref{eqn:28} that 
    \[
    \min_{r\in \mathbb{R}^+}\psi_{1,\lambda}''(r)=1-\bar{\phi},
    \]
    the requirement that $J(\cdot)\geq 0$, and hence that $(J\star 1)_h\geq 0$, in Lemma \ref{lem:discex}, which implies the convexity of the primal form associated to \eqref{eqn:lagrangian}, could be relaxed to the requirement that $(J\star 1)_h\geq -(1-\bar{\phi})$. We also observe that, when $J(\cdot)\geq 0$ is not verified, we could split $J(\cdot)=J_+(\cdot)+J_-(\cdot)$, with $J_+(\cdot)\geq 0$ and $J_-(\cdot)\leq 0$, and consider in \eqref{pbm:ph} the splitting \\ $\epsilon (J_+\star 1)_h\phi_h^n+\epsilon (J_-\star 1)_h\phi_h^{n-1}$, as done e.g. in \cite{guan} for the approximation of a nonlocal Cahn--Hilliard equation with constant mobility and smooth potential. 
\end{rem}
We now obtain $\lambda$-independent stability bounds for the solution $(\phi_{h,\lambda}^n,\mu_{h,\lambda}^n)$ of System \eqref{eqn:phlamb}, employed to pass to the limit as $\lambda\rightarrow 0$ in system \eqref{eqn:phlamb} and thus to obtain an existence and stability result for the unregularized System \eqref{pbm:ph}. Before proceeding, we need to introduce the following subsets of $\Omega$:
\[
\Omega_{m,*}(\phi_{h}^{n-1}):=\left\{\bigcup_{K\in \Omega_{m}(\phi_{h}^{n-1})}\bar{K}:\phi_{h}^{n-1}\bigr|_{K}\not\equiv 0\right\},
\]
and, for any $q_h\in K^h$,
\[
\Omega_{m,+}(q_h,\phi_{h}^{n-1}):=\left\{\bigcup_{j\in I_{m}(\phi_{h}^{n-1})}\text{supp}(\chi_j):q_{h}(\mathbf{x}_j)> 0\right\}.
\]
The following convergence results holds.
\begin{lem}
\label{lem:chepswhepsconv}
Let $J(\cdot )\in C(\bar{\Omega})$, with $J(\cdot )\geq 0$ and $J(\mathbf{x})=J(-\mathbf{x})$ for all $\mathbf{x}\in \Omega$. Let moreover $0\leq \phi_h^{n-1}<1$. Then, for every sequence $\lambda_k \to 0$, there exist a subsequence $\lambda_k' \to 0$ and a $\phi_{h}^n\in K^h$, with $0\leq \phi_{h}^n<1$, such that, for $k \to \infty$,
\begin{equation}
\label{eqn:chepsconvlem}
\phi_{h,\lambda_k '}^n\rightarrow \phi_{h}^n.
\end{equation}
For every sequence $\lambda_k \to 0$, there exist a subsequence $\lambda_k' \to 0$ and a $\mu_{h}^n\in S^h$ such that, for $k \to \infty$,
\begin{equation}
\label{eqn:whepsconvlem}
\mu_{h,\lambda_k '}^n\rightarrow \mu_{h}^n \quad \text{on} \quad \Omega_{m,*}(\phi_{h}^{n-1})\cup \Omega_{m,+}(\phi_{h}^{n},\phi_{h}^{n-1}) \quad \text{and} \quad \nabla \mu_{h,\lambda_k '}^n\rightarrow \nabla \mu_{h}^n \quad \text{on} \; \Omega_{m,*}(\phi_{h}^{n-1}),
\end{equation}
for each $m=1,\dots,M$.
\end{lem}
\begin{pf}
We start by proving stability bounds for the regularized Problem \eqref{eqn:phlamb}.
Choosing $\chi=\mu_{h,\lambda}^{n}$ in the first equation of \eqref{eqn:phlamb} and $\xi=\phi_{h}^{n-1}$ in the second equation of \eqref{eqn:phlamb}, we get
\begin{align}
\label{eqn:stabtemp1}
& \frac{1}{\epsilon}(\psi'_{1,\lambda}(\phi_{h,\lambda}^{n})+\bar{\psi}'_{2}(\psi_{h}^{n-1}),\phi_{h,\lambda}^{n}-\phi_{h}^{n-1})^h+\epsilon \left((J\star 1)_h\phi_{h,\lambda}^{n},\phi_{h,\lambda}^{n}-\phi_{h}^{n-1}\right)^h\\
& \notag -\epsilon \left(J\star \phi_{h}^{n-1},\phi_{h,\lambda}^{n}-\phi_{h}^{n-1}\right)^{h^2}
+\Delta t||[b(\phi_{h}^{n-1})]^{1/2}\nabla \mu_{h,\lambda}^{n}||^{2}\leq 0.
\end{align}
Thanks to the symmetry of the kernel $J$, we have that
\begin{align*}
& \frac{1}{2}\left(J\star (\phi_{h,\lambda}^{n}-\phi_{h}^{n-1}),\phi_{h,\lambda}^{n}-\phi_{h}^{n-1}\right)^{h^2}=\frac{1}{2}\left(J\star \phi_{h,\lambda}^{n},\phi_{h,\lambda}^{n}\right)^{h^2}+\frac{1}{2}\left(J\star \phi_{h}^{n-1},\phi_{h}^{n-1}\right)^{h^2}\\
&-\left(J\star \phi_{h}^{n-1},\phi_{h,\lambda}^{n}\right)^{h^2}.
\end{align*}
Hence,
\begin{align*}
& \frac{1}{2}\left(J\star \phi_{h}^{n-1},\phi_{h}^{n-1}\right)^{h^2}-\frac{1}{2}\left(J\star \phi_{h,\lambda}^{n},\phi_{h,\lambda}^{n}\right)^{h^2}+\frac{1}{2}\left(J\star (\phi_{h,\lambda}^{n}-\phi_{h}^{n-1}),\phi_{h,\lambda}^{n}-\phi_{h}^{n-1}\right)^{h^2}\\
&=-\left(J\star \phi_{h}^{n-1},\phi_{h,\lambda}^{n}-\phi_{h}^{n-1}\right)^{h^2}.
\end{align*}
Using the latter relation in \eqref{eqn:stabtemp1}, using moreover the identity $2s(s-r)=s^2-r^2+(s-r)^2$, $\forall r,s \in \mathbb{R}$, and the convexity and the concavity properties of $\psi_{1,\lambda}(\cdot)$ and $\bar{\psi}_{2}(\cdot)$, it follows that
\begin{align}
\label{eqn:estimateeps}
&\notag \frac{1}{\epsilon}(\psi_{1,\lambda}(\phi_{h,\lambda}^{n})+\bar{\psi}_{2}(\phi_{h,\lambda}^{n}),1)^h+\frac{\epsilon}{2} \left((J\star 1)_h\phi_{h,\lambda}^{n},\phi_{h,\lambda}^{n}\right)^h+\frac{\epsilon}{2}\left((J\star 1)_h(\phi_{h,\lambda}^{n}-\phi_{h}^{n-1}),(\phi_{h,\lambda}^{n}-\phi_{h}^{n-1})\right)^h\\
& \notag -\frac{\epsilon}{2}\left(J\star \phi_{h,\lambda}^{n},\phi_{h,\lambda}^{n}\right)^{h^2} +\frac{\epsilon}{2}\left(J\star (\phi_{h,\lambda}^{n}-\phi_{h}^{n-1}),(\phi_{h,\lambda}^{n}-\phi_{h}^{n-1})\right)^{h^2}+\Delta t||[b(\phi_{h}^{n-1})]^{1/2}\nabla \mu_{h,\lambda}^{n}||^{2}\\
& \leq \frac{1}{\epsilon}(\psi_{1,\lambda}(\phi_{h}^{n-1})+\bar{\psi}_{2}(\phi_{h}^{n-1}),1)^h+\frac{\epsilon}{2} \left((J\star 1)_h\phi_{h}^{n-1},\phi_{h}^{n-1}\right)^h-\frac{\epsilon}{2}\left(J\star \phi_{h}^{n-1},\phi_{h}^{n-1}\right)^{h^2}.
\end{align}
Using \eqref{eqn:jhcomp} in \eqref{eqn:estimateeps} and the hypothesis on $J$ and on $\phi_h^{n-1}$, we finally obtain that
\begin{align}
\label{eqn:estimateeps2}
&\notag \frac{1}{\epsilon}(\psi_{1,\lambda}(\phi_{h,\lambda}^{n})+\bar{\psi}_{2}(\phi_{h,\lambda}^{n}),1)^h+\underbrace{\frac{\epsilon}{2} \sum_{i,j\in \mathcal{I}}J(\mathbf{x}_i-\mathbf{x}_j)\left(\phi_{h,\lambda}^{n}(\mathbf{x}_i)-\phi_{h,\lambda}^{n}(\mathbf{x}_j)\right)^2(1,\chi_i)(1,\chi_j)}_{\geq 0}\\
& \notag +\underbrace{\frac{\epsilon}{2} \sum_{i,j\in \mathcal{I}}J(\mathbf{x}_i-\mathbf{x}_j)\left(\phi_{h,\lambda}^{n}(\mathbf{x}_i)-\phi_{h}^{n-1}(\mathbf{x}_i)+\phi_{h,\lambda}^{n}(\mathbf{x}_j)-\phi_{h}^{n-1}(\mathbf{x}_j)\right)^2(1,\chi_i)(1,\chi_j)}_{\geq 0}\\
& +\Delta t||[b(\phi_{h}^{n-1})]^{1/2}\nabla \mu_{h,\lambda}^{n}||^{2}\leq C.
\end{align}
From \eqref{eqn:29} and \eqref{eqn:31}, observing \eqref{eqn:29+31} and the fact that $\psi_{1,\lambda}(r)+\bar{\psi}_{2}(r)$ is bounded from below by a negative constant when $r\in [0,1]$, we have the property, for $\lambda$ sufficiently small, that there exist $C_1,C_2>0$ independent on $\lambda$ such that
\begin{equation}
    \label{eqn:estimateeps3}
    \psi_{1,\lambda}(r)+\bar{\psi}_{2}(r)\geq C_1r^2-C_2, \quad \forall r\in \mathbb{R}^+.
\end{equation}
Hence, we conclude from \eqref{eqn:estimateeps3} and \eqref{eqn:estimateeps2} that
\begin{equation}
\label{eqn:estimateeps3bis}
(\phi_{h,\lambda}^{n},\phi_{h,\lambda}^{n})^h\leq C,
\end{equation}
and since the lumped scalar product \eqref{eqn:lump} induces a norm on $S_h$, from the Bolzano--Weierstrass Theorem it follows that there exists a subsequence $\{\phi_{h,\lambda_k '}^n\}$ and a $\phi_{h}^n\in K^h$ such that \eqref{eqn:chepsconvlem} holds. 

We now prove that $\phi_{h,\lambda}^{n}<1-\lambda$ for sufficiently small $\lambda$. We start by observing from \eqref{eqn:29} that
\begin{equation}
    \label{eqn:estimateeps4}
    \psi_{\lambda}(r)+\frac{1}{3}\geq \frac{1-\bar{\phi}}{2\lambda^2}([r-1]_+)^2, \quad \forall r\in \mathbb{R}, \forall \lambda \leq \lambda_0,
\end{equation}
with $\lambda_0$ sufficiently small. Hence, from \eqref{eqn:estimateeps2} we deduce that
\begin{equation}
\label{eqn:chepsmin1}
([\phi_{h,\lambda}^{n}-1]_{+}^2,1)^{h}\leq C\lambda^2.
\end{equation}
From \eqref{eqn:chepsmin1}, \eqref{eqn:interp1} and \eqref{eqn:interp5} it also follows that
\begin{equation}
\label{eqn:chepsmin1inf}
||[\phi_{h,\lambda}^{n}-1]_{+}||_{L^{\infty}(\Omega)}\leq Ch^{-d/2}\lambda.
\end{equation}
Let us suppose that there exist nodes $x_l$, $l\in \mathcal{I}$, such that $\phi_{h,\lambda}^{n}(x_l)\geq 1-\lambda$. From \eqref{eqn:29}, \eqref{eqn:31}, \eqref{eqn:estimateeps2}, \eqref{eqn:estimateeps3bis} and \eqref{eqn:chepsmin1} we have that
\[
-(1-\bar{\phi})\log (\lambda)\left(\sum_{l\in \mathcal{I}}\chi_l,1\right)+(\psi_{1,\lambda}(\phi_{h,\lambda}^{n}),1)^h\biggl|_{\phi_{h,\lambda}^{n}< 1-\lambda}\leq Ch^{-d/2}-(\bar{\psi}_{2}(\phi_{h,\lambda}^{n}),1)^h\leq Ch^{-d/2}+C.
\]
Calling $\Omega_{\lambda}$ the support of the basis functions corresponding to nodes on which $\phi_{h,\lambda}^{n}\geq 1-\lambda$, using \eqref{eqn:interp1} and the fact that $\psi_{1}(\phi_{h,\lambda}^{n})\geq 0$ for all $0\leq \phi_{h,\lambda}^{n}\geq 1-\lambda$, we have that
\begin{equation}
\label{eqn:chepsminunfbound}
||-(1-\bar{\phi})\log \lambda||_{L^{\infty}(\Omega_{\lambda})}\leq C\,h^{-3d/2},
\end{equation}
where the bound is independent on $\lambda$, which is absurd. Therefore, we deduce that there exists a value $\lambda_0$ sufficiently small such that, for each $\lambda \leq \lambda_0$,
\begin{equation}
\label{eqn:unfboundpsi1}
\phi_{h,\lambda}^{n}< 1-\lambda.
\end{equation}
We next show \eqref{eqn:whepsconvlem}. 
We introduce the set 
\[
I_m^{*}(\phi_h^{n-1}) \subset I_m(\phi_h^{n-1})
\]
given by the nodes of $I_m(\phi_h^{n-1})$ which are vertices of elements in $\Omega_{m,*}(\phi_{h}^{n-1})$, and the quantities
\begin{align}
\label{eqn:averagestar}
	&\Sigma_m^*(\phi_{h}^{n-1}) := \sum_{j\in I_m^*(\phi_{h}^{n-1})}\chi_j \\
	& \notag \dashint_{\Omega_{m,*}(\phi_{h}^{n-1})}\mu_{h,\lambda}^n := \frac{(v^h,\Sigma_{m}^*(\phi_{h}^{n-1}))^h}{(1,\Sigma_{m}^*(\phi_{h}^{n-1}))^h}.
\end{align}
Using the Poincar\'e's inequality on $\Omega_{m,*}(\phi_{h}^{n-1})$, \eqref{eqn:interp5}, \eqref{eqn:estimateeps} and the fact that $\phi_h^{n-1}<1$ leads to
\begin{align}
\label{eqn:wepsnpoinc}
&\left(\left(\left[\left(I-\dashint _{\Omega_{m,*}(\phi_{h}^{n1})}\right)\mu_{h,\lambda}^n\right]\Sigma_m^*(\phi_{h}^{n-1})\right)^2,1\right)^{h}  \leq C\int _{\Omega_{m,*}(\phi_{h}^{n-1})}|\nabla \mu_{h,\lambda}^n|^2d\mathbf{x}\\
&\notag \leq C[b_{\text{min}}(\phi_{h}^{n-1})]^{-1}\int _{\Omega_{m,*}(\phi_{h}^{n-1})}b(\phi_{h}^{n-1})\nabla |\mu_{h,\lambda}^n|^2d\mathbf{x}\leq C((\Delta t)^{-1})[b_{min}(\phi_{h}^{n-1})]^{-1},
\end{align}
where 
\[b_{\text{min}}(\phi_{h}^{n-1}):=\min_{K\subset \Omega_{m,*}(\phi_{h}^{n-1})}\frac{1}{|K|}\int_{K}b(\phi_{h}^{n-1})dx.\]
We now bound $\dashint _{\Omega_{m,*}(\phi_{h}^{n-1})}\mu_{h,\lambda}^n$.
Let us take
\[
 K^h \ni \xi=\phi_{h,\lambda}^n+\Sigma_{m}^*(\phi_{h}^{n-1})
\]
in the second equation of system \eqref{eqn:phlamb}. We get
\begin{align*}
&(\mu_{h,\lambda}^n,\Sigma_{m}^*(\phi_{h}^{n-1}))^h\leq \left(\frac{1}{\epsilon}\psi'_{1,\lambda}(\phi_{h,\lambda}^n),\Sigma_{m}^*(\phi_{h}^{n-1})\right)^h+\left(\frac{1}{\epsilon}\bar{\psi}'_{2}(\phi_{h}^{n-1}),\Sigma_{m}^*(\phi_{h}^{n-1})\right)^h\\
&+\left(\epsilon(J\star 1)_h\phi_{h,\lambda}^n,\Sigma_{m}^*(\phi_{h}^{n-1})\right)^h
-\epsilon\left(J\star \phi_{h}^{n-1}, \Sigma_{m}^*(\phi_{h}^{n-1})\right)^{h^2}.
\end{align*}
Given \eqref{eqn:unfboundpsi1} and the hypothesis $0\leq \phi_h^{n-1}<1$, which give the boundedness of the terms $\psi'_{1,\lambda}(\phi_{h,\lambda}^n)$ and $\bar{\psi}'_{2}(\phi_{h}^{n-1})$, given moreover the boundedness of $\Sigma_{m}^*(c_{h}^{n-1})$, the hypothesis that $J(\cdot)\in C(\bar{\Omega})$ and \eqref{eqn:estimateeps3bis}, we obtain that
\begin{equation}
\label{eqn:approach21}
|(\mu_{h,\lambda}^n,\Sigma_{m}^*(\phi_{h}^{n-1}))^h|\leq C+C||\Sigma_{m}^*(\phi_{h}^{n-1})||_{L^{\infty}(\Omega)}+C(\phi_{h,\lambda}^n,\phi_{h,\lambda}^n)^h\leq C.
\end{equation}
Now, combining \eqref{eqn:wepsnpoinc} with \eqref{eqn:approach21}, recalling the definition \eqref{eqn:averagestar} and using the Poincar\'e's inequality, we obtain
\begin{equation}
\label{eqn:wepsnl2conv}
(\mu_{h,\lambda}^n\Sigma_{m}^*(\phi_{h}^{n-1}),\mu_{h,\lambda}^n\Sigma_{m}^*(\phi_{h}^{n-1}))^h\leq C+C((\Delta t)^{-1})[b_{min}(\phi_{h}^{n-1})]^{-1}.
\end{equation}
We also introduce the set 
\[
I_m^{+}(\phi_{h}^n,\phi_h^{n-1}) \subset I_m(\phi_h^{n-1})
\]
given by the nodes of $I_m(\phi_h^{n-1})$ belonging to $\Omega_{m,+}(\phi_{h}^{n},\phi_h^{n-1})$, and the quantity
\begin{align}
\label{eqn:averageplus}
	&\Sigma_m^+(\phi_{h}^{n},\phi_h^{n-1}) := \sum_{j\in I_m^+(\phi_{h}^{n},\phi_h^{n-1})}\chi_j.
\end{align}
We observe from the convergence property \eqref{eqn:chepsconvlem} that, given $j\in I_m^{+}(\phi_{h}^n,\phi_{h}^{n-1})$, $\phi_h^n(\mathbf{x}_j)>0$ implies that $\phi_{h,\lambda}^n(\mathbf{x}_j)>0$ for sufficiently small $\lambda$.
Combining \eqref{eqn:eqmu} for each $j\in I_m^{+}(\phi_{h}^n,\phi_{h}^{n-1})$ and testing the resulting equation by $\pi^h[\mu_{h,\lambda}^n\Sigma_m^+(\phi_{h}^{n},\phi_h^{n-1})]$, with similar calculations as those employed in \eqref{eqn:approach21} we obtain that 
\begin{align*}
&(\mu_{h,\lambda}^n,\mu_{h,\lambda}^n\Sigma_{m}^+(\phi_{h}^{n},\phi_h^{n-1}))^h\leq \\
&C+C||\Sigma_{m}^+(\phi_{h}^{n},\phi_h^{n-1})||_{L^{\infty}(\Omega)}+C(\phi_{h,\lambda}^n,\phi_{h,\lambda}^n)^h+\frac{1}{2}(\mu_{h,\lambda}^n,\mu_{h,\lambda}^n\Sigma_{m}^+(\phi_{h}^{n},\phi_h^{n-1}))^h,
\end{align*}
hence
\begin{equation}
    \label{eqn:wepsnl2convplus}
    (\mu_{h,\lambda}^n,\mu_{h,\lambda}^n\Sigma_{m}^+(\phi_{h}^{n},\phi_h^{n-1})^h\leq C.
\end{equation}

Finally, from \eqref{eqn:wepsnl2conv}, \eqref{eqn:wepsnl2convplus}, \eqref{eqn:estimateeps}, \eqref{eqn:interp5} and the Bolzano-Weierstrass Theorem it follows that there exist a subsequence $\{\mu_{h,\lambda_k '}^n\}$ and a $\mu_{h}^n\in S^h$ such that \eqref{eqn:whepsconvlem} holds. 

We finally prove that the limit point $\phi_h^n$ satisfies the property $\phi_h^n<1$. We observe from the fact that $0\leq \phi_{h,\lambda}^n<1-\lambda$ and \eqref{eqn:29} that $\psi_{1,\lambda}(\phi_{h,\lambda}^n)\geq 0$. Also, from the fact that $\chi_j\in K^h$ for all $j\in \mathcal{I}$, we have that $\pi_h[\psi_{1,\lambda}(\phi_{h,\lambda}^n)]\geq 0$. Hence, from the Fatou's Lemma and from \eqref{eqn:estimateeps2} we have that
\begin{equation}
\label{eqn:fatouh}
\int_{\Omega} \lim_{\lambda \to 0}\pi_h[\psi_{1,\lambda}(\phi_{h,\lambda}^n)]\leq \lim_{\lambda \to 0}\int_{\Omega}\pi_h[\psi_{1,\lambda}(\phi_{h,\lambda}^n)]=\lim_{\lambda \to 0}(\psi_{1,\lambda}(\phi_{h,\lambda}^n),1)^h\leq C.
\end{equation}
From the convergence property \eqref{eqn:chepsconvlem} and from \eqref{eqn:29} we have that
\[
\lim_{\lambda \to 0}\pi_h[\psi_{1,\lambda}(\phi_{h,\lambda}^n)]=
\begin{cases}
\pi_h[\psi_{1}(\phi_{h}^n)] \quad \text{if} \; \phi_h^{n}<1,\\
\infty \quad \text{elsewhere},
\end{cases}
\]
which implies, together with \eqref{eqn:fatouh}, that 
\begin{equation}
    \label{eqn:unifbound2}
    \phi_h^{n}<1.
\end{equation}
The proof is complete.
\end{pf}
\subsection{Well-posedness and gradient stability of Problem $\mathbf{P}^{h}$.}
We now prove the well-posedness and gradient stability of Problem $P^{h}$.
\begin{thm}
\label{th:wellposednessph}
Let the hypothesis of Lemma \ref{lem:chepswhepsconv} be satisfied.
Then there exists a solution $(\phi_{h}^{n},\mu_{h}^{n})$ to Problem \eqref{pbm:ph} for any $n=1,\dots,N$. The solution $\{\phi_{h}^{n}\}_{n=1}^{N}$ is unique, while the solution $\mu_{h}^{n}$ is unique on $\Omega_{m}(\phi_{h}^{n-1})$, for $m=1, \dots, M$ and $n=1, \dots, N$. Moreover, the solution satisfies the stability bound
\begin{equation}
    \label{eqn:stab1}
    \max_{n=1\to N}||\phi_h^n||^2+\Delta t\sum_{n=1}^N||(b(\phi_h^{n-1}))^{1/2}\nabla \mu_h^n||^2+\Delta t \sum_{n=1}^N(b_{\text{max}}^{n-1})^{-1}\left|\left|\hat{\mathcal{G}}^h\left(\frac{\phi_h^n-\phi_h^{n-1}}{\Delta t}\right)\right|\right|^2\leq C(||\phi_h^0||^2),
\end{equation}
where $b_{\text{max}}^{n-1}:=\max_{n=1\to N}||b(\phi_h^{n-1})||_{L^{\infty}(\Omega)}$. Finally, the function
\begin{align}
    \label{eqn:stab2}
&\frac{1}{\epsilon}(\psi_{1}(\phi_{h}^{n})+{\psi}_{2}(\phi_{h}^{n}),1)^h+\frac{\epsilon}{2} \left((J\star 1)_h\phi_{h}^{n},\phi_{h}^{n}\right)^h -\frac{\epsilon}{2}\left(J\star \phi_{h}^{n},\phi_{h}^{n}\right)^{h^2},
    \end{align}
is a decreasing Lyapunov function for the discrete solution.
\end{thm}
\begin{pf}
Recalling the definition \eqref{eqn:enfuncconv} and noting the convexity of $\psi_{1,\lambda}(\cdot)$, we can introduce a regularized lower semi continuous convex energy functional defined as
\begin{equation}
\label{eqn:enfuncconveps}
F_{1,\lambda}[\phi_{h,\lambda}^n]=\int_{\Omega}\left\{\frac{\epsilon}{2}(J\star 1)_h(\phi_{h,\lambda}^n)^2+
\frac{1}{\epsilon}\psi_1 (\phi_{h,\lambda}^n)+\chi_{\mathbb{R}^+}(\phi_{h,\lambda}^n)\right\}dx,
\end{equation}
and, endowing the space $S^h$ with the lumped scalar product \eqref{eqn:lump}, we can rewrite system \eqref{eqn:phlamb}, analogously to \eqref{eqn:phsubdiff}, as
\begin{equation}
\label{eqn:phlambsubdiff}
\begin{cases}
\displaystyle \biggl(\frac{\phi_{h,\lambda}^n-\phi_{h}^{n-1}}{\Delta t},\chi \biggr)^h+(b(\phi_{h}^{n-1})\nabla \mu_{h,\lambda}^n,\nabla \chi)=0,
\\
\displaystyle \left(\mu_{h,\lambda}^n-\frac{1}{\epsilon}\bar{\psi}'_{2}(\phi_{h}^{n-1})+\epsilon (J\star \phi_{h}^{n-1})_h,\xi -\phi_{h,\lambda}^n\right)^h+F_{1,\lambda}[\phi_{h,\lambda}^n]\leq F_{1,\lambda}[\xi],
\end{cases}
\end{equation}
for all $\chi,\xi \in S^h$. 

We now pass to the limit in \eqref{eqn:phlambsubdiff} as $\lambda \to 0$, considering the convergence properties \eqref{eqn:chepsconvlem} and \eqref{eqn:whepsconvlem}.
For any $\chi\in S^h $, we have
\begin{align}
\label{eqn:limeps}
\lim_{\lambda \rightarrow 0}\biggl(\frac{\phi_{h,\lambda}^n-\phi_{h}^{n-1}}{\Delta t},\chi \biggr)^h & =\biggl(\frac{\phi_{h}^n-\phi_{h}^{n-1}}{\Delta t},\chi \biggr)^h;\\
\label{eqn:limeps1} \lim_{\lambda \rightarrow 0}(b(\phi_{h}^{n-1})\nabla \mu_{h,\lambda}^n,\nabla \chi) & =(b(\phi_{h}^{n-1})\nabla \mu_{h}^n,\nabla \chi).
\end{align}
Since from \eqref{eqn:29} and \eqref{eqn:unfboundpsi1} we have that $\psi_{1,\lambda}(\phi_{h,\lambda}^n)\geq 0$ for $\lambda$ sufficiently small, and since, due to the convergence property \eqref{eqn:chepsconvlem}, \[\liminf_{\lambda\to 0}\psi_{1,\lambda}(\phi_{h,\lambda}^n)=\psi_{1}(\phi_{h}^n)\] if $\phi_{h}^n<1$, from the Fatou's Lemma and the semi continuity property of the indicator function $\chi_{\mathbb{R}^+}(\cdot)$, we deduce that
\begin{align*}
 \lim_{\lambda \rightarrow 0}F_{1,\lambda}[\phi_{h,\lambda}^n] &\geq F_{1}[\phi_{h}^n].
\end{align*}

\noindent
Taking $\xi=2\phi_{h,\lambda}^n-\phi_h^n$, where $\phi_h^n$ is the limit point in \eqref{eqn:chepsconvlem}, we obtain that
\[
\left(\mu_{h,\lambda}^n,\phi_{h,\lambda}^n-\phi_h^n\right)^h\leq F_{1,\lambda}[2\phi_{h,\lambda}^n-\phi_h^n]-F_{1,\lambda}[\phi_{h,\lambda}^n]+\left(\frac{1}{\epsilon}\bar{\psi}'_{2}(\phi_{h}^{n-1})-\epsilon (J\star \phi_{h}^{n-1})_h,\phi_{h,\lambda}^n-\phi_h^n\right)^h,
\]
from which, taking the limit as $\lambda \to 0$ and considering \eqref{eqn:chepsconvlem}, \eqref{eqn:unfboundpsi1} and \eqref{eqn:unifbound2}, we obtain that
\begin{equation}
    \label{eqn:mulimit1}
    \lim_{\lambda \to 0}\left(\mu_{h,\lambda}^n,\phi_{h,\lambda}^n-\phi_h^n\right)^h\leq 0.
\end{equation}
Also, taking $\xi=\phi_h^n$ in \eqref{eqn:chepsconvlem}, we obtain in a similar way that 
\begin{equation}
    \label{eqn:mulimit2}
    \lim_{\lambda \to 0}\left(\mu_{h,\lambda}^n,\phi_{h,\lambda}^n-\phi_h^n\right)^h\geq 0.
\end{equation}
Combining the latter results, we have that
\begin{equation}
    \label{eqn:mulimit3}
    \lim_{\lambda \to 0}\left(\mu_{h,\lambda}^n,\phi_{h,\lambda}^n-\phi_h^n\right)^h= 0.
\end{equation}
Now in order to pass to the limit in the first term in the left hand side of \eqref{eqn:phlambsubdiff} using the convergence property \eqref{eqn:whepsconvlem} we need to consider the subset $\tilde{S}^h(\phi_h^n)\subseteq S^h$ defined as $\tilde{S}^h(\phi_h^n):=\{q_h\in S^h:q_h(\mathbf{x}_j)=0 \; \text{if} \; \phi_h^n(\mathbf{x}_j)=0, j\in \mathcal{I}\}$. We observe that, as a consequence of \eqref{eqn:whepsconvlem} and \eqref{eqn:mulimit3}, we have that
\begin{equation}
    \label{eqn:mulimit4}
    \lim_{\lambda \to 0}\left(\mu_{h,\lambda}^n,\tilde{\xi}-\phi_{h,\lambda}^n\right)^h\rightarrow (\mu_{h}^n,\tilde{\xi}-\phi_{h}^n)^h,
\end{equation}
for any $\tilde{\xi}\in \tilde{S}^h(\phi_h^n)$. We conclude that the limit point $(\phi_{h}^n,\mu_{h}^n)$ satisfies, for each $(\chi,\tilde{\xi})\in S^h\times \tilde{S}^h(\phi_h^n)$,
\begin{equation}
\label{eqn:phsubdiff1}
\begin{cases}
\displaystyle \biggl(\frac{\phi_{h}^n-\phi_{h}^{n-1}}{\Delta t},\chi \biggr)^h+(b(\phi_{h}^{n-1})\nabla \mu_{h}^n,\nabla \chi)=0,
\\
\displaystyle \left(\mu_{h}^n-\frac{1}{\epsilon}{\psi}'_{2}(\phi_{h}^{n-1})+\epsilon (J\star \phi_{h}^{n-1})_h,\tilde{\xi}-\phi_{h}^n\right)^h+F_{1}[\phi_{h}^n]\leq F_{1}[\tilde{\xi}].
\end{cases}
\end{equation}
Note that $\bar{\psi_2}'(\phi_h^{n-1})\equiv {\psi_2}'(\phi_h^{n-1})$ since $\phi_h^{n-1}<1$. 
Finally, since $\phi_{h}^{n}<1$ (see \eqref{eqn:unifbound2}) and $\psi_{1}(\phi_{h}^{n})$ is convex and smooth for $\phi_{h}^{n}<1$, system \eqref{eqn:phsubdiff1} is equivalent to system \eqref{pbm:ph} (see \eqref{eqn:phsubdiff}), valid for any $\xi\in K^h$, with $\xi(\mathbf{x}_j)=0$ if $\phi_h^n(\mathbf{x}_j)=0$. Recalling Remark \ref{rmk:ineq}, note that if \eqref{pbm:ph} is valid for each $(\chi,\xi)\in S^h\times K^h$, with $\xi(\mathbf{x}_j)=0$ if $\phi_h^n(\mathbf{x}_j)=0$, it is also valid for each $(\chi,\xi)\in S^h\times K^h$. Hence the limit point $(\phi_{h}^n,\mu_{h}^n)$ is a solution of Problem $P^{h}$. 

The uniqueness of the solution $(\phi_{h}^n,\mu_{h}^n)$ of Problem $P^{h}$ can be inferred similarly as in Lemma \ref{lem:discex}. For what concerns the stability estimate \eqref{eqn:stab1}, 
choosing $\chi=\mu_{h}^{n}$ in the first equation of \eqref{pbm:ph} and $\xi=\phi_{h}^{n-1}$ in the second equation of \eqref{pbm:ph}, similarly to \eqref{eqn:estimateeps} and \eqref{eqn:estimateeps2} we obtain that
\begin{align}
\label{eqn:estimatelim}
&\notag \frac{1}{\epsilon}(\psi_{1}(\phi_{h}^{n})+{\psi}_{2}(\phi_{h}^{n}),1)^h+\frac{\epsilon}{2} \left((J\star 1)_h\phi_{h}^{n},\phi_{h}^{n}\right)^h-\frac{\epsilon}{2}\left(J\star \phi_{h}^{n},\phi_{h}^{n}\right)^{h^2}+\Delta t||[b(\phi_{h}^{n-1})]^{1/2}\nabla \mu_{h}^{n}||^{2}\\
& \notag +\frac{\epsilon}{2} \sum_{i,j\in \mathcal{I}}J(\mathbf{x}_i-\mathbf{x}_j)\left(\phi_{h}^{n}(\mathbf{x}_i)-\phi_{h}^{n-1}(\mathbf{x}_i)+\phi_{h}^{n}(\mathbf{x}_j)-\phi_{h}^{n-1}(\mathbf{x}_j)\right)^2(1,\chi_i)(1,\chi_j)\\
& \leq \frac{1}{\epsilon}(\psi_{1}(\phi_{h}^{n-1})+{\psi}_{2}(\phi_{h}^{n-1}),1)^h+\frac{\epsilon}{2} \left((J\star 1)_h\phi_{h}^{n-1},\phi_{h}^{n-1}\right)^h-\frac{\epsilon}{2}\left(J\star \phi_{h}^{n-1},\phi_{h}^{n-1}\right)^{h^2},
\end{align}
from which we get that \eqref{eqn:stab2} is a decreasing Lyapunov function for the discrete solution. Summing \eqref{eqn:estimatelim} from $n=1\to N$, using \eqref{eqn:jhcomp}, \eqref{eqn:interp5} and \eqref{eqn:estimateeps3} and observing that $0\leq \phi_h^0<1$, we obtain the first two bounds in \eqref{eqn:stab1}. Choosing now $\chi=\hat{\mathcal{G}}^h\left(\frac{\phi_{h}^{n}-\phi_{h}^{n-1}}{\Delta t}\right)$ in the first equation of system \eqref{pbm:ph}, using \eqref{eqn:greendiscr} and Cauchy-Schwarz and Young inequalities, we get
\begin{align*}
&\biggl(\frac{\phi_{h}^{n}-\phi_{h}^{n-1}}{\Delta t},\hat{\mathcal{G}}^h\left(\frac{\phi_{h}^{n}-\phi_{h}^{n-1}}{\Delta t}\right)\biggr)^h =\left|\left|\nabla \hat{\mathcal{G}}^h\left(\frac{\phi_{h}^{n}-\phi_{h}^{n-1}}{\Delta t}\right)\right|\right|^{2}\\
&=-\biggr(b(\phi_{h}^{n-1})\nabla \mu_{h}^{n},\nabla \hat{\mathcal{G}}^h\left(\frac{\phi_{h}^{n}-\phi_{h}^{n-1}}{\Delta t}\right)\biggl) \leq \frac{1}{2}||b(\phi_{h}^{n-1})\nabla \mu_{h}^{n}||^{2}+\frac{1}{2}\left|\left|\nabla \hat{\mathcal{G}}^h\left(\frac{\phi_{h}^{n}-\phi_{h}^{n-1}}{\Delta t}\right)\right|\right|^{2}\\
&\leq Cb_{\text{max}}^{n-1}||[b(\phi_{h}^{n-1})]^{1/2}\nabla \mu_{h}^{n}||^{2}+\frac{1}{2}\left|\left|\nabla \hat{\mathcal{G}}^h\left(\frac{\phi_{h}^{n}-\phi_{h}^{n-1}}{\Delta t}\right)\right|\right|^{2}.
\end{align*}
Using \eqref{eqn:estimatelim} and summing from $n=1\rightarrow N$ we get the bound for the last term in \eqref{eqn:stab1}.
\end{pf}
\subsection{Numerical algorithm}
We propose a numerical algorithm which solves the discrete variational inequality \eqref{pbm:ph} by directly solving the KKT conditions associated to the minimization problem \eqref{eqn:lagrangian} (without regularization). Let us introduce the following matrices (expressed by components)
\begin{align*}
    &M_{ij}:=(\chi_j,\chi_i)^h;\\
    &A_{\phi_h^{n-1},ij}:=(b(\phi_h^{n-1})\nabla \chi_j, \nabla \chi_i);\\
    &J_{1,ij}:=\left[\sum_{k\in \mathcal{I}}J(\mathbf{x}_i-\mathbf{x}_k)(1,\chi_i)(1,\chi_k)\right]\delta_{ij};\\
    &J_{2,ij}:=J(\mathbf{x}_i-\mathbf{x}_j)(1,\chi_i)(1,\chi_j),
\end{align*}
with $i,j\in \mathcal{I}$ and where $\delta_{ij}$ is the Kronecker symbol. We also introduce the symmetry-preserving modified matrices on the passive nodes
\begin{align*}
    &\hat{M}_{ij}:=
    \begin{cases}
     M_{i,j} \quad \text{if}\; i,j\notin J_0(\phi_h^{n-1}),\\
     1 \quad \text{if}\; i=j\in J_0(\phi_h^{n-1}),\\
     0 \quad \text{if}\; i\in J_0(\phi_h^{n-1}),j\neq i\\
     0 \quad \text{if}\; j\in J_0(\phi_h^{n-1}),i\neq J\\
    \end{cases}\quad 
    \hat{A}_{\phi_h^{n-1},ij}:=
    \begin{cases}
     {A}_{\phi_h^{n-1},ij}\quad \text{if}\; i,j\notin J_0(\phi_h^{n-1}),\\
     1 \quad \text{if}\; i=j\in J_0(\phi_h^{n-1}),\\
     0 \quad \text{if}\; i\in J_0(\phi_h^{n-1}),j\neq i\\
     0 \quad \text{if}\; j\in J_0(\phi_h^{n-1}),i\neq J,\\
    \end{cases}
\end{align*}
with $i,j\in \mathcal{I}$. We observe that ${A}_{\phi_h^{n-1},ij}$ is invertible. Let us introduce the quantities $N_0:=\#\left(J_0(\phi_h^{n-1})\right)$, $N_+:=\#\left(\mathcal{I}\setminus J_0(\phi_h^{n-1})\right)=N_t-N_0$, where $N_t:=\#(\mathcal{I})$. We introduce the matrices $K\in \mathbb{R}^{N_+\times N_t}$, whose rows are vectors in $\mathbb{R}^{N_t}$ given by the one-hot encoding representations of the active nodes, and $L\in \mathbb{R}^{N_0\times N_t}$, whose rows are vectors in $\mathbb{R}^{N_t}$ given by the one-hot encoding representations of the passive nodes. We indicate with the notation $\bar{v}_h$ the vector of components $v_h(\mathbf{x}_j)$ for any $v_h\in S^h$, $j\in \mathcal{I}$. From \eqref{pbm:ph}$_1$ we have that 
\begin{equation}
    \label{eqn:mualg}
    \mu_h^n=-\frac{1}{\Delta t}\left(\hat{A}_{\phi_h^{n-1}}\right)^{-1}\hat{M}\left(\bar{\phi}_h^n-\bar{\phi}_h^{n-1}\right).
\end{equation}
Substituting \eqref{eqn:mualg} in \eqref{pbm:ph}$_2$, and introducing the matrix
\[
Q:=\epsilon J_1+\frac{1}{\Delta t}\left(\hat{A}_{\phi_h^{n-1}}\right)^{-1}\hat{M},
\]
we obtain a discrete variational inequality which can be rephrased, considering also the constraints on the passive nodes, as the complementarity problem
\begin{equation}
    \label{eqn:compl}
    \begin{cases}
    Q\bar{\phi}_h^n+\frac{1}{\epsilon}{M}\left(\psi_1'(\bar{\phi}_h^n)+\psi_2'(\bar{\phi}_h^{n-1})\right)-\epsilon J_2\bar{\phi}_h^{n-1}-\frac{1}{\Delta t}{M}\left(\hat{A}_{\phi_h^{n-1}}\right)^{-1}\hat{M}\bar{\phi}_h^{n-1}-L^T\bar{\rho}-K^T\bar{\eta}=0,\\
    0\leq \bar{\eta} \perp K\bar{\phi}_h^n\geq 0,\\
    L\bar{\phi}_h^n=0,
    \end{cases}
\end{equation}
where $\bar{\eta}\in \mathbb{R}^{N_+}$ and $\bar{\rho}\in \mathbb{R}^{N_0}$ are the Lagrange multipliers of the inequality constraint on the active nodes and of the equality constraint on the passive nodes respectively. We solve the KKT System \eqref{eqn:compl} employing a null space method, i.e. reducing the KKT system onto the null space of the matrix $L$, via a preconditioned accelerated gradient method. Observing that $LK^T=\boldsymbol{0}$, and defining a reduced vector $\tilde{v}_h\in \mathbb{R}^{N_+}$ such that $\bar{v}_h=K^T\tilde{v}_h$, we formulate the following algorithm
\begin{algorithmic}
\Require $\alpha_0>0$ (an acceleration parameter), $\phi_h^{n-1},\mu_h^{n-1}$;\\
\textbf{Step 1}
\For{$k\geq 0$}  \\
\textbf{Initialization}\[\phi_h^{n,0}=\phi_h^{n-1},\mu_h^{n,0}=\mu_h^{n-1};\]
Find $\phi_h^{n,k+1}\in K^h$ such that:
\If {$j \in J_{0}(\phi_{h}^{n-1})$}
    \State $\phi_{h}^{n,k+1}(\mathbf{x}_j)\gets \phi_{h}^{n-1}(\mathbf{x}_j)$
\Else
        \begin{align}
        \label{eqn:varinequal}
        & \notag \tilde{\phi}_{h,j}^{n,k+1}=\max \biggl(0,\tilde{\phi}_{h,j}^{n,k}-\alpha_k(\text{diag}(KQK^T))^{-1}\biggl[KQK^T\tilde{\phi}_h^n\\
        & +\frac{1}{\epsilon}K{M}K^T\left(\psi_1'(\tilde{\phi}_h^n)+\psi_2'(\tilde{\phi}_h^{n-1})\right)-\epsilon KJ_2K^T\tilde{\phi}_h^{n-1}-\frac{1}{\Delta t}K{M}\left(\hat{A}_{\phi_h^{n-1}}\right)^{-1}\hat{M}K^T\tilde{\phi}_h^{n-1}\biggr]\biggr)          
        \end{align}
        \EndIf  \\

\If {$||\phi_{h}^{n,K+1}-\phi_{h}^{n,K}||_{L^{\infty}(\Omega)}<10^{-6}$}
    \State $\phi_h^n\gets \phi_{h}^{n,K+1}$; \textbf{break}.
\EndIf
\EndFor \\
\textbf{Step 2}
Find $\mu_h^n\in S^h$ such that
\[
\mu_h^n=-\frac{1}{\Delta t}\left(\hat{A}_{\phi_h^{n-1}}\right)^{-1}\hat{M}\left(\bar{\phi}_h^n-\bar{\phi}_h^{n-1}\right).
\]
\end{algorithmic}
The acceleration parameter $\alpha_k$ is dinamically chosen by a projected search in such a way that the functional associated to \eqref{eqn:compl} is decreased at each iterative step $k$. Note that, since the operator acting on $\tilde{\phi}_h^n$ in the square bracket in \eqref{eqn:varinequal} is continuous and strictly monotone, the
projection map defined in \eqref{eqn:varinequal} has a unique fixed point (see e.g. \cite[Chapter 2]{temam} for details). Since the support of the discrete solution can move at most for a length $h_K$ locally at each time step (see Remark \ref{rem:deggreen}), where $h_K$ is the diameter of an element $K$ which contains a passive node, at each time step we adapt the value of $\Delta t$ such that the CFL-type condition 
\[
\Delta t<\Delta t_{CFL}:=\frac{\min_{K\in \mathcal{T}_h}h_K}{\max_{\Omega}|\mathbf{v}_h^n|_{\infty}}
\]
is guaranteed, where $\mathbf{v}_h^n$ may be defined from \eqref{eqn:22} as the velocity of the phase $\phi_h^n$ as
\[
\mathbf{v}_h^n:=-\left(\phi_h^{n-1}\right)^{\alpha-1}(1-\phi_h^{n-1})^2\nabla \mu_h^n.
\]
\begin{rem}
    \label{rem:j12}
    Note that the matrices $J_1$ and $J_2$ in \eqref{eqn:compl}, associated to the nonlocal terms, are assembled once and for all as an initialization at the beginning of the simulation before entering the time step cycles. We practically calculate their entries only for those couples of dofs $i,j\in \mathcal{I}$ such that $J(\mathbf{x}_i-\mathbf{x}_j)>\text{tol}$, where $\text{tol}$ is a small tolerance which we take as $\text{tol}=10^{-6}$. Hence, the computation of the nonlocal terms is computationally fast and efficient. {We remark that in three space dimensions the computation of these terms for fine meshes is still not affordable. We will design a proper parallelization technique based on domain decomposition to efficiently calculate them for three dimensional applications as a future development. For now, we will show simulation results in two space dimensions.}
\end{rem}

\section{Numerical simulations} 
\label{Simulations}
In this section we report the results of numerical simulations in two space dimensions in order to show the stability and the qualitative properties of the solution of the proposed numerical scheme. In particular, we investigate the phase separation dynamics
and the coarsening dynamics associated to the model \eqref{eqn:22}-\eqref{eqn:23}, starting from different initial configurations given by a small uniformly distributed random perturbation around constant values
belonging to the metastable regime of the single-well potential \eqref{eqn:4}.

The values of the parameters are taken as $\epsilon=0.014$, $\bar{\phi}=0.6$ and $M=1$. We consider the following form for the kernel $J$, which satisfies \eqref{eqn:5}:
\[
J(\mathbf{x},\mathbf{y})=\frac{1}{\epsilon^4}\mathrm{e}^{-\frac{x^2+y^2}{2\epsilon^2}}.
\]
Moreover, we consider a domain $\Omega=[-1,1]\times[-1,1]$, with a uniform triangulation of dimension $80\times80$, and $\Delta t=\min\left(0.1\epsilon^2,\Delta t_{CFL}\right)$. 
\begin{rem}
    \label{rem:minjstar1}
    We observe that $J$ is symmetric and that $J(\cdot)>0$, hence the hypothesis of Theorem \ref{th:wellposednessph} are satisfied. Moreover,
    by direct computation with the chosen values of the parameters we obtain that
    \[
    \epsilon \inf_{\Omega}(J\star 1)_h\sim 37.008>\frac{1-\bar{\phi}}{\epsilon}\sim 28.5714>\frac{2+(1-\bar{\phi})-3\sqrt[3]{1-\bar{\phi}}}{\epsilon}\sim 13.5415.
    \]
    Hence, the hypothesis of Theorem \ref{thm:3} are satisfied: at the continuous level, we expect that $\delta < \phi(\mathbf{x},t)<1-\delta$ for some $\delta>0$ for a.e. $(\mathbf{x},t)\in \Omega\times [0,T]$.
\end{rem}
We consider test cases with $\alpha=1$ and $\alpha=2$ in the form of the mobility \eqref{eqn:23a}, in order to investigate if the discrete dynamics preserves the strict separation property of the continuous solution predicted by Theorem \ref{thm:3} (cfr. Remark \ref{rem:minjstar1}).

Hence, we consider numerical tests with $\phi_0=0.05\pm 0.025 \iota$, $\phi_0=0.3\pm 0.15 \iota$ and $\phi_0=0.36\pm 0.072 \iota$, where $\iota$ is a random perturbation uniformly distributed in the interval $[0, 1]$. In these cases we expect the system dynamics to undergo an initial transitory regime of spinodal decomposition, after which a stationary state with phase-separated domains characterized by different topologies, depending on the initial datum, is reached. We also compare the results for the non-local model \eqref{eqn:22}-\eqref{eqn:23}, with the given parameters, with the results for the corresponding local discrete model, with the same parameters, introduced and analyzed in \cite{agosti1}, in order to depict the main characteristics of the non-local dynamics with respect to the local one. 

\begin{figure}[h!]
\includegraphics[width=1.0\linewidth]
{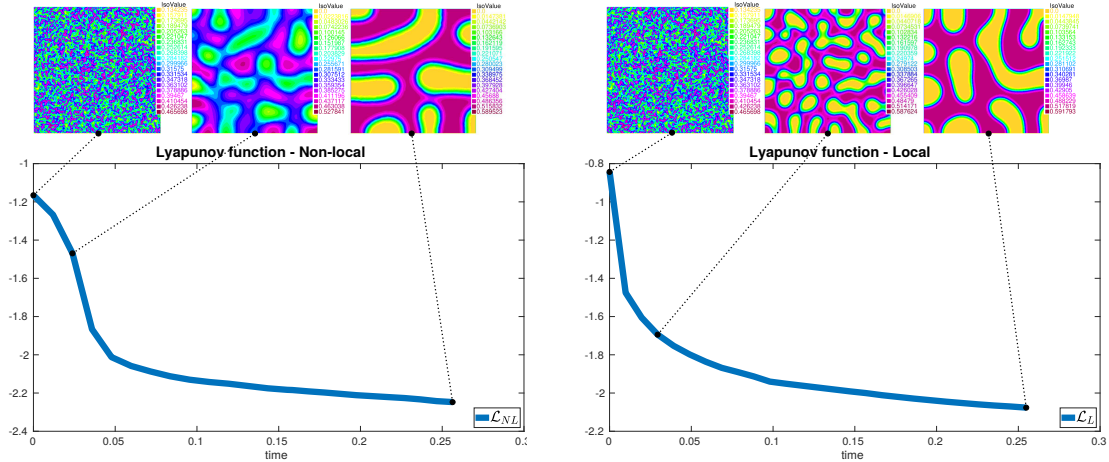}
\centering
\caption{Plots of the Lyapunov functionals $\mathcal{L}_{NL}$ and $\mathcal{L}_L$ vs time for the non-local and the local model with $\alpha=1$ and $\phi_0=0.3\pm 0.15\iota$, together with inserts showing the plots of $\phi_h^n$ at different time points.}
\label{fig:2}
\end{figure}

We start by verifying the gradient stability of the numerical scheme \eqref{pbm:ph}. In Figure \ref{fig:2} we report the plots of the Lyapunov functional \eqref{eqn:stab2}, which we call $\mathcal{L}_{NL}$, and of the Lyapunov functional for the local model, i.e.
\[
\mathcal{L}_{L}:=\frac{\epsilon}{2}||\nabla \phi_h^n||^2+\frac{1}{\epsilon}(\psi_1(\phi_h^n)+\psi_2(\phi_h^n)),
\]
versus time, solving the non-local and the local models in the case $\alpha=1$ and $\phi_0=0.3\pm 0.15\iota$.
We observe from Figure \ref{fig:2} that both the non-local and the local discrete schemes are gradient stable. In the non-local case, the spinodal decomposition is slower. As shown in the insert with the plots of $\phi_h^n$ at time $t=0.03$, i.e. during the spinodal decomposition when the Lyapunov functions undergo a steep decrease, in the non-local case the domains of pure phases are mixing on longer spatial ranges than the local case due to the non-local interaction.

\begin{figure}[hb!]
\includegraphics[width=0.7\linewidth]
{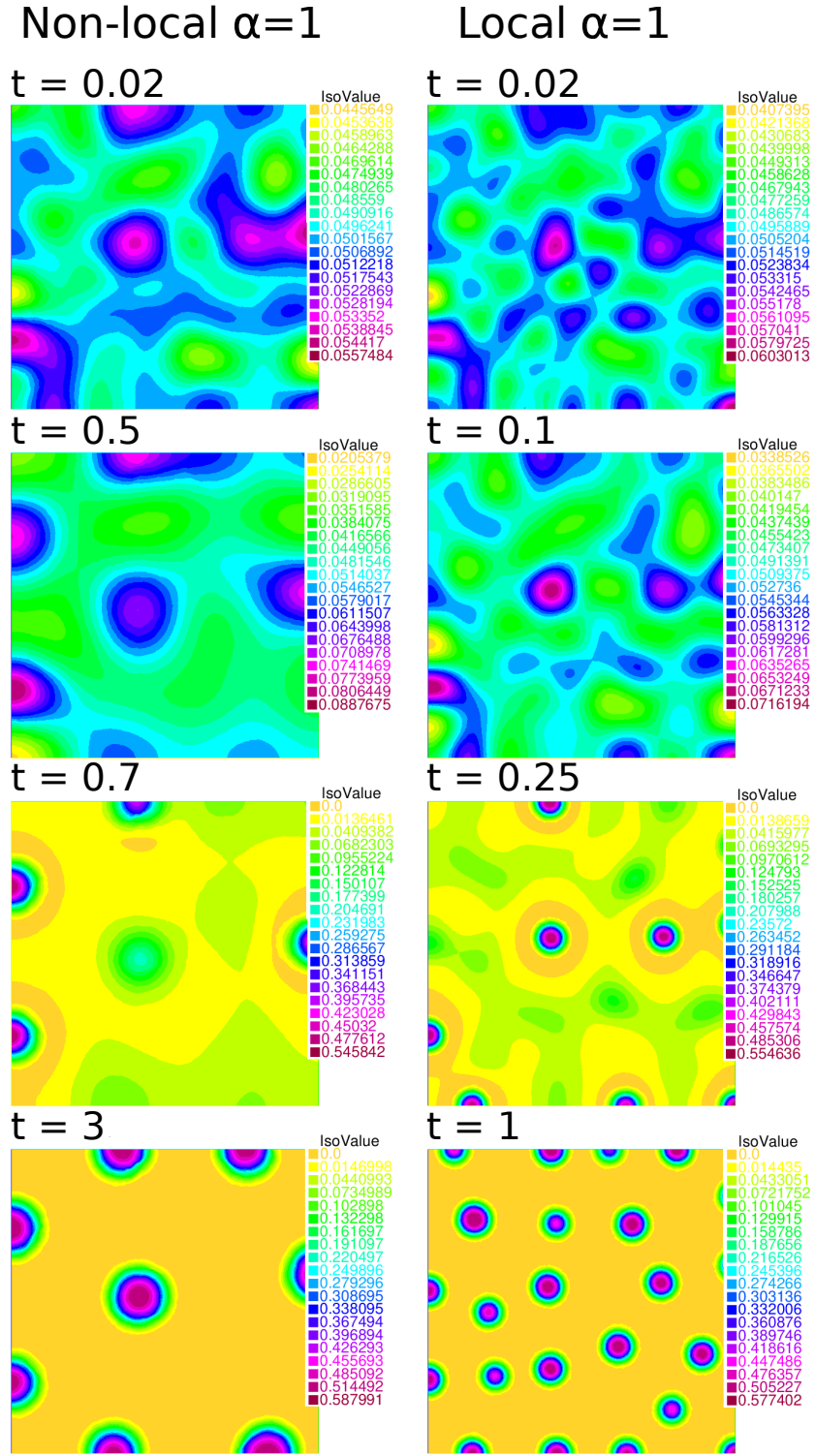}
\centering
\caption{Plot of $\phi_h^n$ at different time points in the case $\phi_0=0.05\pm 0.025\iota$ and $\alpha=1$, for the non-local model (first column) and the local model (second column).}
\label{fig:3}
\end{figure}

In Figure \ref{fig:3} we report the numerical results in the case $\phi_0=0.05\pm 0.025\iota$, showing the plots of $\phi_h^n$ at different time points throughout the phase separation dynamics, up to late times at
which we can observe the coarsening dynamics of the separated domain subregions, considering $\alpha=1$ and comparing between the non-local and the local models.
We observe from Figure \ref{fig:3} that, as expected (see e.g. \cite{bray,agosti1}), the phase separation dynamics consists in the formation of circular clusters with $\phi\sim \bar{\phi}$ (the expected value of $\phi$ inside the clusters in the stationary state can be obtained intercepting the graph of the single well potential with its tangent line passing at $\phi=0$ \cite{bray}) immersed in a bath with $\phi\sim 0$. Comparing the results between the non-local and the local models, we observe that in the non-local case the dynamics is slower and that the interface regions of domains of separated phases are wider than in the local case. Also, the number of clusters is lower and the clusters are bigger in the non-local case. Indeed, due to the non-local interaction between the phases, smaller neighboring clusters tend to merge into bigger ones in the non-local case, while in the local case they do not interact.
{As explained in the Introduction, these non local dynamics may be more appropriate than its local counterpart to represent long-range interactions in cells-cells and cells-matrix adhesion during the invasive and metastatic phenotype of cancer progression \cite{armstrong,hillen}.}

\begin{figure}[ht!]
\includegraphics[width=0.8\linewidth]
{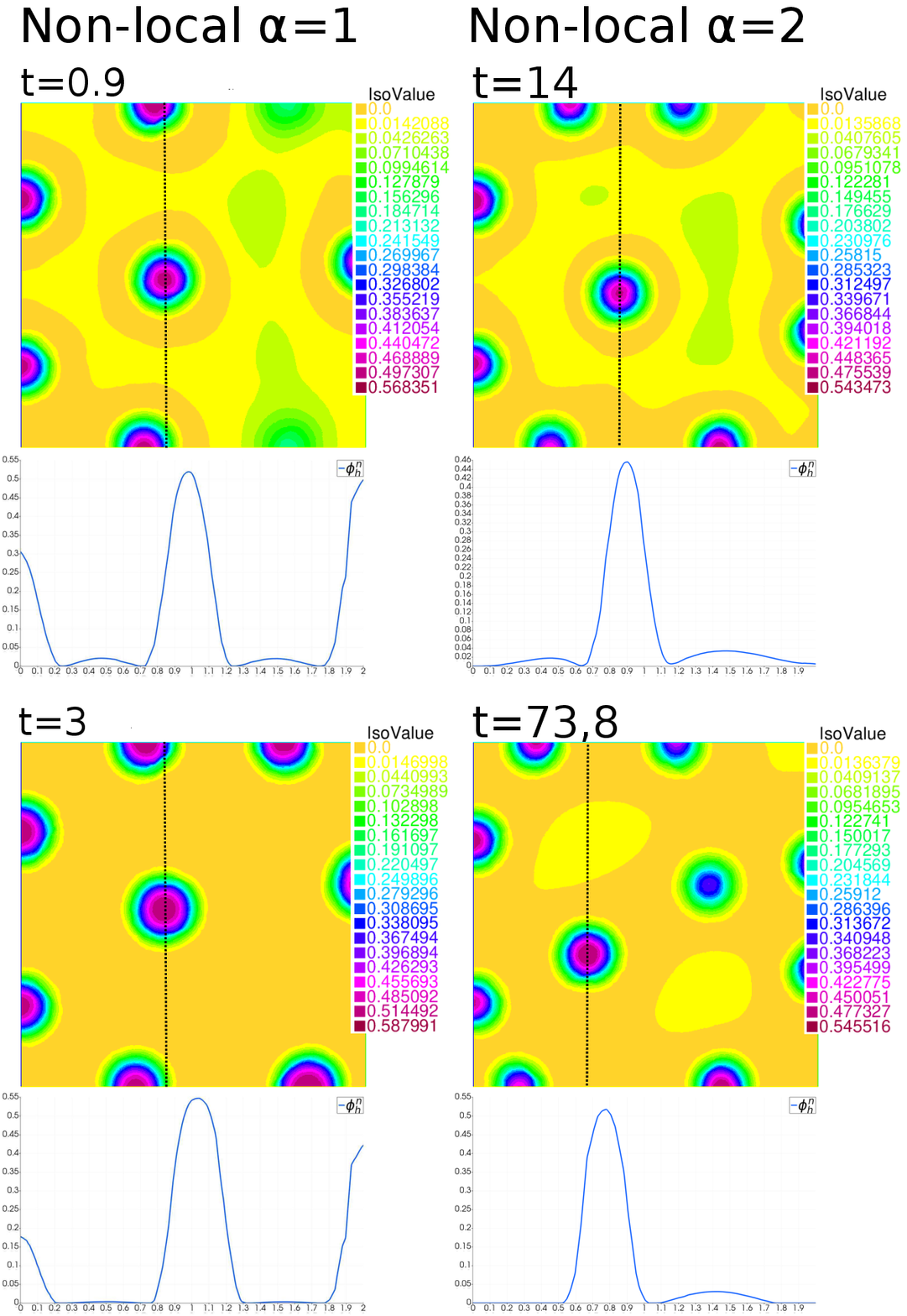}
\centering
\caption{Plot of $\phi_h^n$ at different time points at late times for the non-local models with $\alpha=1$ (left column) and $\alpha =2$ (right column), in the case $\phi_0=0.05\pm 0.025\iota$, together with plots over vertical lines of $\phi_h^n$. The vertical lines are indicated as dotted lines in the two dimensional plots.}
\label{fig:4}
\end{figure}

In Figure \ref{fig:4} we compare the plots of $\phi_h^n$ at different time points at late times for the non-local models with $\alpha=1$ and $\alpha =2$. We also show the plots over vertical lines of $\phi_h^n$, where the vertical lines are shown as dotted lines in the corresponding two dimensional plots and are chosen so that they intersect at least one cluster throughout the domain.
We observe from Figure \ref{fig:4} that the separation and coarsening dynamics are much slower in the case $\alpha=2$ with respect to the case $\alpha=1$. Moreover, the line plots of $\phi_h^n$ show that in the case $\alpha=1$ the phase variable attains the value zero at interface regions during cluster formation, while in the case $\alpha=2$ the value zero is not attained during cluster formation. At the very late times also in the case $\alpha=2$ the value zero is attained. This is an indication that the strict separation property predicted by Theorem \ref{thm:3} (see also Remark \ref{rem:minjstar1}) is valid also for the discrete solution, up to numerical errors.

\begin{figure}[h!]
\includegraphics[width=0.7\linewidth]
{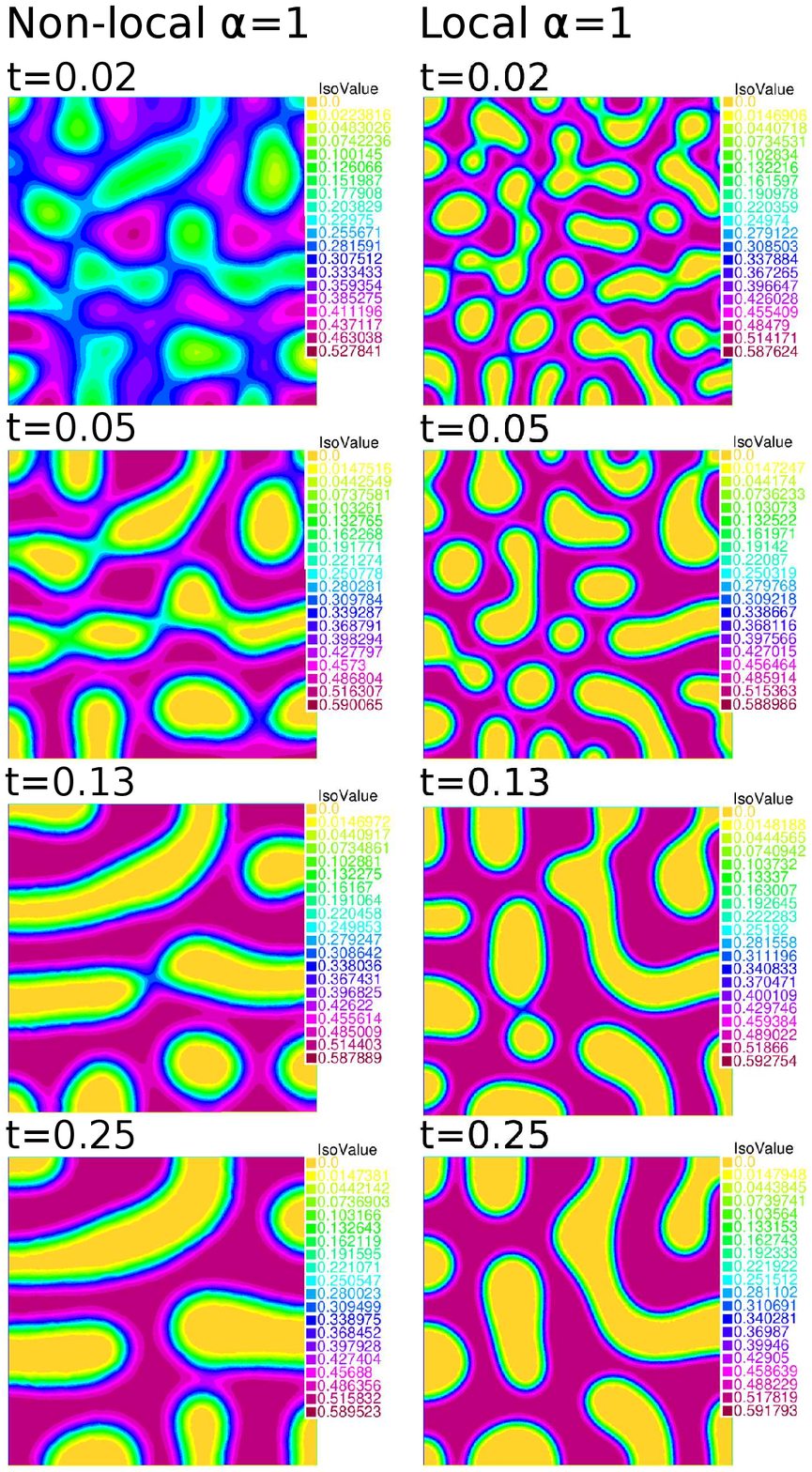}
\centering
\caption{Plot of $\phi_h^n$ at different time points in the case $\phi_0=0.3\pm 0.15\iota$ and $\alpha=1$, for the non-local model (first column) and the local model (second column).}
\label{fig:5}
\end{figure}

In Figure \ref{fig:5} we report the numerical results in the case $\phi_0=0.3\pm 0.15\iota$, showing the plots of $\phi_h^n$ at different time points, considering $\alpha=1$ and comparing between the non-local and the local models.
We observe from Figure \ref{fig:5} that, as expected (see e.g. \cite{bray,agosti1}), the phase separation dynamics consists in the formation of strips with $\phi\sim \bar{\phi}$ immersed in a bath with $\phi\sim 0$. Comparing the results between the non-local and the local models, we observe that in the non-local case the interface regions of the strips are wider than in the local case. Also, due to the non-local interaction between the phases neighboring strips tend to interact and merge during the coarsening dynamics, while in the local case they do not interact.

\begin{figure}[h!]
\includegraphics[width=0.8\linewidth]
{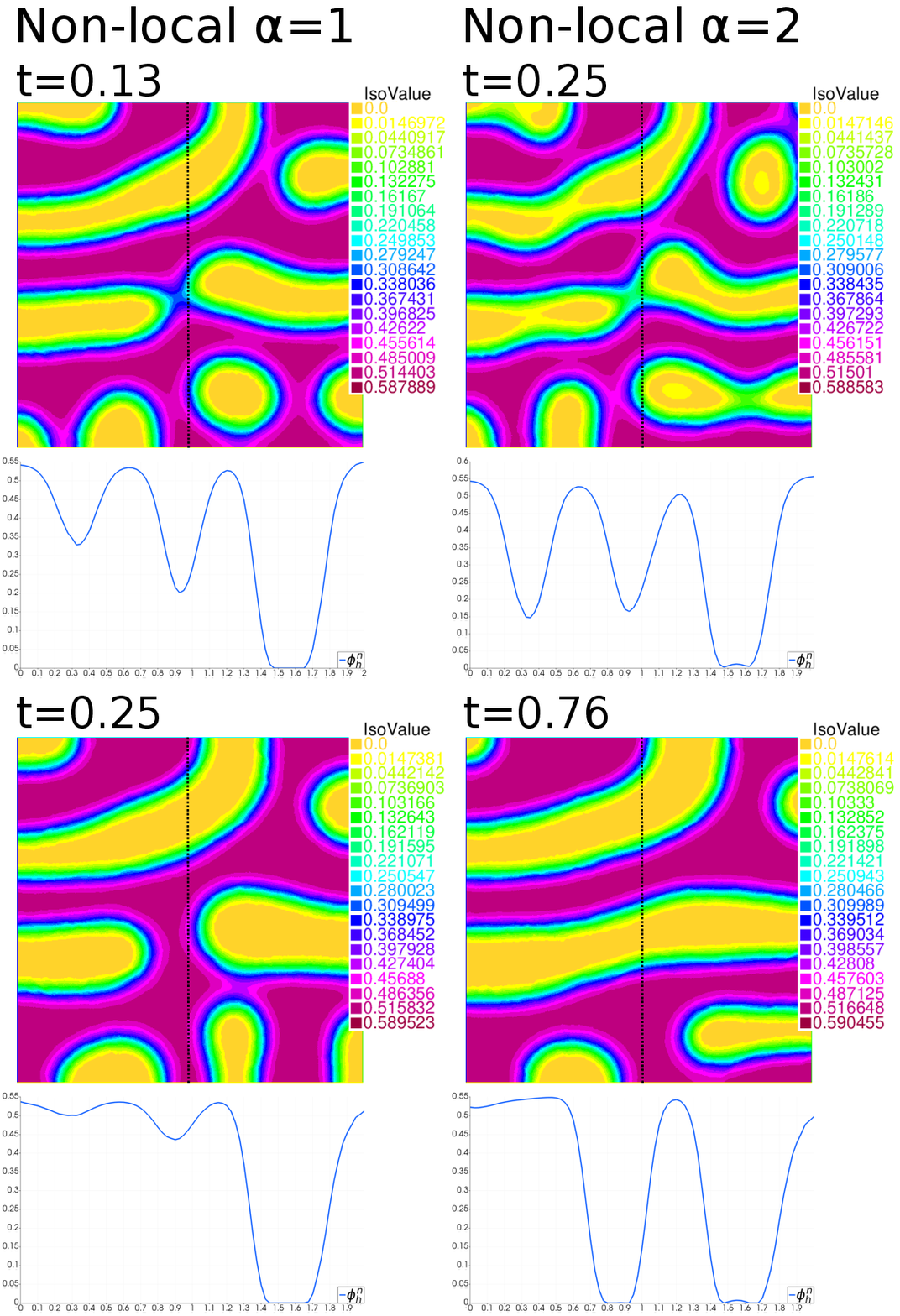}
\centering
\caption{Plot of $\phi_h^n$ at different time points at late times for the non-local models with $\alpha=1$ (left column) and $\alpha =2$ (right column), in the case $\phi_0=0.3\pm 0.15\iota$, together with plots over vertical lines of $\phi_h^n$. The vertical lines are indicated as dotted lines in the two dimensional plots.}
\label{fig:6}
\end{figure}

As in Figure \ref{fig:4}, in Figure \ref{fig:6} we compare the plots of $\phi_h^n$ at different time points at late times for the non-local models with $\alpha=1$ and $\alpha =2$, showing also the plots over vertical lines of $\phi_h^n$.

We observe from Figure \ref{fig:6} that the separation and coarsening dynamics are slower in the case $\alpha=2$ with respect to the case $\alpha=1$. Moreover, the line plots of $\phi_h^n$ show that in the case $\alpha=1$ the phase variable attains the value zero at interface and in the core regions of the strips, while in the case $\alpha=2$ the value zero is not attained in the strips core regions, while it is touched at their interface. This again is an indication that the strict separation property predicted by Theorem \ref{thm:3} is valid also for the discrete solution, up to numerical errors.

\begin{figure}[h!]
\includegraphics[width=0.7\linewidth]
{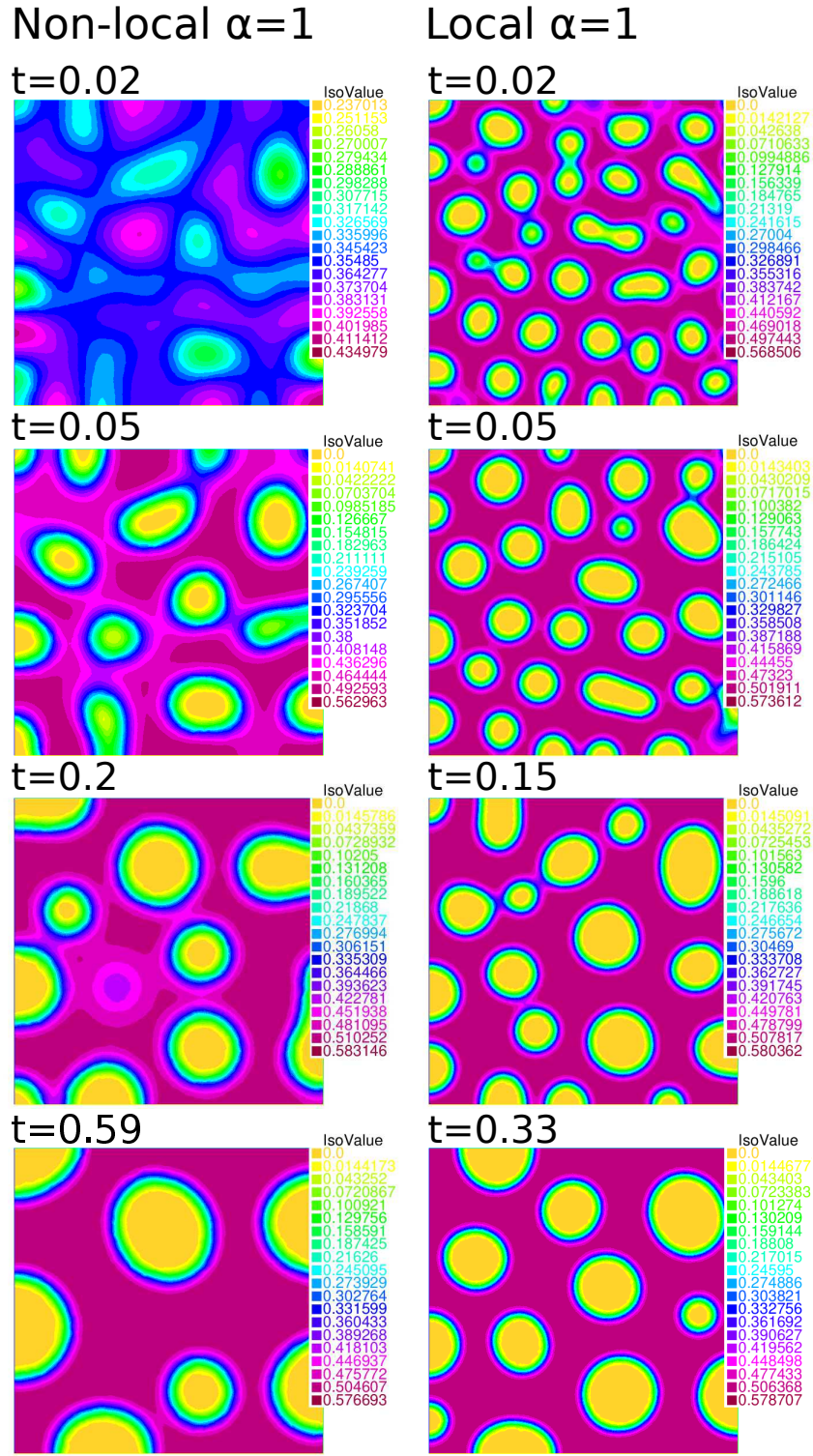}
\centering
\caption{Plot of $\phi_h^n$ at different time points in the case $\phi_0=0.36\pm 0.072\iota$ and $\alpha=1$, for the non-local model (first column) and the local model (second column).}
\label{fig:7}
\end{figure}

Finally, in Figure \ref{fig:7} we report the numerical results in the case $\phi_0=0.36\pm 0.072\iota$, showing the plots of $\phi_h^n$ at different time points, considering $\alpha=1$ and comparing between the non-local and the local models.
We observe from Figure \ref{fig:7} that, as expected (see e.g. \cite{bray,agosti1}), the phase separation dynamics consists in the formation of circular clusters with $\phi\sim 0$ immersed in a bath with $\phi\sim \bar{\phi}$. Comparing the results between the non-local and the local models, we observe that in the non-local case the dynamics is slower and that the interface regions of domains of separated phases are wider than in the local case. Also, the number of clusters is lower and the clusters are bigger in the non-local case. Indeed, due to the non-local interaction between the phases, smaller neighboring clusters tend to merge into bigger ones in the non-local case, while in the local case they do not interact.

\begin{figure}[h!]
\includegraphics[width=0.8\linewidth]
{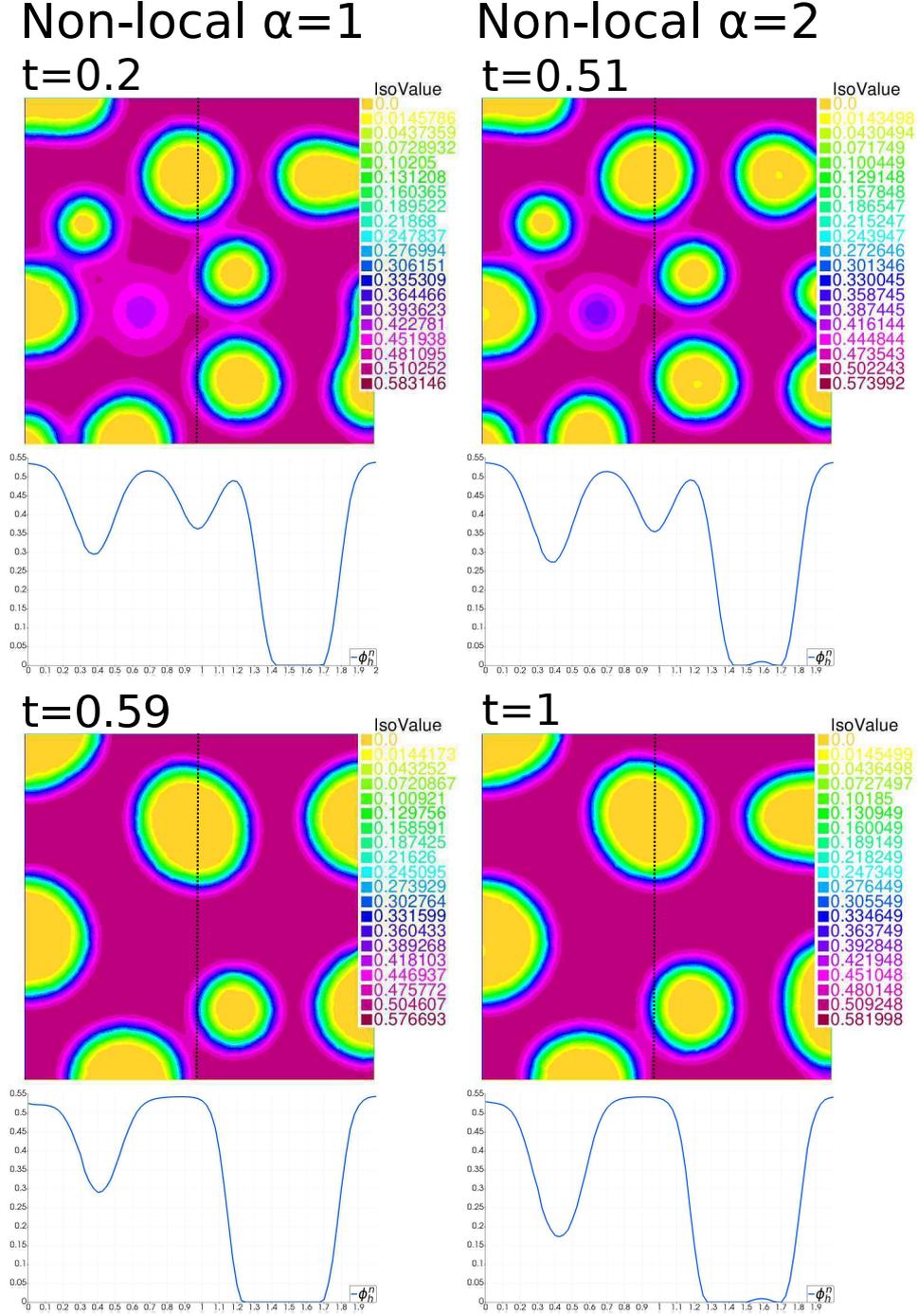}
\centering
\caption{Plot of $\phi_h^n$ at different time points at late times for the non-local models with $\alpha=1$ (left column) and $\alpha =2$ (right column), in the case $\phi_0=0.36\pm 0.072\iota$, together with plots over vertical lines of $\phi_h^n$. The vertical lines are indicated as dotted lines in the two dimensional plots.}
\label{fig:8}
\end{figure}

In Figure \ref{fig:8} we compare the plots of $\phi_h^n$ at different time points at late times for the non-local models with $\alpha=1$ and $\alpha =2$. We also show the plots over vertical lines of $\phi_h^n$, where the vertical lines are chosen so that they intersect at least one cluster throughout the domain.
We observe from Figure \ref{fig:8} that the separation and coarsening dynamics are slower in the case $\alpha=2$ with respect to the case $\alpha=1$. Moreover, the line plots of $\phi_h^n$ show that in the case $\alpha=1$ the phase variable attains the value zero at interface and in the core regions of the clusters, while in the case $\alpha=2$ the value zero tends to be not attained in the core clusters' regions, while it is touched at their interface. This again is an indication that the strict separation property predicted by Theorem \ref{thm:3} is valid also for the discrete solution, up to numerical errors.

\section{Conclusions}
In this paper we presented a general mixture model for tumor growth consisting in  a binary, closed and incompressible mixture of tumor and healthy cells, evolving in a bounded domain. We took into account the viscous properties of the mixture, due to the filtration processes between the phases, the short-range adhesion properties between the cells, expressed in terms of a singular single-well potential of the Lennard--Jones type, and the long-range interactions between the phases, expressed in terms of a convolution operator with a localized symmetric kernel. We also took into account cells proliferation and death and chemotactic phenomena through the coupling with a massless nutrient species. Employing a generalized form of the least dissipation principle, we obtained a coupled system of non linear PDEs given by a
Darcy-type equation for the mixture velocity field coupled with a convective reaction–diffusion type equation for the nutrient concentration and a non–local convective and degenerate Cahn--Hilliard equation for the tumor phase. For biological applications the use of a single-well potential to describe the cells-cells adhesion properties is more physically appropriate than the use of double-well potentials, which is typically employed in literature. 
Since the degeneracy set of the mobility and the set of singularities of the single-well potential do not coincide, classical analytical results about the existence and the strict separation property for weak solutions available in literature for double-well potentials are not available in the present case.
We then obtained a global existence result {for $d\leq 3$} for a weak solution of a simplified version of the model, which considers only the non-local and degenerate Cahn--Hilliard equation without growth and advection. We also proved a strict separation property {in 3D} for the weak solution under certain hypothesis on the convolution kernel and depending on the degree of degeneracy of the cells mobility. As a consequence of the strict separation property, we proved the uniqueness of the weak solution. 
We also proposed a well posed and gradient stable continuous
finite element approximation of the model {for $d\leq 3$}, which takes the form of a discrete variational inequality. The proposed discrete scheme is computationally efficient in the calculation of the non-local contributions and it preserves the physical properties of the continuos solution. Finally, we showed numerical results in
2D for different test cases which verify the gradient stability properties of the proposed scheme and which show qualitatively how the topology of stationary states depends on the initial configuration according to the standard phase-ordering dynamics. Comparing the numerical results between the non-local and the local models, we observed that in the non-local case the long-range interactions make the interface regions wider and increase the degree of interaction between separated domains of pure phases. Also, comparing between test cases with different degrees of degeneracy of the phase mobility, we observed that the discrete solutions preserve the predicted strict separation property of the continuous solution up to numerical errors.

Future developments of the present work will concern the analysis of the full coupled non local model for tumor growth proposed in the paper, together with its application to real 3D patient-specific geometries. {In the latter context, we will compare the predictive ability to forecast tumor progression between the local and the non local versions of the model}.

\section*{Acknowledgments}
\noindent
{
This research has been performed in the framework of the MIUR-PRIN Grant 2020F3NCPX 
``Mathematics for industry 4.0 (Math4I4)''. The present paper also benefits from the support of 
the GNAMPA (Gruppo Nazionale per l'Analisi Matematica, la Probabilit\`a e le loro Applicazioni)
of INdAM (Istituto Nazionale di Alta Matematica).
{LS gratefully acknowledges support 
from the project CUP\_E55F22000270001 from GNAMPA
and the MUR grant  ``Dipartimento di
eccellenza 2023--2027'' for Politecnico di Milano.}
}

%
\bibliographystyle{plain}
\bibliography{biblio_NLCH} 
\end{document}